\documentclass[12pt,twoside]{article}
\usepackage[hmargin=0.8in,vmargin=0.9in]{geometry}
\usepackage[section]{placeins}
\usepackage{xcolor}
\usepackage{multirow}
\usepackage{fancyhdr}
\usepackage{graphicx}
\usepackage{subfigure}
\usepackage{amssymb}
\usepackage{amsmath}
\usepackage{amsfonts}
\usepackage{theorem}
\usepackage{mathrsfs}
\usepackage{mathtools}
\usepackage{bm}
\usepackage{color}
\usepackage{setspace}
\usepackage{exscale}
\usepackage{relsize}
\usepackage{epstopdf}
\usepackage{float}
\usepackage{cite}
\usepackage{cleveref}
\usepackage{color}
\usepackage{mwe}
\usepackage{subfig}
\usepackage{nicefrac}
\usepackage[normalem]{ulem}

\usepackage{booktabs,multirow} 
\usepackage{array} 
\usepackage{paralist} 
\usepackage{verbatim} 
\usepackage{subfigure} 

\theoremstyle{plain}			
\newtheorem{thm}{Theorem}[section]

\newtheorem{rmk}[thm]{Remark}
{\theorembodyfont{\rmfamily}}

\allowdisplaybreaks[1]

\numberwithin{equation}{section}
\numberwithin{figure}{section}
\numberwithin{table}{section}

\newcommand\eref[1]{(\ref{#1})}

\newcommand*\xbar[1]{%
  \hbox{%
    \vbox{%
      \hrule height 0.5pt 
      \kern0.4ex
      \hbox{%
        \kern-0.05em
        \ensuremath{#1}%
        \kern-0.00em
      }%
    }%
  }%
}

\newcommand{\mF}{\bm{F}}

\newcommand{\mU}{\bm{U}}

\newcommand{\mo}{\bm{0}}

\newcommand{\dx}{\Delta x}

\newcommand{\hf}{{\frac{1}{2}}}

\newcommand{\jph}{{j+\frac{1}{2}}}
\newcommand{\jmh}{{j-\frac{1}{2}}}

\usepackage{tabularx}
\newcolumntype{L}[1]{>{\raggedright\arraybackslash}p{#1}}
\newcolumntype{C}[1]{>{\centering\arraybackslash}p{#1}}
\newcolumntype{R}[1]{>{\raggedleft\arraybackslash}p{#1}}

\title{Central-Upwind Scheme for the Phase-Transition Traffic Flow Model}
\author{Shaoshuai Chu\thanks{Department of Mathematics, RWTH Aachen University, 52056 Aachen, Germany; Department of Mathematics and
Shenzhen International Center for Mathematics, Southern University of Science and Technology, Shenzhen, 518055, China;
{\tt chu@igpm.rwth-aachen.de}}, Alexander Kurganov\thanks{Department of Mathematics and Shenzhen International Center for Mathematics,
Southern University of Science and Technology, Shenzhen, 518055, China; {\tt alexander@sustech.edu.cn}}, Saeed
Mohammadian\thanks{School of Civil Engineering, the University of Queensland, Brisbane Qld, 4072, Australia;
{\tt s.mohammadian@uq.edu.au}}, and Zuduo Zheng\thanks{School of Civil Engineering, the University of Queensland, Brisbane Qld, 4072,
Australia; {\tt zuduo.zheng@uq.edu.au}} }
\date{}
\begin{document}

\maketitle
\begin{abstract}
Phase-transition models are an important family of non-equilibrium continuum traffic flow models, offering properties like replicating
complex traffic phenomena, maintaining anisotropy, and promising potentials for accommodating automated vehicles. However, their complex
mathematical characteristics such as discontinuous solution domains, pose numerical challenges and limit their exploration in traffic flow
theory. This paper focuses on developing a robust and accurate numerical method for phase-transition traffic flow models: We propose a
second-order semi-discrete central-upwind scheme specifically designed for phase-transition models. This novel scheme
incorporates the projection onto appropriate flow domains, ensuring enhanced handling of discontinuities and maintaining physical
consistency and accuracy. We demonstrate the efficacy of the proposed scheme through extensive and challenging numerical tests, showcasing
their potential to facilitate further research and application in phase-transition traffic flow modeling. The ability of phase-transition
models to embed the ``time-gap''---a crucial element in automated traffic control---as a conserved variable aligns seamlessly with the
control logic of automated vehicles, presenting significant potential for future applications, and the proposed numerical scheme now
substantially facilitates exploring such potentials.
\end{abstract}

\noindent
{\bf Key words:} Phase-transition traffic flow model; free traffic flow; congested traffic flow; finite-volume methods; central-upwind
schemes.

\noindent
{\bf AMS subject classification:} 76M12, 65M08, 76A30, 76T99, 35L65.

\section{Introduction}
Continuum models treat traffic flow as a compressible fluid and analyze its dynamics using aggregated state variables like flow and density,
proving beneficial in real-world traffic management and control \cite{mohammadian2021benchmarking}. Over time, a variety of continuum models
have been developed, each incorporating different empirical and behavioral characteristics of traffic flow; see, e.g., the
recent works in \cite{mohammadian2023continuum,ZQYZX2022,ZHANG2023128556,QZWZW2017,ZQLZX2023}. These models are generally classified into
two main categories: equilibrium and non-equilibrium models. Equilibrium models are based on flow conservation principles, consistently
linking speed and density without time or space differentiation. In contrast, non-equilibrium models use separate dynamic partial
differential equations for speed and density, with their universal relationship only applicable in steady-state equilibrium, irrespective of
time and space.

Non-equilibrium continuum traffic flow models are crucial for analyzing traffic dynamics, including the development and spread of persistent
waves from initial disturbances to near-equilibrium conditions. To be well-defined, these models must meet the following three criteria:
\begin{itemize}
\item Adherence to causality laws, implying that drivers respond to stimuli ahead of them, and thus, characteristic traffic waves, generated
by the model, must not exceed average traffic speeds due to this directional response \cite{daganzo1995requiem};
\item Considering the fact that the negative correlation between driver speed and intervehicular spacing restricts macroscopic traffic
conditions between specific paths on the flow-density diagram, the model must ensure that traffic states stay within these viable
trajectories and automatically satisfy the consistency conditions such as zero speed at maximum density \cite{treiber2010three};
\item Distinct representation of the dynamics in free-flow and congested phases, ensuring that phenomena unique to each phase are accurately
simulated and exclusive to that phase \cite{kerner2016failure}.
\end{itemize}

To the best of our knowledge, phase-transition models of traffic flow, initially proposed in \cite{Colombo02} and further developed in
\cite{blandin2011general,CG14}, represent the only class of non-equilibrium models that endogenously meet all of the aforementioned three
criteria. The models in this family consider distinct hyperbolic equations for the free-flow and congestion phases, respectively. In the
free-flow phase, a simplified equilibrium model with traffic speed being a direct function of density is used, while for the congestion
phase, a $2\times2$ hyperbolic system is employed, linking the density and speed through a conserved variable, which is equivalent to the
inverse of the average time gap.

The phase-transition models are promising from a physical perspective in several regards. An empirical study comparing empirical vehicle
trajectories has demonstrated the efficacy of phase-transition models in replicating various complex traffic phenomena, including traffic
hysteresis and the convective propagation of congested states \cite{blandin2013phase}. Furthermore, phase-transition traffic flow models
introduced in \cite{Colombo02a} and \cite{blandin2011general} are excellent foundational models for extension to accommodate the governing
principles of Connected and Automated Vehicles (CAVs). This suitability stems from the fact that the time gap, an essential decision
variable used in the control laws of CAVs, is already directly incorporated as a conserved variable in these models.

Nevertheless, the significant potential of phase transition models has been overlooked in the traffic flow theory literature, especially
concerning numerical simulations and analysis. This research gap can be attributed to the inherent mathematical complexity of these models,
which makes them very challenging to simulate numerically using the existing numerical schemes. We can only refer the reader to
\cite{chalons2008godunov}, where a Godunov-type method for a phase-transition model was proposed. In this method, mesh cells along phase
boundaries are adjusted for accurate projections and a Glimm-type sampling technique post-projection is employed to preserve the numerical
solution's integrity. We stress that the Godunov scheme and other Godunov-type upwind schemes rely on (approximately) solving the
(generalized) Riemann problems arising at each cell interface. This introduces significant complexity, primarily due to extensive
conditional logic for determining the nature of the solution---shock waves, rarefaction waves, or contact discontinuities---based on local
conditions. This detailed case-by-case approach significantly increases computational complexity and overhead, rendering Godunov-type upwind
schemes less practical for large-scale applications.

This paper aims to bridge the aforementioned research gap by developing a numerical scheme that overcomes the aforementioned difficulties.
The contributions of this study are two-fold. First, we aim to develop a novel, accurate, and robust numerical scheme and demonstrate its
performance in solving the Riemann problems across a comprehensive range of challenging scenarios. Second, we utilize the proposed scheme to
shed light on the phase-transition model's behavior in a wide range of physically tangible real-world traffic scenarios.

The proposed numerical scheme is based on further development of the second-order semi-discrete central-upwind (CU) schemes, which were
introduced in \cite{Kurganov01,Kurganov00,Kurganov07} (see also \cite{CCHKL_22,CKX_24,KX_22} for recent low-dissipation modifications of the CU schemes) as a ``black-box'' solver for general hyperbolic systems of conservation laws. In this
paper, we incorporate novel mechanisms to tailor the CU scheme for phase-transition models, ensuring physical consistency and
accuracy by projecting computed results onto the correct flow domains. To demonstrate the performance of the proposed scheme, extensive and
challenging numerical tests are conducted. Furthermore, after sufficiently demonstrating the scheme's performance for Riemann problems, we
implement the phase-transition traffic flow model in several complex scenarios, aiming to elaborate on the physical interpretation of the
proposed model.

The rest of the paper is organized as follows. In \S\ref{sec2}, we describe the phase-transition traffic flow model. In \S\ref{sec3}, we
introduce the designed numerical method for the studied model. Finally, in \S\ref{sec4}, we apply the developed scheme to a number of
numerical examples. We demonstrate that the proposed scheme can capture the solution sharply and in a non-oscillatory manner.

\section{Phase-Transition Traffic Flow Model}\label{sec2}
The phase-transition model, which was initially proposed in \cite{Colombo02} (see also \cite{Colombo02a,CG14}), couples two different
traffic flow models. In the free-flow domain, where the traffic density $\rho$ is smaller than a given threshold $\rho^f_{\rm cr}$, the
following scalar equation for $\rho$ is used:
\begin{equation}
\rho_t+(\rho V_f(\rho))_x=0,
\label{2.1}
\end{equation}
where $V_f$ is the traffic speed in free flow. In the congested-flow domain, where $\rho>\rho^f_{\rm cr}$, the following second-order
traffic flow model is utilized:
\begin{equation}
\left\{\begin{aligned}
&\rho_t+(\rho V_c(\rho,q))_x=0,\\
&q_t+((q-q^*)V_c(\rho,q))_x=0,
\end{aligned}\right.
\label{2.2}
\end{equation}
where $V_c$ is the traffic speed in congested regions. In \eref{2.1} and \eref{2.2}, $x$ is the spatial variable, $t$ is time, $q$ is a
flow-type variable equal to the inverse of drivers’ average time gap, and $q^*$ is a constant parameter (we define the corresponding time
gap derived from $q^*$ as the equilibrium time gap). Note that $q$ should not be confused with the traffic flow rate, which is obtained as
$\rho V_f$ in free flow and $\rho V_c$ in congestion. The functions \( V_f \) and \( V_c \) are given by
\begin{equation}
V_f(\rho)=V_{\max},\quad V_c(\rho,q)=\Big(1-\frac{\rho}{\rho_{\max}}\Big)\frac{q}{\rho}
\label{2.3}
\end{equation}
with $V_{\max}$ and $\rho_{\max}$ being the maximum speed and density.

We notice that $q$ is neither defined nor used in \eref{2.1}, but in order to treat the model transition regions, we need to introduce $q$
in the free-flow domain. This is done by substituting $V_c(\rho,q)=V_{\max}$ into \eref{2.3}, which results in
\begin{equation}
q(\rho)=\frac{V_{\max}}{\frac{1}{\rho}-\frac{1}{\rho_{\max}}}.
\label{2.4}
\end{equation}

\begin{figure}[ht!]
\centerline{\includegraphics[trim=0.2cm 0.1cm 0.2cm 0.3cm, clip, width=8cm]{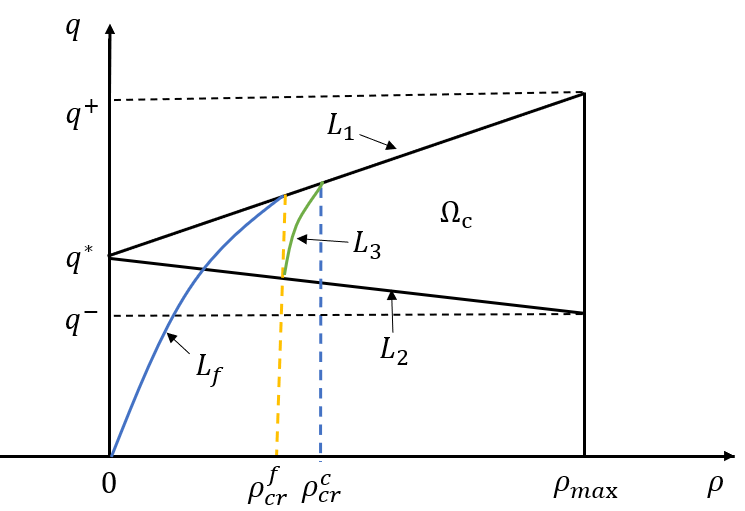}}
\caption{\sf$(\rho,q)$-diagram of the phase-transition traffic flow model.\label{fig1}}
\end{figure}

In Figure \ref{fig1}, we illustrate the phase transition traffic flow model in the $(\rho,q)$-plane. The figure divides the plane into two
distinct domains: the \emph{free-flow curve} $L_f$ and the \emph{congested-flow domain} $\Omega_c$:

\medskip
$\bullet$ The free-flow curve $L_f$ is valid for $\rho\in[0,\rho_{\rm cr}^f]$ and represents traffic conditions, in which the vehicles
travel at maximum speed without interaction;

$\bullet$ The congested-flow domain $\Omega_c$ is bounded by the vertical line $\rho=\rho_{\max}$ and the curves 
\begin{equation}
L_1:~q=q^*+\frac{q^+-q^*}{\rho_{\max}}\,\rho,\quad L_2:~q=q^*+\frac{q^--q^*}{\rho_{\max}}\,\rho,\quad\mbox{and}\quad 
L_3:~q=\frac{\rho\rho_{\max}}{\rho_{\max}-\rho}\,V_{c+},
\label{2.5f}
\end{equation}
which represent different constraints on the traffic flux under high-density conditions. In \eref{2.5f}, $q^\pm$ are the maximum/minimum
admissible values of $q$ and $V_{c+}$ is the maximum speed in congested domain. Note that the line $L_1$ and the curve $L_3$ intersect at
the point with 
\begin{equation*}
\rho=\rho^c_{\rm cr}:=\frac{2\rho_{\max}q^*}{\rho_{\max}V_{c+}+2q^*-q^++\sqrt{(\rho_{\max}V_{c+}+2q^*-q^+)^2+4(q^+-q^*)q^*}}.
\end{equation*}

We stress that the phase transition occurs only along the boundaries between the free and congested domains, specifically near the critical
point $\rho=\rho_{\rm cr}^f$, which marks the threshold at which the flow shifts between free and congested states depending on the
surrounding traffic conditions. As it was shown in \cite{Colombo02a}, once a solution enters either the free or congested region, it remains
there unless it reaches a transition boundary. Therefore, Figure \ref{fig1} not only presents the structural partitioning of the
$(\rho,q)$-plane but also illustrates the allowed directions and constraints of phase transitions in traffic dynamics.

The studied system \eref{2.1}--\eref{2.3} can be put into the following vector form:
\begin{equation}
\mU_t+\mF(\mU)_x=\mo,
\label{2.5}
\end{equation}
where
\begin{equation}
\mU=\begin{pmatrix}\rho\\q\end{pmatrix}\quad{\rm and}\quad F(\mU)=\begin{cases}
(\rho V_{\max},q V_{\max})^\top&\mbox{if}~(\rho,q)\in L_f,\\
(\rho V_c,(q-q^*)V_c)^\top&\mbox{if}~(\rho,q)\in\Omega_c.
\end{cases}
\label{2.6}
\end{equation}
We stress that for the free-flow domain, the equation for $q$ will not be numerically solved, but the vector form \eref{2.5}--\eref{2.6}
will be convenient for the presentation of the numerical method in the next section.

\section{Numerical Method}\label{sec3}
In this section, we develop a second-order semi-discrete CU scheme for the studied phase-transition traffic flow model.

We cover the computational domain with the uniform cells $C_j:=[x_\jmh,x_\jph]$ of size $\dx$ centered at $x_j=(x_\jmh+x_\jph)/2$ and denote
by $\,\xbar\mU_j(t)$ the computed cell averages
\begin{equation}
\xbar\mU_j(t):\approx\frac{1}{\dx}\int\limits_{C_j}\mU(x,t)\,{\rm d}x.
\end{equation}
We suppose that at a certain time $t\ge0$, the cell averages $\,\xbar\mU_j(t)$ are available and they all belong to either $L_f$ or
$\Omega_c$; see Figure \ref{fig1}. Note that $\,\xbar\mU_j$ and many other indexed quantities, which will be introduced below, are
time-dependent and from here on we will suppress this dependence for the sake of brevity.

In the semi-discrete framework, the numerical solution of \eref{2.5}--\eref{2.6} is evolved in time by solving the following system of
ordinary differential equations (ODEs):
\begin{equation}
\frac{{\rm d}\xbar\mU_j}{{\rm d}t}=-\frac{\bm{{\cal F}}_\jph-\bm{{\cal F}}_\jmh}{\dx},
\label{3.2}
\end{equation}
where the CU numerical flux $\bm{{\cal F}}_\jph$ is given by \cite{Kurganov07}
\begin{equation}
\bm{{\cal F}}_\jph=\frac{a^+_\jph\bm F(\mU^-_\jph)-a^-_\jph\bm F(\mU^+_\jph)}{a^+_\jph-a^-_\jph}+
\frac{a^+_\jph a^-_\jph}{a^+_\jph-a^-_\jph}\left(\mU^+_\jph-\mU^-_\jph-\bm Q_\jph\right).
\label{3.3}
\end{equation}
Here, $\mU^\pm_\jph$ are the right- and left-sided point values of $\bm U$ at the cell interface $x=x_\jph$ obtained with the help of a
conservative, second-order accurate, and non-oscillatory piecewise linear reconstruction.

\begin{figure}[ht!]
\centerline{\includegraphics[trim=0.cm 0.3cm 0.cm 0.4cm, clip, width=12cm]{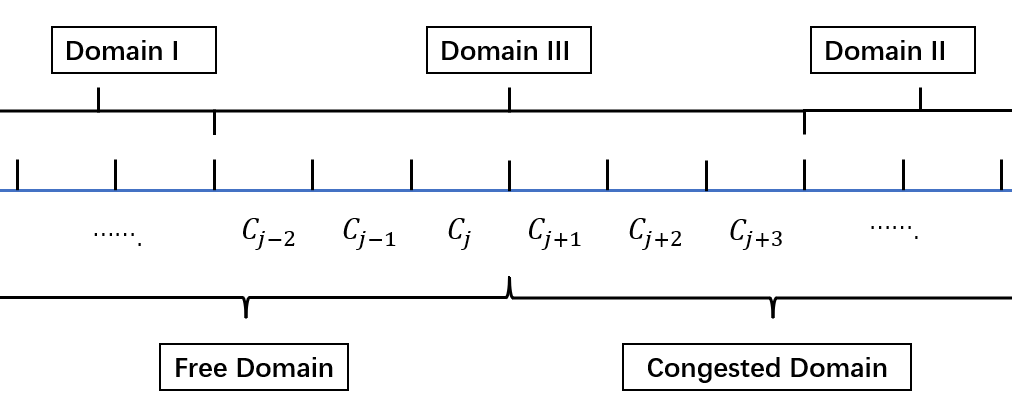}}
\caption{\sf Domains I, II, and III of the phase-transition traffic flow model.\label{figaaa}}
\end{figure}

We use different reconstructions in the free domain away from the phase-transition areas (Domain I), in the congested domain away from the
phase-transition areas (Domain II), and in the phase-transition areas consisting of the six cells (three on the left and three on the right)
around each interface between the free- and congested-flow domains (Domain III); see Figure \ref{figaaa}.

In Domain I, we only need to compute the one-sided point values $\rho^\pm_\jph$. To this end, we apply the generalized minmod reconstruction
(see Appendix \ref{appa}) to the density $\rho$ with $\theta=1.5$ in \eref{equ3.5}. In Domain III, we reconstruct both $\rho^\pm_\jph$ and
$q^\pm_\jph$ using the same minmod limiter from Appendix \ref{appa} but with $\theta=1$ in \eref{equ3.5} now applied to both $\rho$ and $q$.
Note that here we use a smaller value of $\theta$ to minimize possible oscillations in the interface areas. In Domain II, applying the
minmod limiter to the $\rho$- and $q$-component of the computed solution may lead to relatively large oscillations. In order to suppress
them, we apply the minmod limiter with $\theta=1.5$ in \eref{equ3.5} to the local characteristic variables, which are obtained using the
local characteristic decomposition (see, e.g., \cite{don9,Nonomura20,Qiu02,Shu20} and references therein); see Appendix \ref{appb} for
details.
\begin{rmk}
We detect the interface using the following simple approach: if $(\,\xbar\rho_J-\rho^f_{\rm cr})(\,\xbar\rho_{J+1}-\rho^f_{\rm cr})\le0$,
then $x=x_{J+\hf}$ is the interface, and cells $C_{J-2},\ldots,C_{J+3}$ belong to Domain III.
\end{rmk}

After reconstructing the point values $(\rho^\pm_\jph,q^\pm_\jph)$, we modify them unless $(\rho^\pm_\jph,q^\pm_\jph)\in L_f$ or
$(\rho^\pm_\jph,q^\pm_\jph)\in\Omega_c$. The details on the proposed modification (projection) onto $L_f/\Omega_c$ can be founded in
\S\ref{appc}.

In \eref{3.3}, $a^\pm_\jph$ are the one-sided local speeds of propagation, which can be estimated using the largest ($\lambda_2$) and the
smallest ($\lambda_1$) eigenvalues of the Jacobian $A=\frac{\partial\mF}{\partial\mU}$. This can be done as follows (see, e.g.,
\cite{Kurganov01,Kurganov07}):
\begin{equation}
\begin{aligned}
a^+_\jph=\max\left\{\lambda_2\big(\mU^-_\jph\big),\lambda_2\big(\mU^+_\jph\big),0\right\},\quad
a^-_\jph=\min\left\{\lambda_1\big(\mU^-_\jph\big),\lambda_1\big(\mU^+_\jph\big),0\right\}.
\end{aligned}
\label{3.3a}
\end{equation}
Computing $a^\pm_\jph$ according to \eref{3.3a} requires, however, considering four possible cases depending on whether the left and right
reconstructed states $\mU^\pm_\jph$ are in free or congested domain as in the free-flow domain the eigenvalues of
\begin{equation}
\begin{aligned}
A=\begin{pmatrix}
V_{\rm max} & 0\\
0 & V_{\rm max}
\end{pmatrix}
\end{aligned}
\end{equation}
are
$\lambda^f_1=\lambda^f_2=V_{\max}$ and in the congested-flow domain, the eigenvalues of 
\begin{equation}
\begin{aligned}
A=\begin{pmatrix}
-\dfrac{q}{\rho_{\max}}&\dfrac{\rho_{\max}-\rho}{\rho_{\max}}\\[2.5ex]
\dfrac{q(q^*-q)}{\rho^{\,2}}&
\dfrac{(q^*-2 q)(\rho-\rho_{\max})}{\rho\rho_{\max}}
\end{pmatrix}
\end{aligned}
\end{equation}
are 
\begin{equation}
\lambda^c_1=(q-q^*)\Big(\frac{1}{\rho}-\frac{2}{\rho_{\max}}\Big)-\frac{q^*}{\rho_{\max}},\quad\lambda^c_2=V_c(\rho,q).
\end{equation}

We then obtain:

$\bullet$ If $(\rho^-_\jph,q^-_\jph)\in\Omega_c$ and $(\rho^+_\jph,q^+_\jph)\in\Omega_c$, we take
\begin{equation}
a^+_\jph=\max\left\{\lambda^c_2\big(\mU^-_\jph\big),\lambda^c_2\big(\mU^+_\jph\big),0\right\},\quad
a^-_\jph=\min\left\{\lambda^c_1\big(\mU^-_\jph\big),\lambda^c_1\big(\mU^+_\jph\big),0\right\};
\end{equation}

$\bullet$ If both $(\rho^-_\jph,q^-_\jph)\in L_f$ and $(\rho^+_\jph,q^+_\jph)\in L_f$, we take
\begin{equation}
a^+_\jph=V_{\max},\quad a^-_\jph=0;
\end{equation}

$\bullet$ If $(\rho^-_\jph,q^-_\jph)\in\Omega_c$ and $(\rho^+_\jph,q^+_\jph)\in L_f$, we take
\begin{equation}
\begin{aligned}
&a^+_\jph=\max\left\{\lambda^c_2\big(\mU^-_\jph\big),\lambda^f_2\big(\mU^+_\jph\big),0\right\}=
\max\left\{\lambda^c_2\big(\mU^-_\jph\big),V_{\max},0\right\}=V_{\max},\\
&a^-_\jph=\min\left\{\lambda^c_1\big(\mU^-_\jph\big),\lambda^f_1\big(\mU^+_\jph\big),0\right\}=
\min\Big\{\lambda^c_1\big(\mU^-_\jph\big),V_{\max},0\Big\}=\min\Big\{\lambda^c_1\big(\mU^-_\jph\big),0\Big\};
\end{aligned}
\label{3.4a}
\end{equation}

$\bullet$ If $(\rho^-_\jph,q^-_\jph)\in L_f$ and $(\rho^+_\jph,q^+_\jph)\in\Omega_c$, we take
\begin{equation}
\begin{aligned}
&a^+_\jph=\max\left\{\lambda^f_2\big(\mU^-_\jph\big),\lambda^c_2\big(\mU^+_\jph\big),0\right\}=
\max\left\{V_{\max},\lambda^c_2\big(\mU^+_\jph\big),0\right\}=V_{\max},\\
&a^-_\jph=\min\left\{\lambda^f_1\big(\mU^-_\jph\big),\lambda^c_1\big(\mU^+_\jph\big),0\right\}=
\min\Big\{V_{\max},\lambda^c_1\big(\mU^+_\jph\big),0\Big\}=\min\left\{\lambda^c_1\big(\mU^+_\jph\big),0\right\}.
\end{aligned}
\label{3.5a}
\end{equation}
Notice that in the computation of $a^\pm_\jph$ in \eref{3.4a} and \eref{3.5a}, we have used the fact that
$\lambda^c_1\big(\mU^\pm_\jph\big)\le\lambda^c_2\big(\mU^\pm_\jph\big)\le V_{c+}\le V_{\max}$.

Finally, the term $\bm Q_\jph$ in \eref{3.3} represents a ``built-in'' anti-diffusion and is given by \cite{Kurganov07}
\begin{equation}
\bm Q_\jph={\rm minmod}\left(\mU^+_\jph-\mU^*_\jph,\mU^*_\jph-\mU^-_\jph\right),
\end{equation}
where
\begin{equation}
\mU^*_\jph=\frac{a^+_\jph\mU^+_\jph-a^-_\jph\mU^-_\jph-\left\{\bm F\big(\mU^+_\jph\big)-\bm F\big(\mU^-_\jph\big)\right\}}
{a^+_\jph-a^-_\jph}.
\end{equation}

The ODE system \eref{3.2}--\eref{3.3} has to be numerically integrated using an appropriate ODE solver. Upon completion of every time step
(or every stage of a Runge-Kutta method), we project the obtained set of $(\,\xbar\rho_j,\,\xbar q_j)$ onto either $L_f$ or $\Omega_c$
according to the procedure introduced in \S\ref{appc}.

\subsection{Projection onto $L_f$ or $\Omega_c$ }\label{appc}
Assume that we have obtained a set of data (either reconstructed point values or evolved cell averages) for $\rho$ and $q$. If a certain
pair $(\widetilde\rho,\widetilde q)$ is neither on $L_f$ or in $\Omega_c$, we have to project it onto either $L_f$ or $\Omega_c$ according
to the following four  possible cases:

{\bf Case I:} If $\widetilde\rho\le\rho^f_{\rm cr}$, which means that $\widetilde\rho$ corresponds to the free flow, we replace
$\widetilde q$ with $q(\widetilde\rho)$ computed using \eref{2.4};

{\bf Case II:} If $\rho^f_{\rm cr}<\widetilde\rho<\rho^c_{\rm cr}$ and
$\widetilde q>\dfrac{\widetilde\rho\rho_{\max}}{\rho_{\max}-{\widetilde\rho}}\,V_{c+}$, that is, the point $(\widetilde\rho,\widetilde q)$
is above the line $L_3$ (see Figure \ref{fig1}), we shift this point vertically down to $L_3$ and replace $\widetilde q$ with
$\dfrac{\widetilde\rho\rho_{\max}}{\rho_{\max}-{\widetilde\rho}}\,V_{c+}$;

{\bf Case III:} If $\rho^c_{\rm cr}<\widetilde\rho$ and $\widetilde q>q^*+\dfrac{q^+-q^*}{\rho_{\max}}\,\widetilde\rho$, that is, the point
$(\widetilde\rho,\widetilde q)$ is above the line $L_1$ (see Figure \ref{fig1}), we shift this point vertically down to $L_1$ and replace
$\widetilde q$ with $q^*+\dfrac{q^+-q^*}{\rho_{\max}}\,\widetilde\rho$;

{\bf Case IV:} If $\rho^f_{\rm cr}<\widetilde\rho$ and $\widetilde q<q^*+\dfrac{q^--q^*}{\rho_{\max}}\,\widetilde\rho$, that is, the point
$(\widetilde\rho,\widetilde q)$ is below the line $L_2$ (see Figure \ref{fig1}), we shift this point vertically up to $L_2$ and replace
$\widetilde q$ with $q^*+\dfrac{q^--q^*}{\rho_{\max}}\,\widetilde\rho$.

Note that in all of these four cases, we modify $q$ only without changing $\rho$ so that vehicles are not getting artificially added to or
removed from the road.

\section{Numerical Examples}\label{sec4}
In this section, we evaluate the performance of the proposed numerical scheme using extensive numerical tests. To this end, we use the
three-stage third-order strong stability preserving (SSP) Runge-Kutta method (see, e.g., \cite{Gottlieb12,Gottlieb11}) in order to
numerically solve the ODE system \eref{3.2}--\eref{3.3} using the CFL number 0.4, which is selected based on our initial investigation
findings.

For clarity and better organization, we categorize our investigation scenarios into two main types based on the similarities in initial and
boundary conditions. The first category includes Riemann problems with simple free boundary conditions (Examples 1 and 2), while the second 
category consists of more realistic, traffic-oriented cases (Examples 3 and 4), where complex Dirichlet boundary conditions are imposed at
the downstream boundary to simulate traffic situations that induce backward-propagating waves. In all of our examples, we employ the
parameters listed in Table \ref{tab0}, which include both the phase-transition model settings and the physicality constraints for the
congested domain. These parameters are chosen to reflect empirical traffic conditions: the free-flow speed is set to 30$\,$m/s, consistent
with typical freeway speed limits; the maximum density ranges from 110 to 170$\,$veh/km based on average vehicle size and spacing; and the
critical density is defined as 15--30\% of the maximum density. The time gap and the bounds $q^+$ and $q^-$ are also selected from typical
empirical ranges, as discussed in \cite{mohammadian2021benchmarking,treiber2010three}.
\begin{table}[ht!]
\centering
\caption{\sf Parameters of the phase-transition model adopted for numerical tests}
\begin{tabular}{cccccccccc}
\toprule
\textbf{Parameter}&$V_{\max}$&$V_{c+}$&$\rho_{\max}$&$q^*$&$\rho^f_{\rm cr}$&$q^+$&$q^-$\\
\hline
\textbf{Value}&30&24&0.16&0.6&0.02&0.93186&0.18856\\
\bottomrule
\end{tabular}\label{tab0}
\end{table}

\paragraph*{Example 1.} In the first example, we consider several Riemann initial data of the form
\begin{equation}
\begin{aligned}
\mU(x,0)=\begin{cases}\mU_L&\mbox{if}~x<x_0,\\\mU_R&\mbox{otherwise},\end{cases}\quad x_0=40000,
\end{aligned}
\end{equation}
prescribed in the computational domain $[0,800000]$ subject to the free boundary conditions. This setting corresponds to a Riemann problem
on a long road with open boundaries.
\begin{table}[ht!]
\centering
\caption{\sf Initial conditions for Tests 1--12}
\begin{tabular}{ccccccccc}
\toprule
\multirow{2}{*}{\textbf{Test}}&\multicolumn{4}{c}{$\mathbf{Left-side}$}&\multicolumn{4}{c}{$\mathbf{Right-side}$}\\
\cline{2-9}
& \textbf{Phase}&\textbf{$\rho_L$}&\textbf{$V_L$}&\multicolumn{1}{c}{$q_L-q^*$}&\textbf{Phase}&\textbf{$\rho_R$}&\textbf{$V_R$}&$q_R -q^*$\\
\hline\hline
1&Free&0.011 &30&-0.2456&Congested&0.0825&4.5113&0.1684\\
2&Free&0.011 &30&-0.2456&Congested&0.0775&4.5945&0.0906\\
3&Free&0.0075&30&-0.3639&Congested&0.0675&5.338 &0.02325\\
4&Free&0.001 &30&-0.5698&Congested&0.0625&4.73  &-0.1149\\
5&Free&0.001 &30&-0.5698&Congested&0.0875&2.9945&-0.02175\\
\hline
6&Congested&0.128 &0.42321&-0.3291&Congested&0.0375&13.838 &0.07778\\
7&Congested&0.0375&13.838 &0.07778&Congested&0.128 &0.42321&-0.3291\\
\hline
8 &Congested&0.0825&4.5113&0.1684  &Free&0.011 &30&-0.2456\\
9 &Congested&0.0775&4.5945&0.0906  &Free&0.011 &30&-0.2456\\
10&Congested&0.0675&5.338 &0.02324 &Free&0.0075&30&-0.3639\\
11&Congested&0.0625&4.73  &-0.1149 &Free&0.001 &30&-0.5698\\
12&Congested&0.0875&2.9945&-0.02175&Free&0.001 &30&-0.5698\\
\bottomrule
\end{tabular}
\label{tab1}
\end{table}
Table \ref{tab1} presents twelve comprehensive test cases used in our first example, along with the essential information about traffic
states on either side of the initial discontinuity. All of these test cases are numerically challenging and can give rise to complex waves,
comprising of multiple wave types such as shock waves, contact discontinuities, and rarefaction waves due to involving transitions between
the free-flow and congested phases.

In Tests 1--5, free-flow traffic upstream of the initial discontinuity encounters congested traffic downstream. In all of these test cases,
vehicles in the upstream will have to decelerate in order to adapt to the traffic conditions downstream, and the solution would involve
shock waves. While these test cases are quite comparable in terms of traffic speed on either side of the initial discontinuity, there are
significant differences in the sign and magnitude of $q-q^*$ on the right-hand side (RHS), which result in different driving behaviors and
are expected to manifest in the corresponding traffic waves. In general, the farther the quantity $q-q^*$ is from zero, the greater the
deviation of the drivers' {\em average time gap}, $1/q$, from the equilibrium time gap, $1/q^*$ (the time gap corresponding to the
equilibrium speed-density relationship in \eref{2.3} with $q=q^*$).  

Therefore, in scenarios where $q-q^*>0$ on the RHS (Tests 1--3), drivers' average time gap on this side of the discontinuity is already
smaller than the equilibrium time gap, which implies more aggressive driving behavior and a tendency toward denser traffic. In these cases,
an intermediate state is expected to form, connecting the upstream free-flow to the congested downstream, as $q$ on the RHS is larger than
equilibrium. In contrast, in scenarios where $q-q^*<0$ on the RHS (Tests 4 and 5), drivers' average time gap is larger than the equilibrium
time gap, and an intermediate state is not expected to arise, as the $q$ on the RHS of the discontinuity is already smaller than the
equilibrium.

Tests 6 and 7 represent scenarios in which there is no phase transition, as the traffic states on both sides of the initial discontinuity
are in a congested phase. From a computational and numerical standpoint, these test cases are less challenging. However, they are included
for their physical significance, namely, to examine the model's solution structure for Riemann problems where one state's time gap is close
to the equilibrium, while the other exhibits a deviation from this equilibrium.

In Tests 8--12, traffic state in the upstream of the initial discontinuity is congested, whereas free-flow traffic exists in the downstream.
As a result, vehicles in the upstream will go through acceleration speed adaptation manoeuvres in order to adapt to the traffic condition in
the downstream. These cases encompass a variety of transitions between traffic phases, characterized by the fluctuation of the quantity
$q-q^*$ across a spectrum of values, from negative to positive. This variation indicates diverse scenarios where drivers' average preferred
time gap deviates from the equilibrium time gap, and are considered in the numerical tests for two aspects: (a) to evaluate and demonstrate
the performance of the proposed scheme under such complex phase-transitions and (b) to shed light on the phase-transition model's solution
structures for such initial discontinuities and provide physical interpretations.

\begin{figure}[ht!]
\centerline{\includegraphics[trim=0.1cm 0.1cm 0.7cm 0.2cm, clip, width=5.4cm]{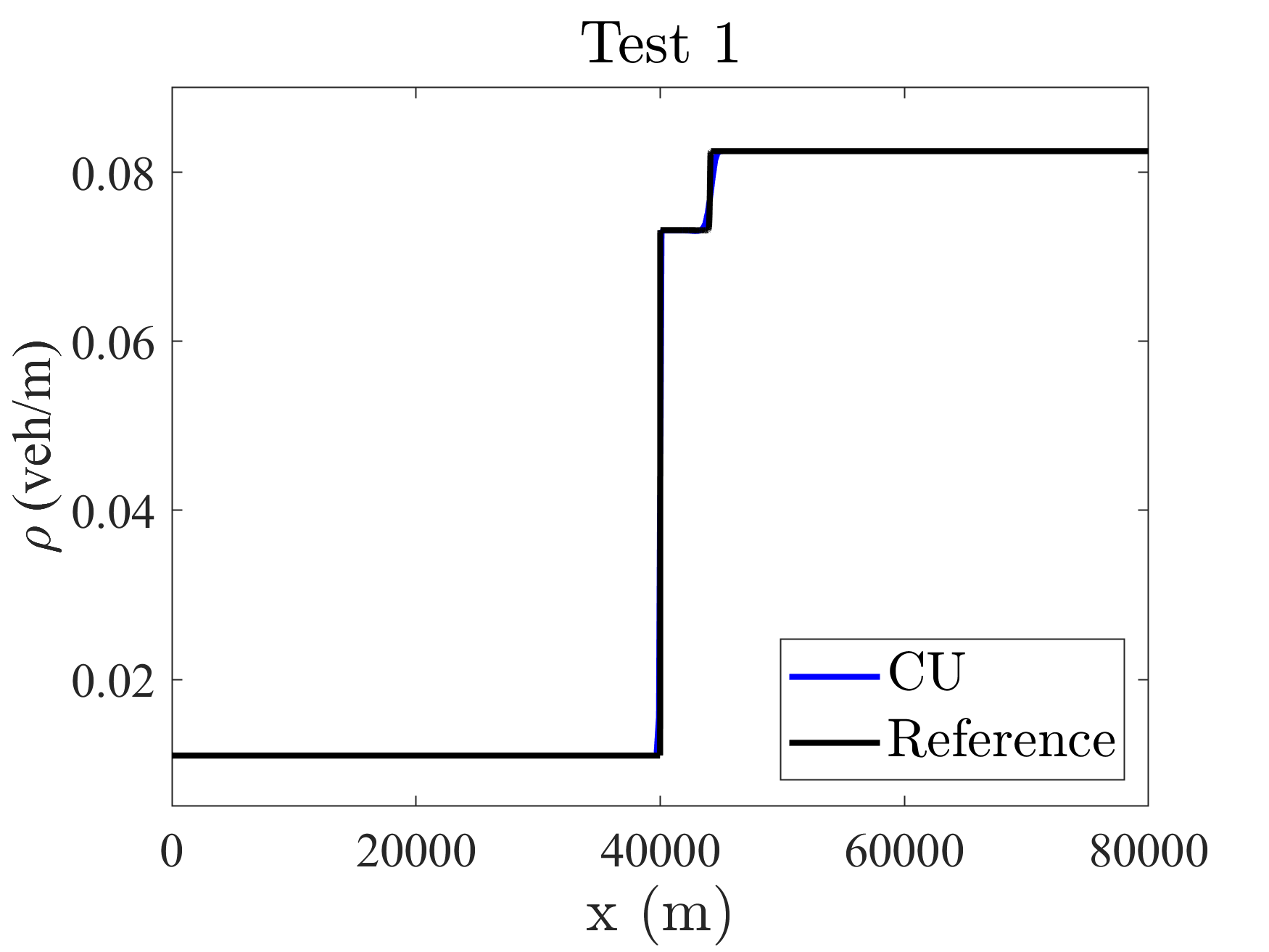}\hspace{1cm}
            \includegraphics[trim=0.1cm 0.1cm 0.7cm 0.2cm, clip, width=5.4cm]{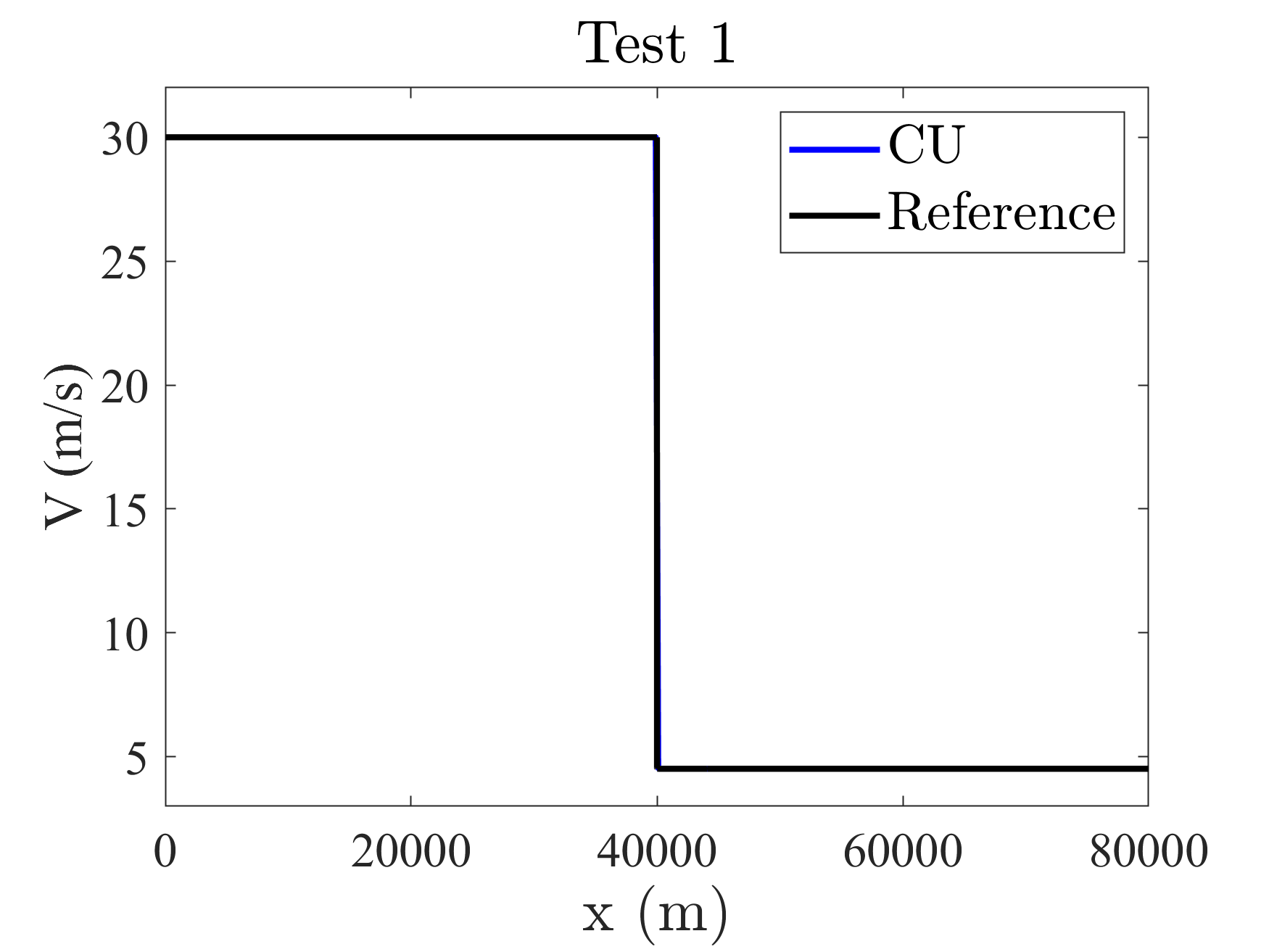}}
\vskip6pt
\centerline{\includegraphics[trim=0.1cm 0.1cm 0.7cm 0.2cm, clip, width=5.4cm]{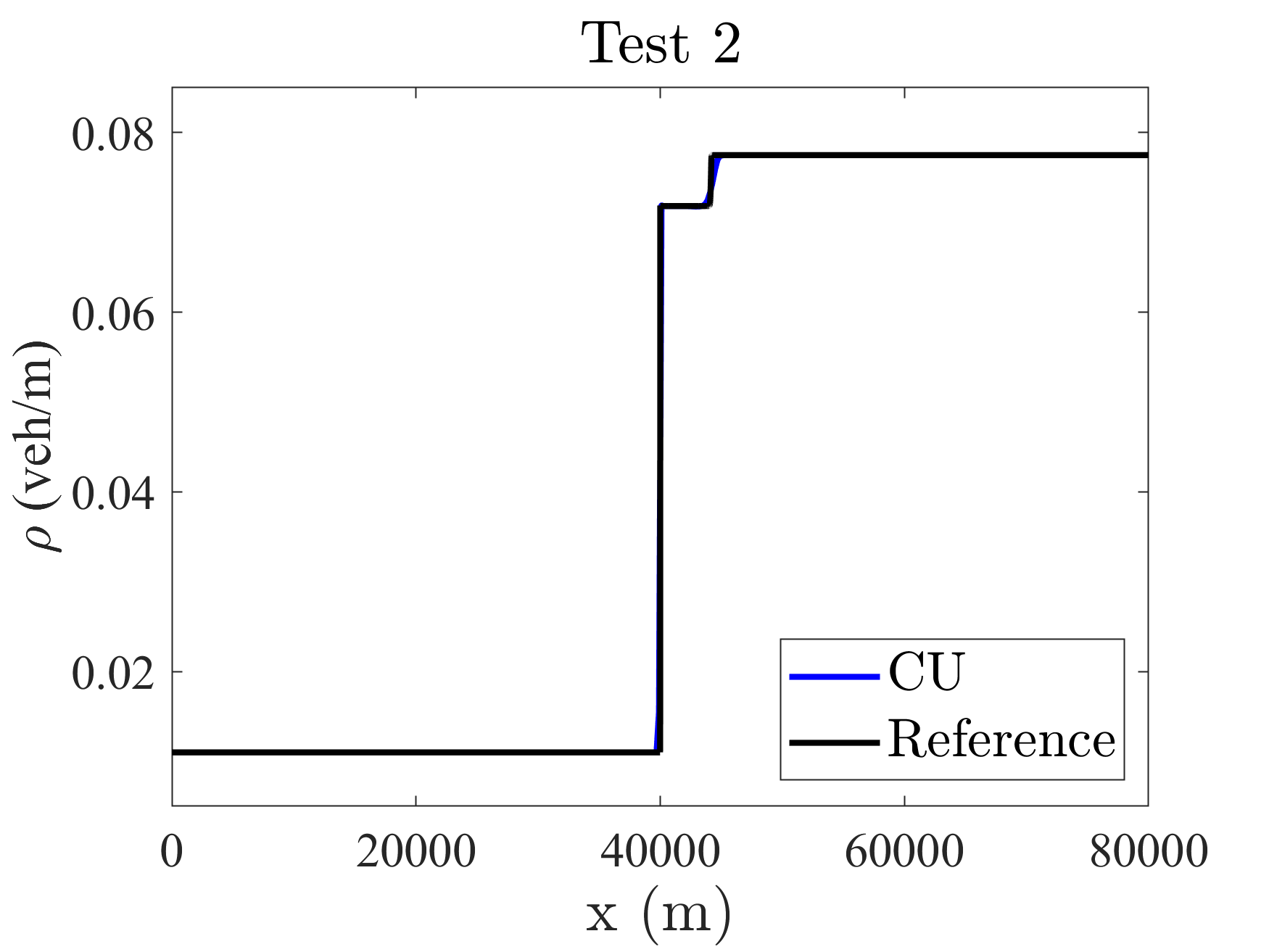}\hspace{1cm}
            \includegraphics[trim=0.1cm 0.1cm 0.7cm 0.2cm, clip, width=5.4cm]{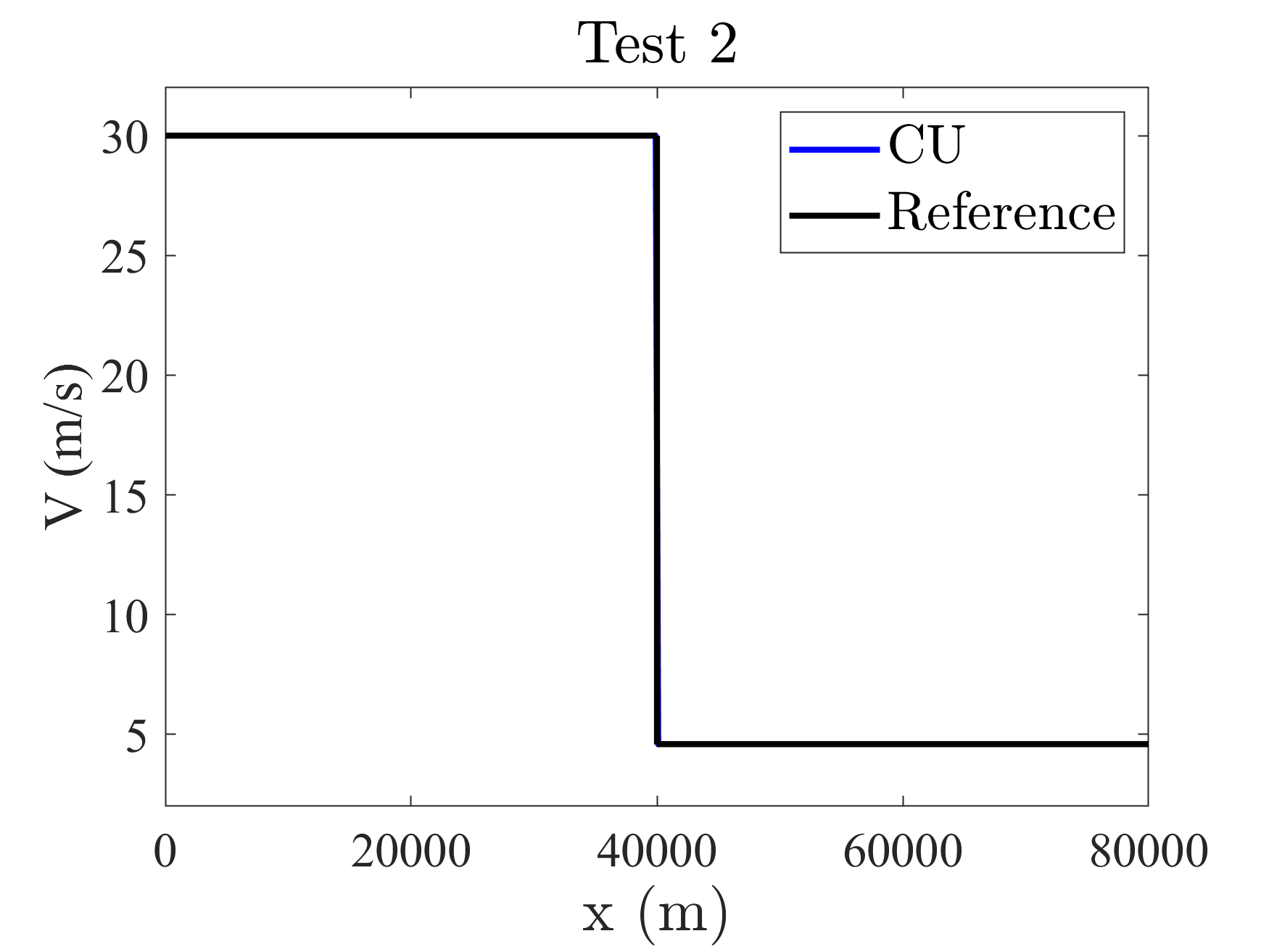}}
\vskip6pt
\centerline{\includegraphics[trim=0.1cm 0.1cm 0.7cm 0.2cm, clip, width=5.4cm]{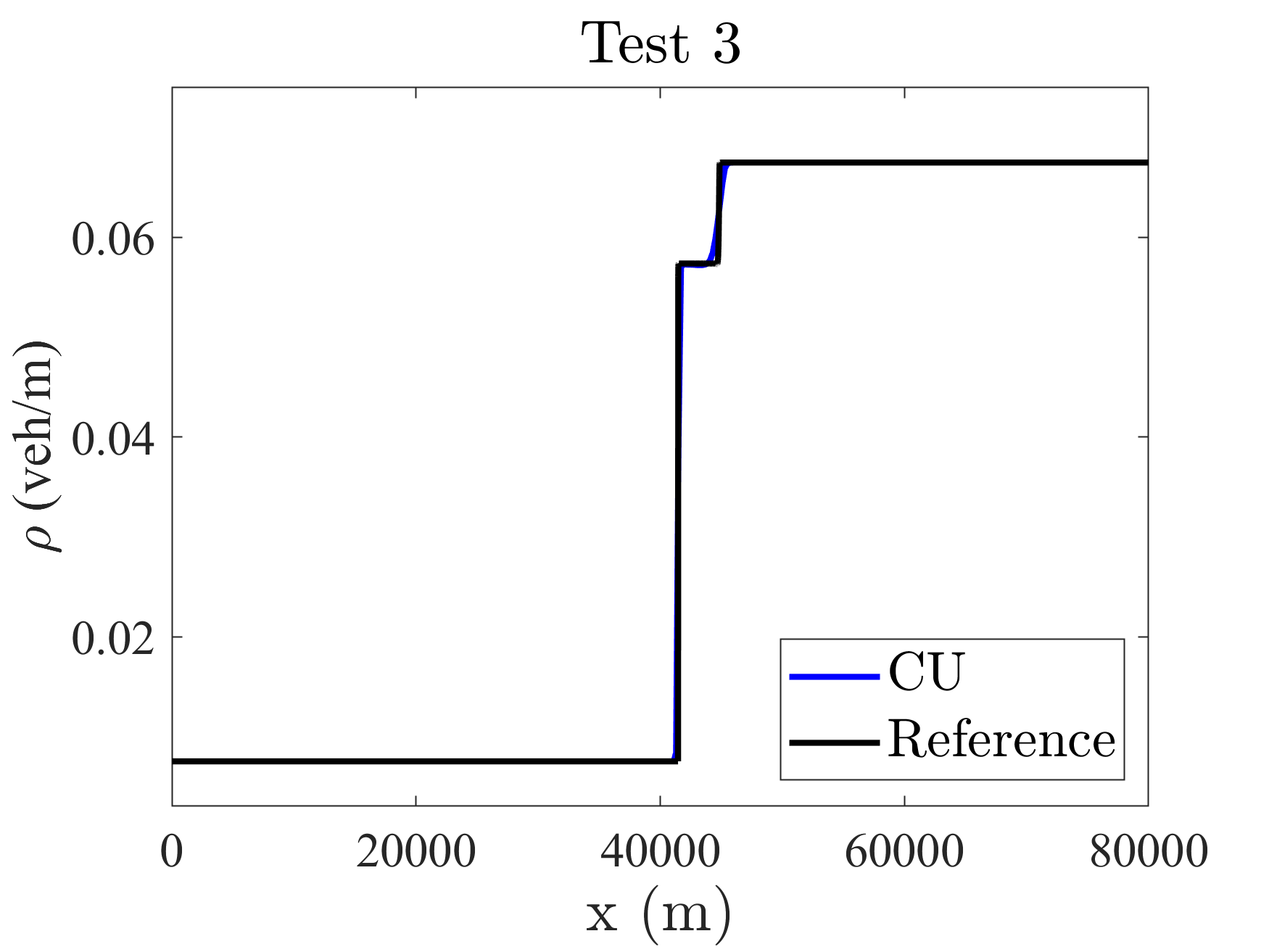}\hspace{1cm}
            \includegraphics[trim=0.1cm 0.1cm 0.7cm 0.2cm, clip, width=5.4cm]{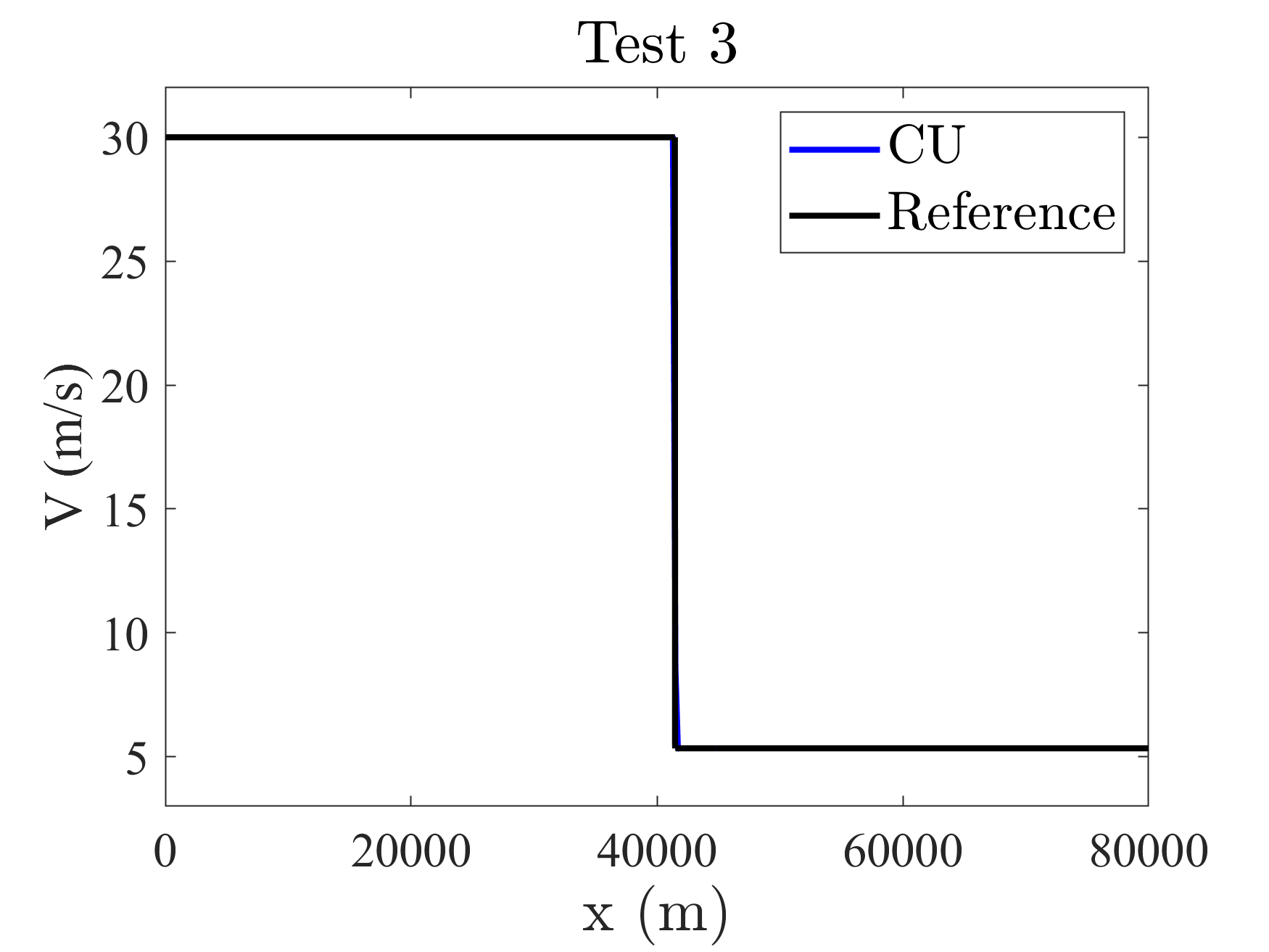}}
\vskip6pt
\centerline{\includegraphics[trim=0.1cm 0.1cm 0.7cm 0.2cm, clip, width=5.4cm]{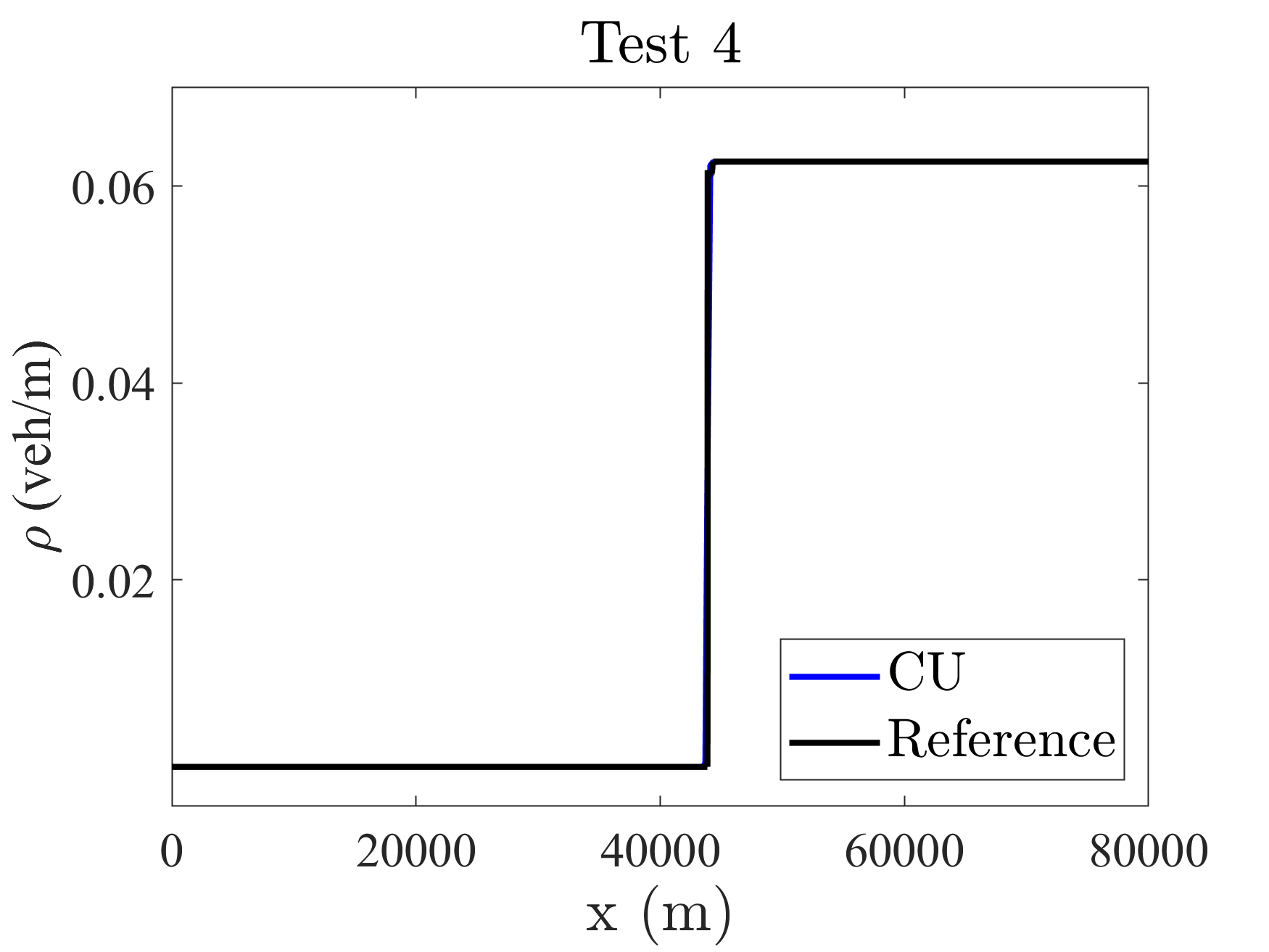}\hspace{1cm}
            \includegraphics[trim=0.1cm 0.1cm 0.7cm 0.2cm, clip, width=5.4cm]{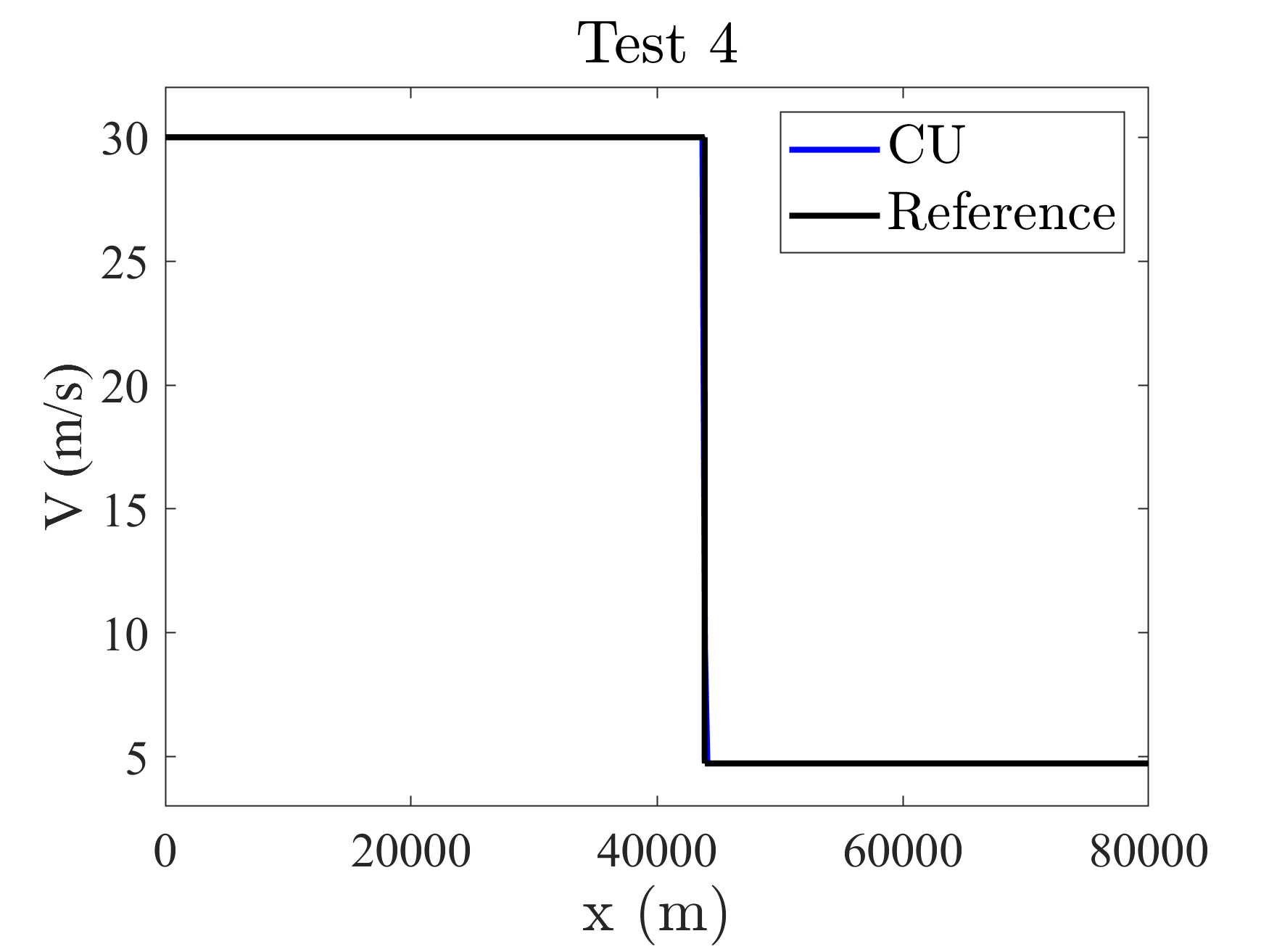}}
\caption{\sf Example 1, Tests 1--4: $\rho$ (left column) and $V$ (right column) at time $T_{\rm final}=900$.\label{fig41}}
\end{figure}

\begin{figure}[ht!]
\centerline{\includegraphics[trim=0.1cm 0.1cm 0.7cm 0.2cm, clip, width=5.4cm]{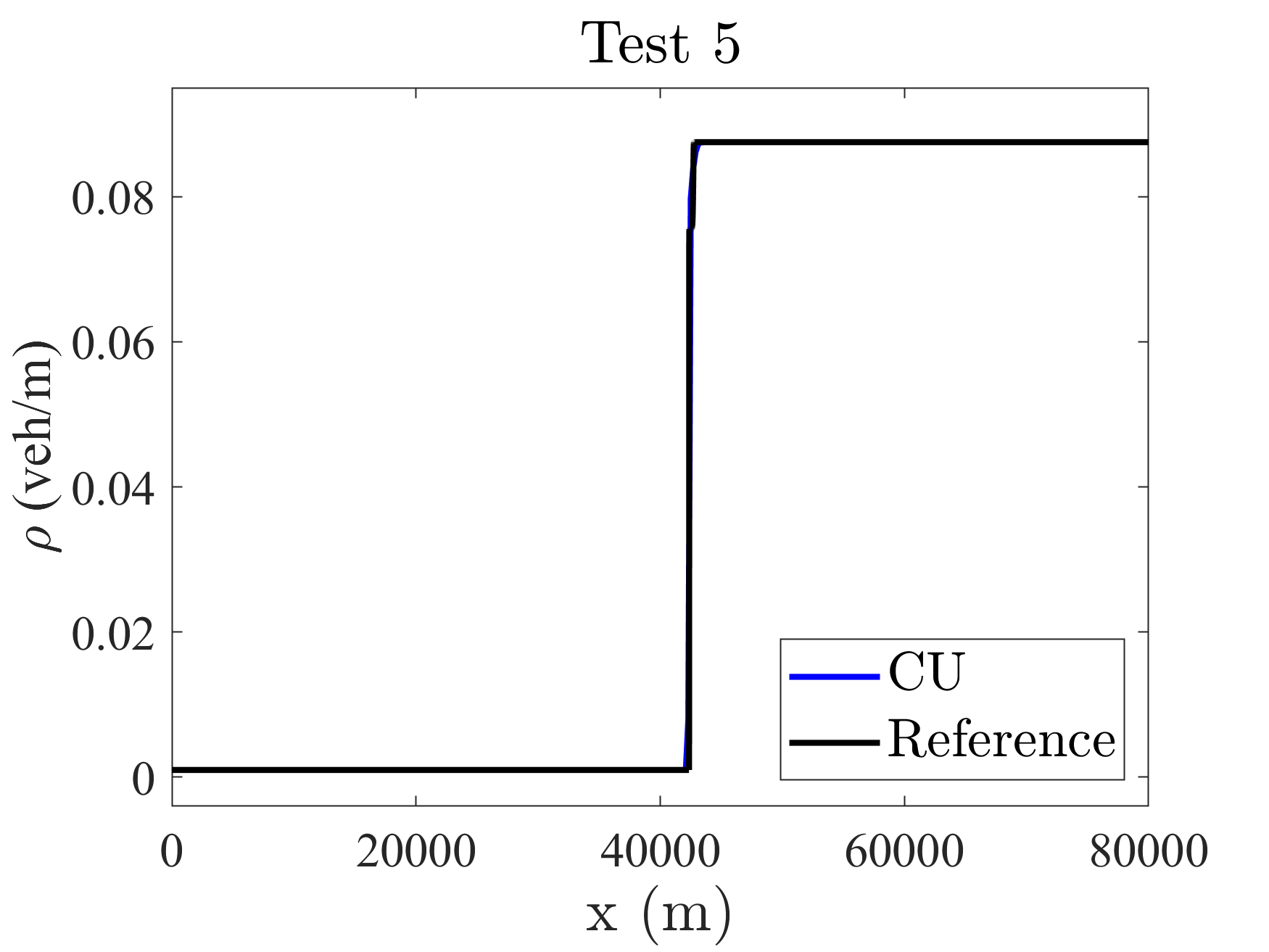}\hspace{1cm}
            \includegraphics[trim=0.1cm 0.1cm 0.7cm 0.2cm, clip, width=5.4cm]{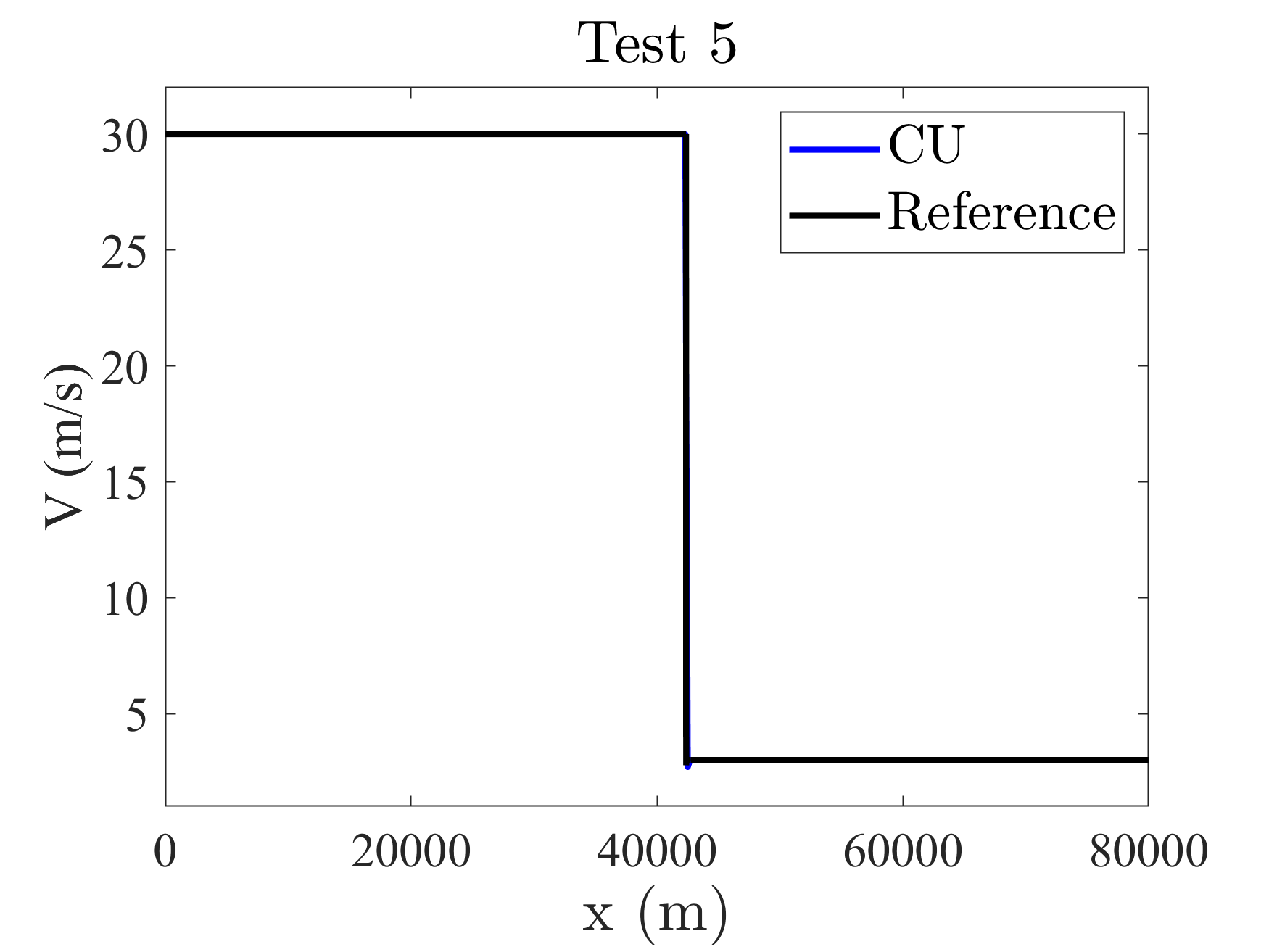}}
\vskip6pt
\centerline{\includegraphics[trim=0.1cm 0.1cm 0.7cm 0.2cm, clip, width=5.4cm]{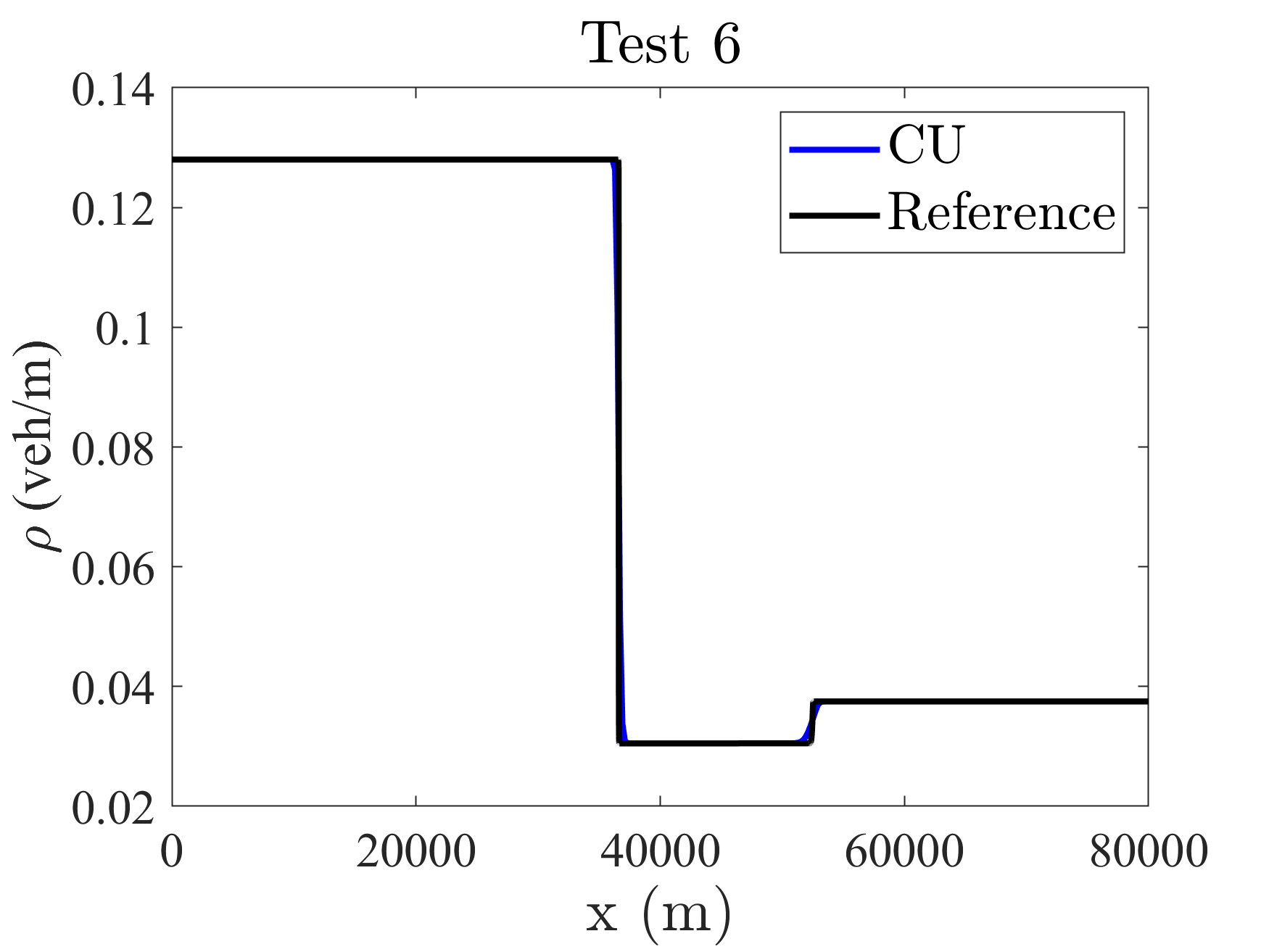}\hspace{1cm}
            \includegraphics[trim=0.1cm 0.1cm 0.7cm 0.2cm, clip, width=5.4cm]{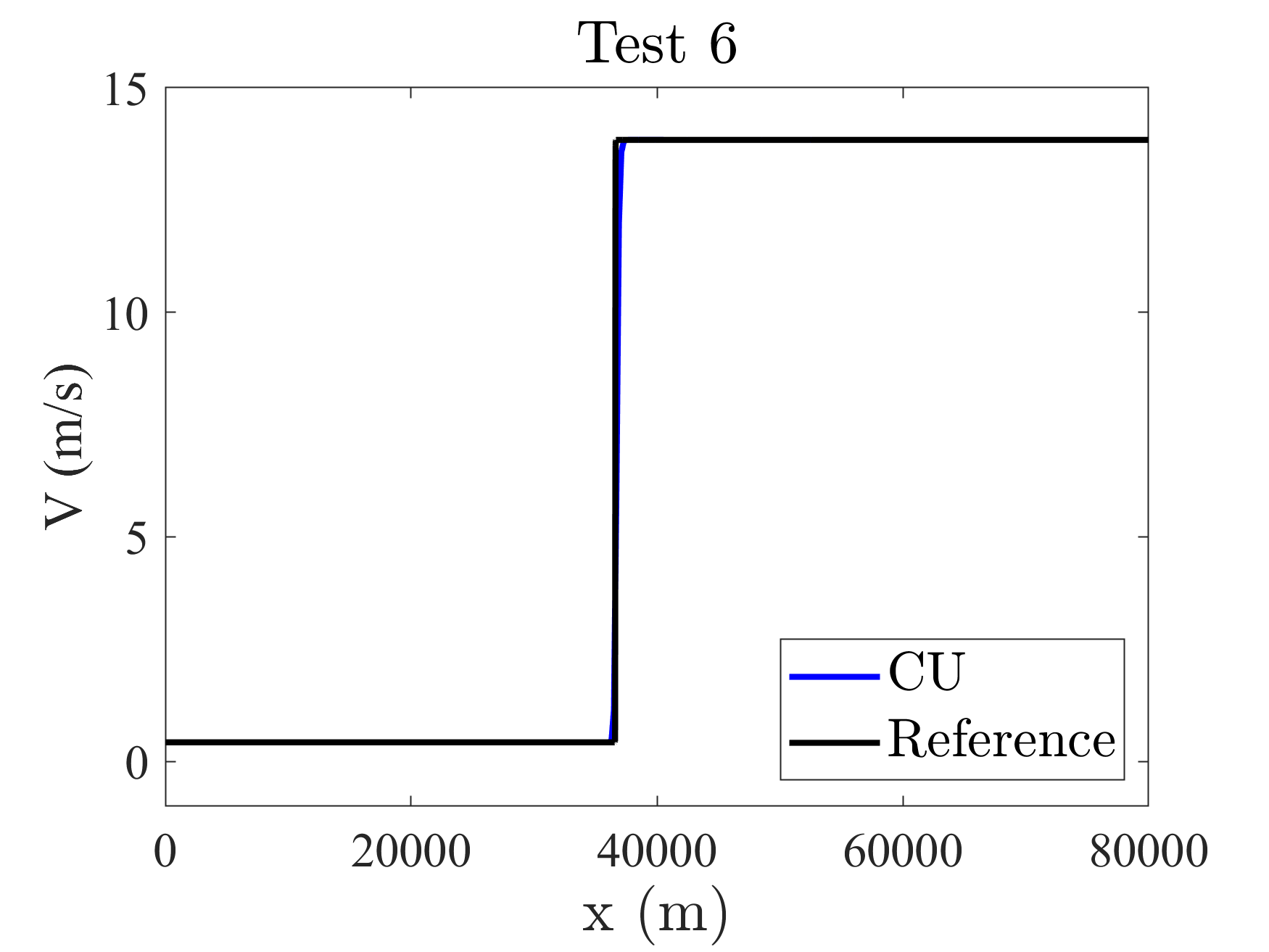}}
\vskip6pt
\centerline{\includegraphics[trim=0.1cm 0.1cm 0.7cm 0.2cm, clip, width=5.4cm]{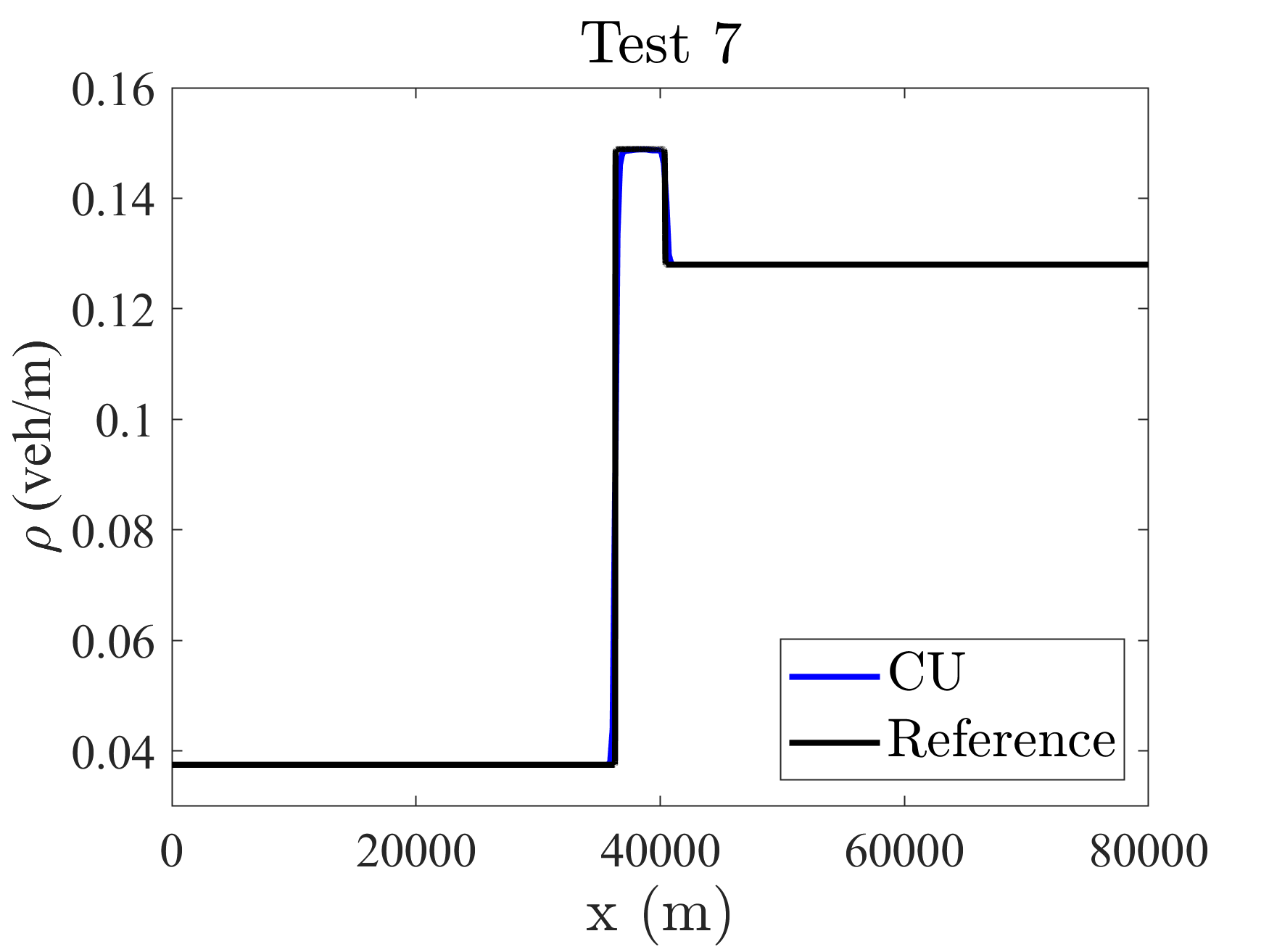}\hspace{1cm}
            \includegraphics[trim=0.1cm 0.1cm 0.7cm 0.2cm, clip, width=5.4cm]{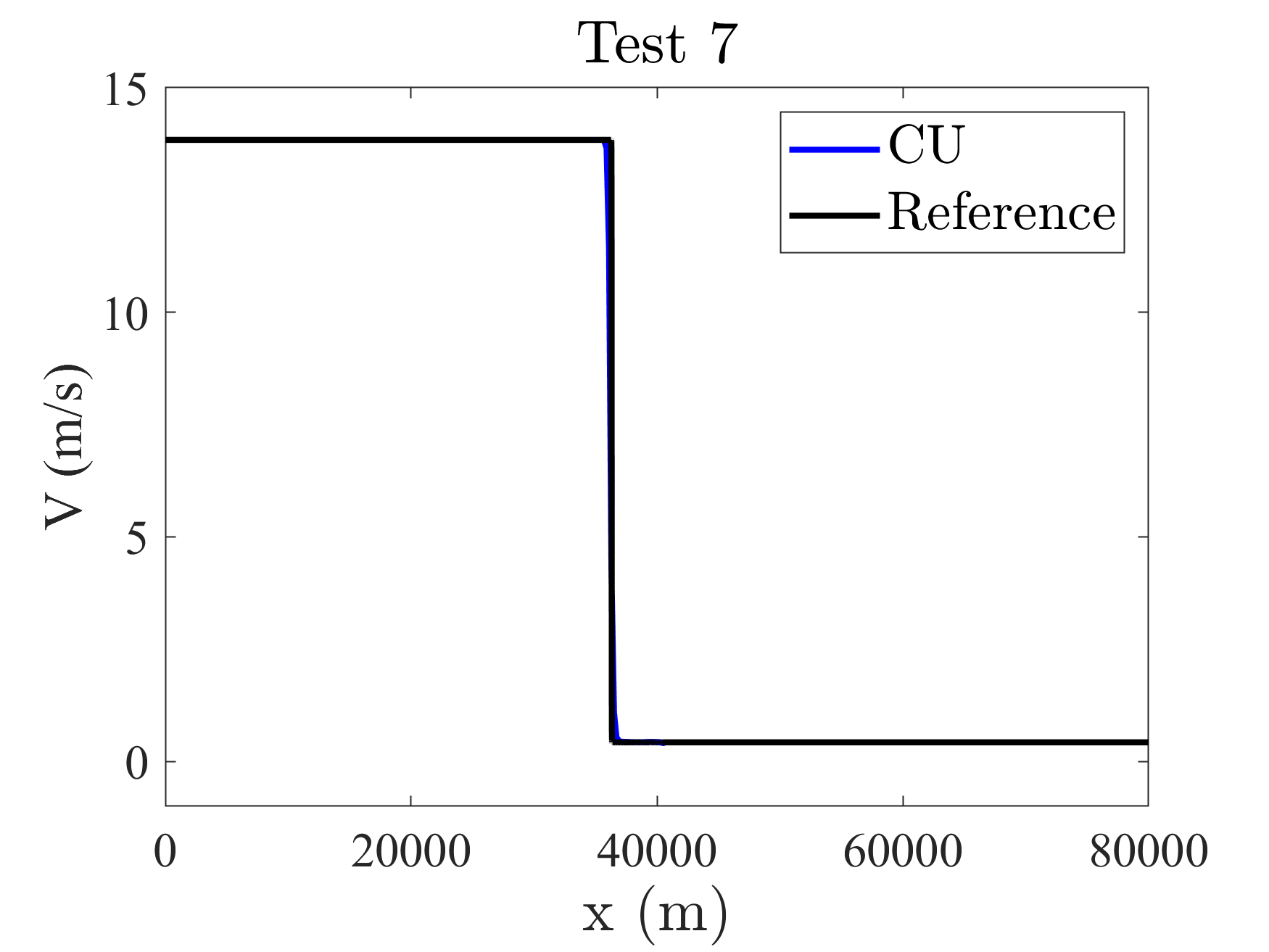}}
\vskip6pt
\centerline{\includegraphics[trim=0.1cm 0.1cm 0.7cm 0.2cm, clip, width=5.4cm]{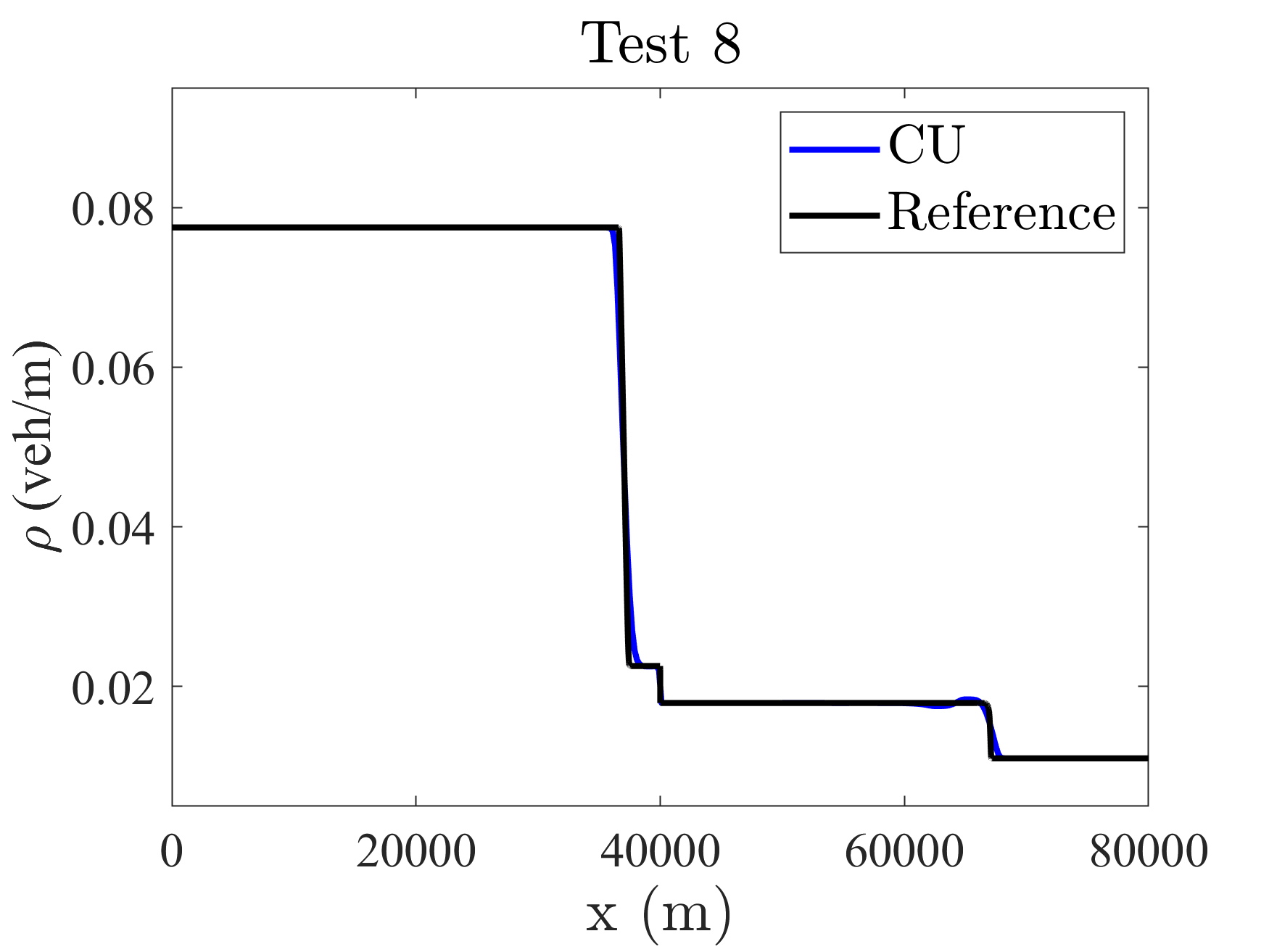}\hspace{1cm}
            \includegraphics[trim=0.1cm 0.1cm 0.7cm 0.2cm, clip, width=5.4cm]{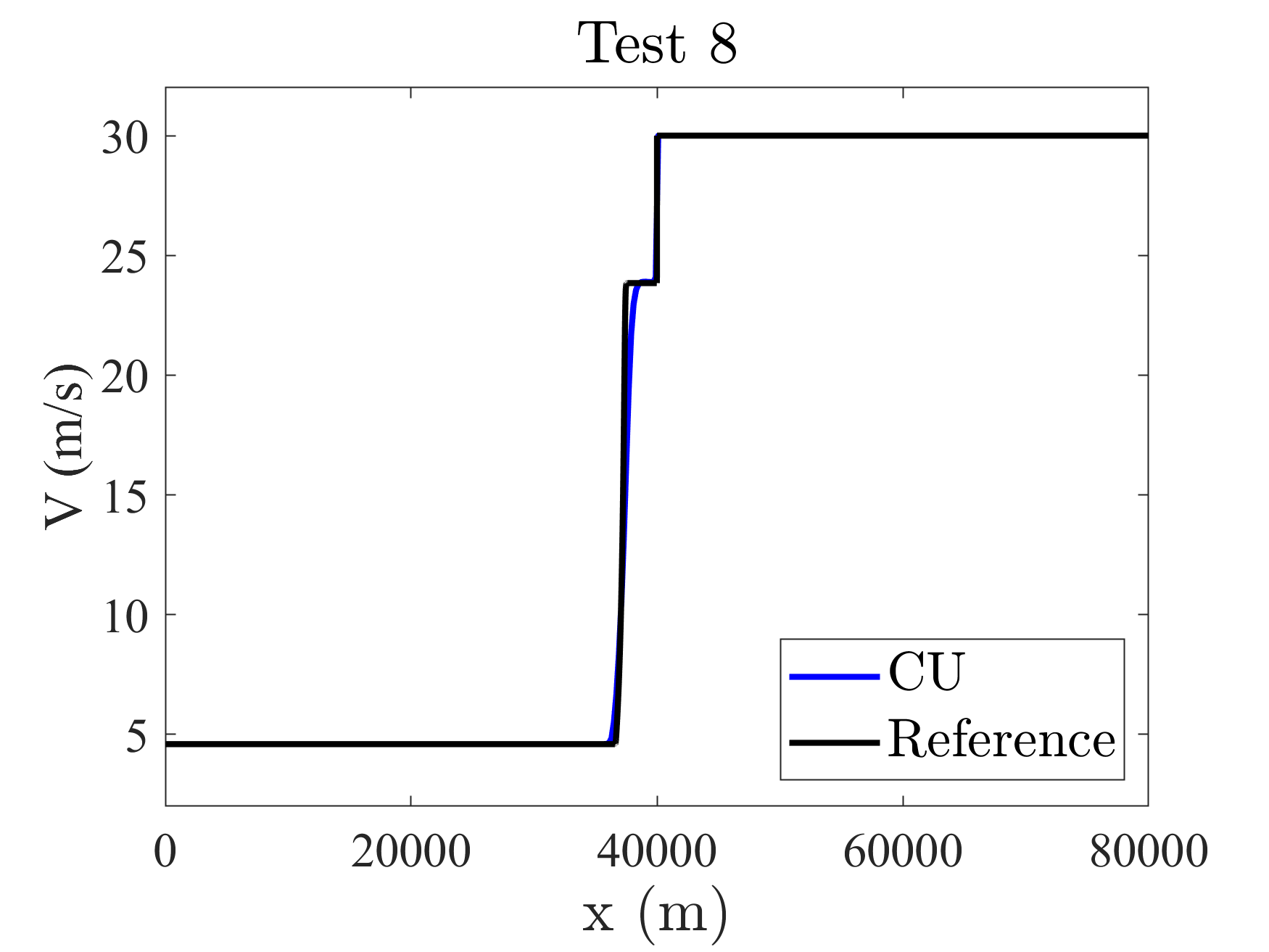}}
\caption{\sf Example 1, Tests 5--8: $\rho$ (left column) and $V$ (right column) at  time $T_{\rm final}=900$.\label{fig42}}
\end{figure}

\begin{figure}[ht!]
\centerline{\includegraphics[trim=0.1cm 0.1cm 0.7cm 0.2cm, clip, width=5.4cm]{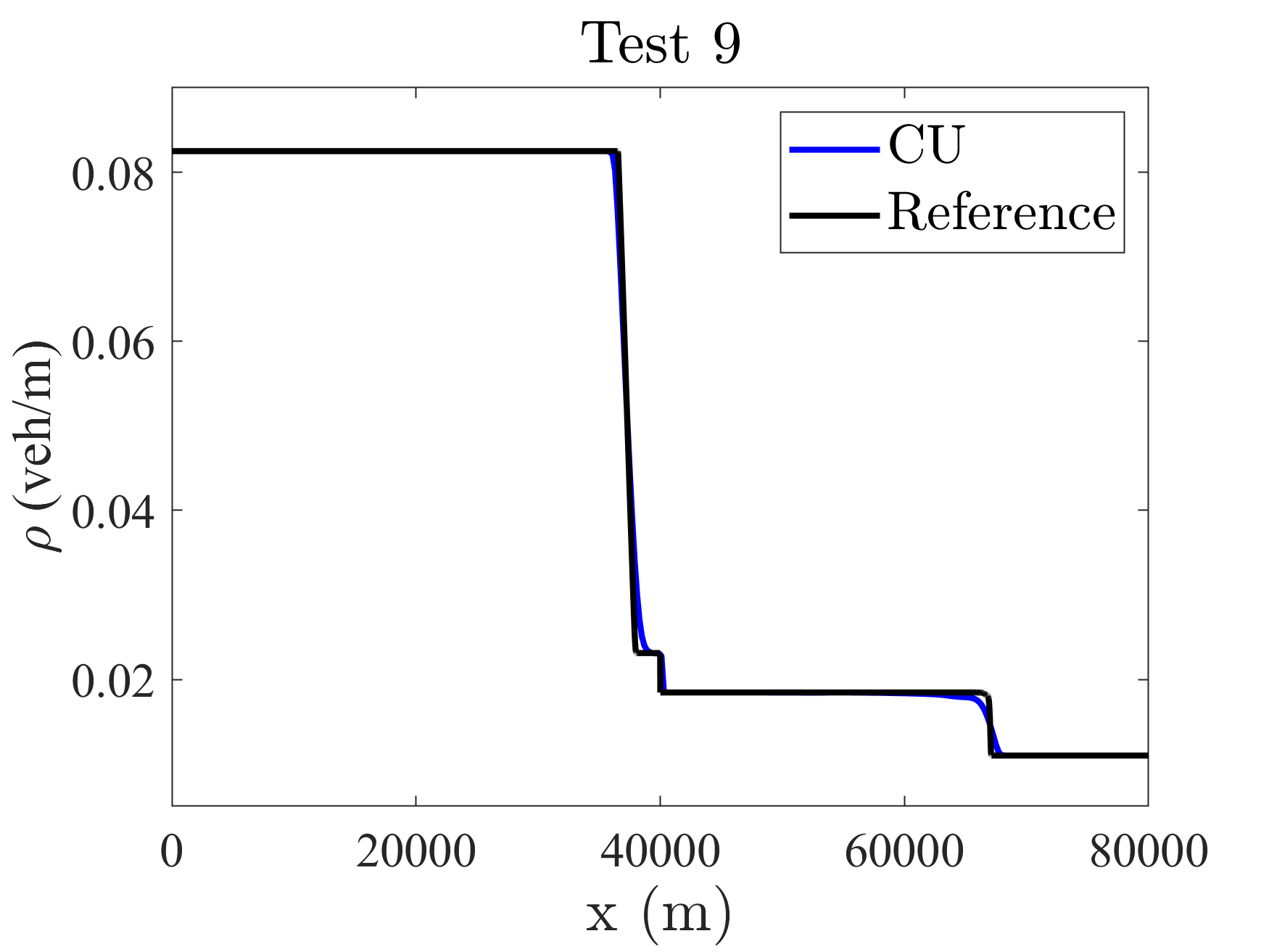}\hspace{1cm}
            \includegraphics[trim=0.1cm 0.1cm 0.7cm 0.2cm, clip, width=5.4cm]{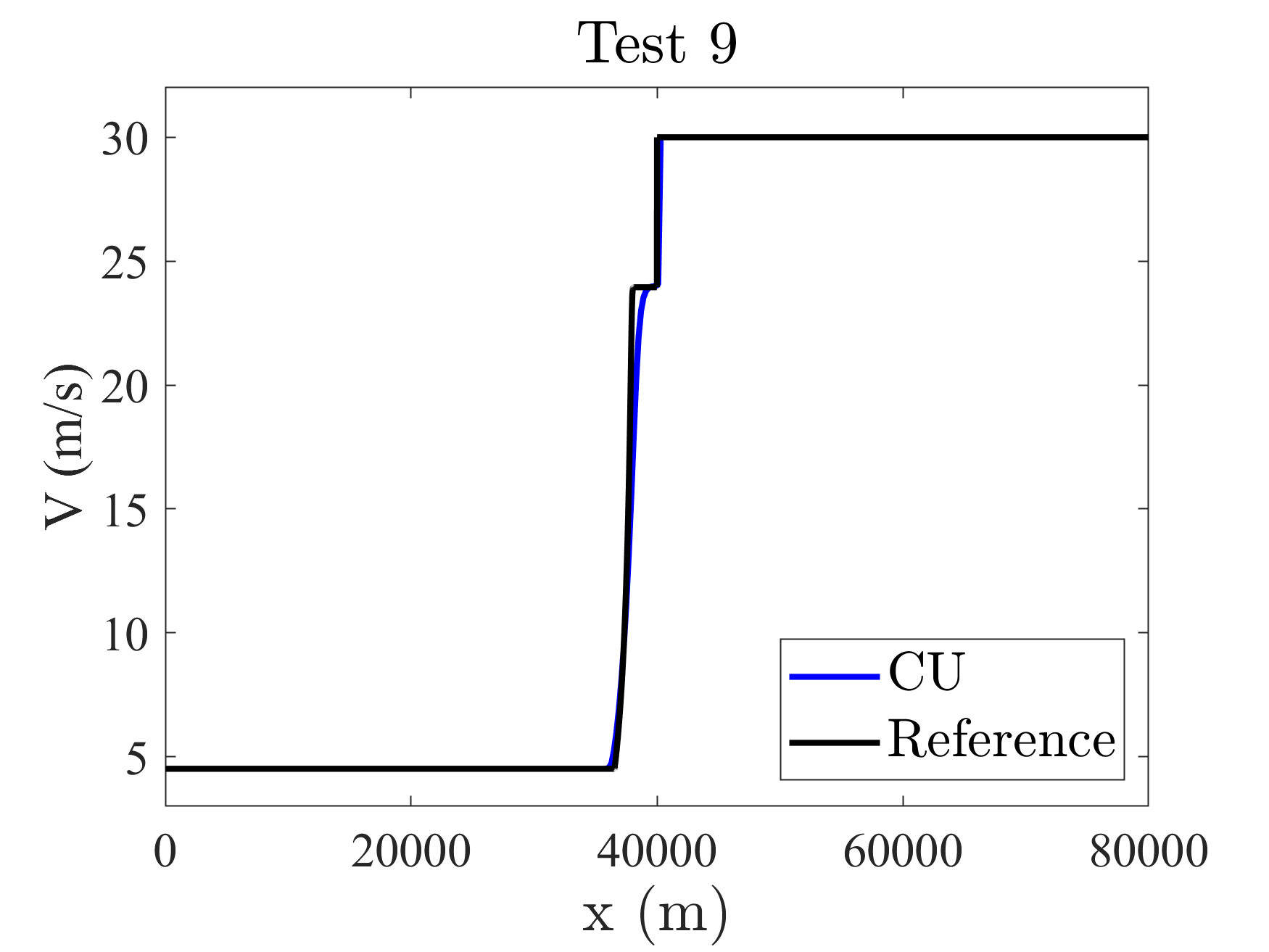}}
\vskip6pt
\centerline{\includegraphics[trim=0.1cm 0.1cm 0.7cm 0.2cm, clip, width=5.4cm]{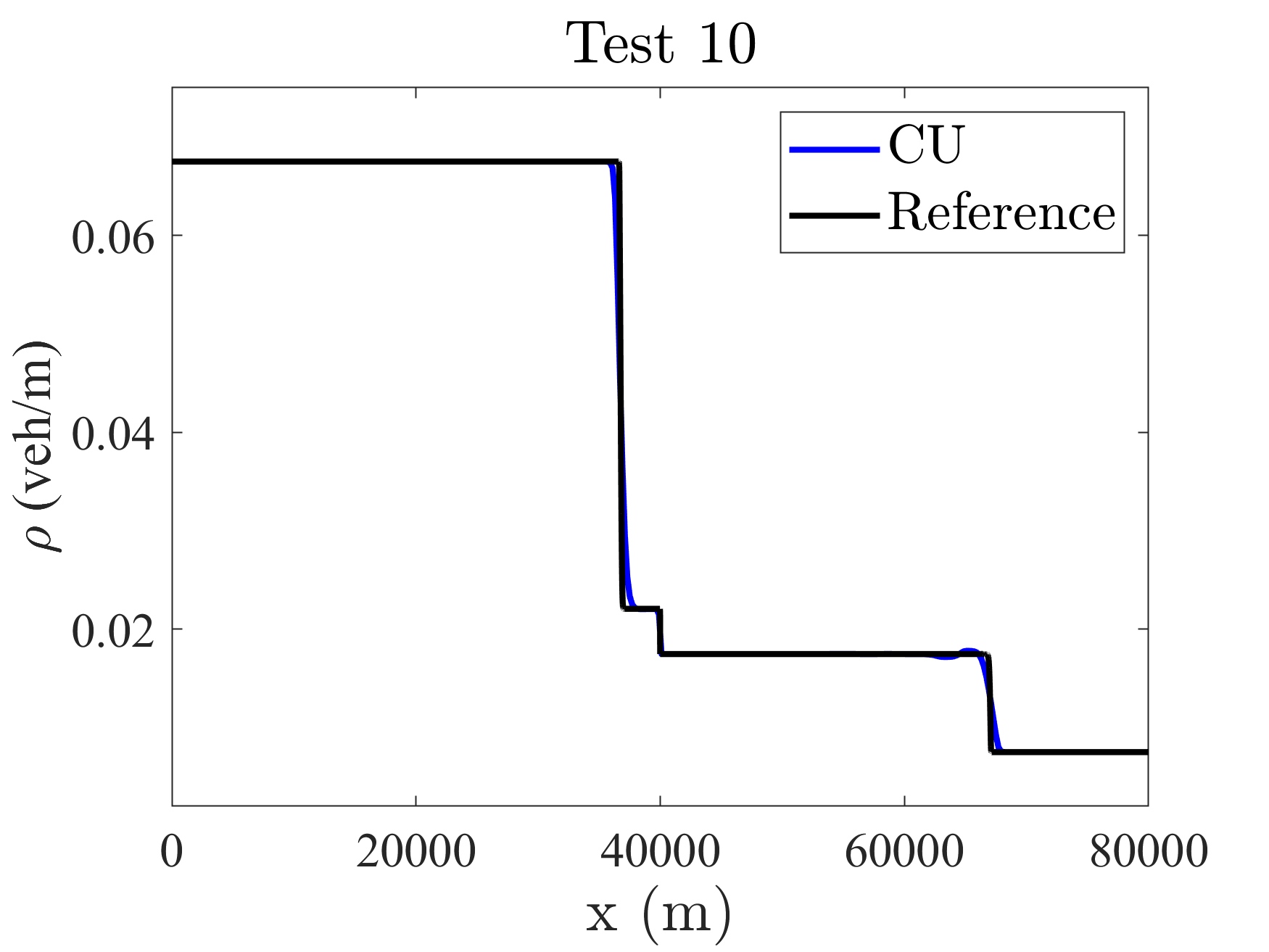}\hspace{1cm}
            \includegraphics[trim=0.1cm 0.1cm 0.7cm 0.2cm, clip, width=5.4cm]{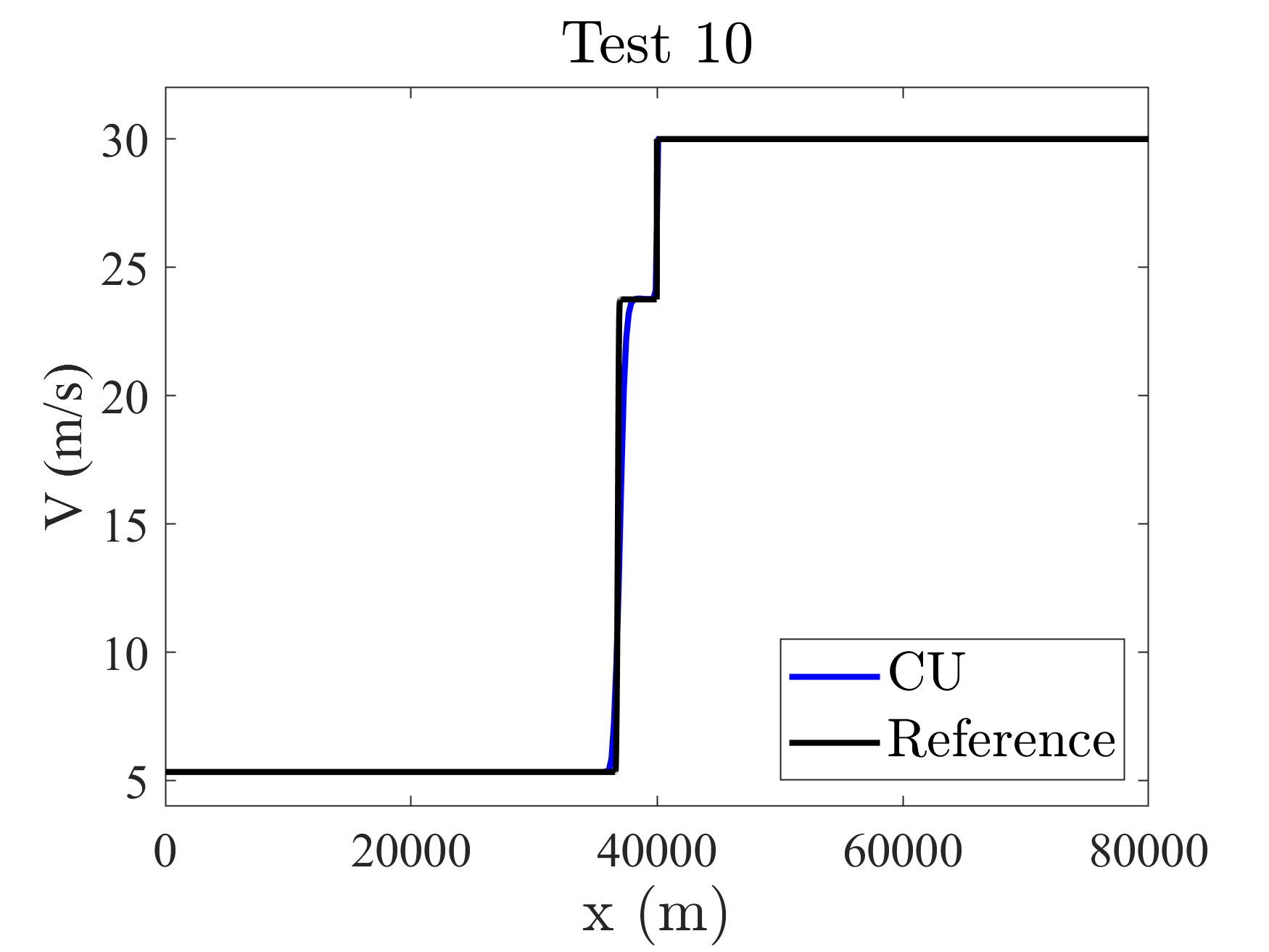}}
\vskip6pt
\centerline{\includegraphics[trim=0.1cm 0.1cm 0.7cm 0.2cm, clip, width=5.4cm]{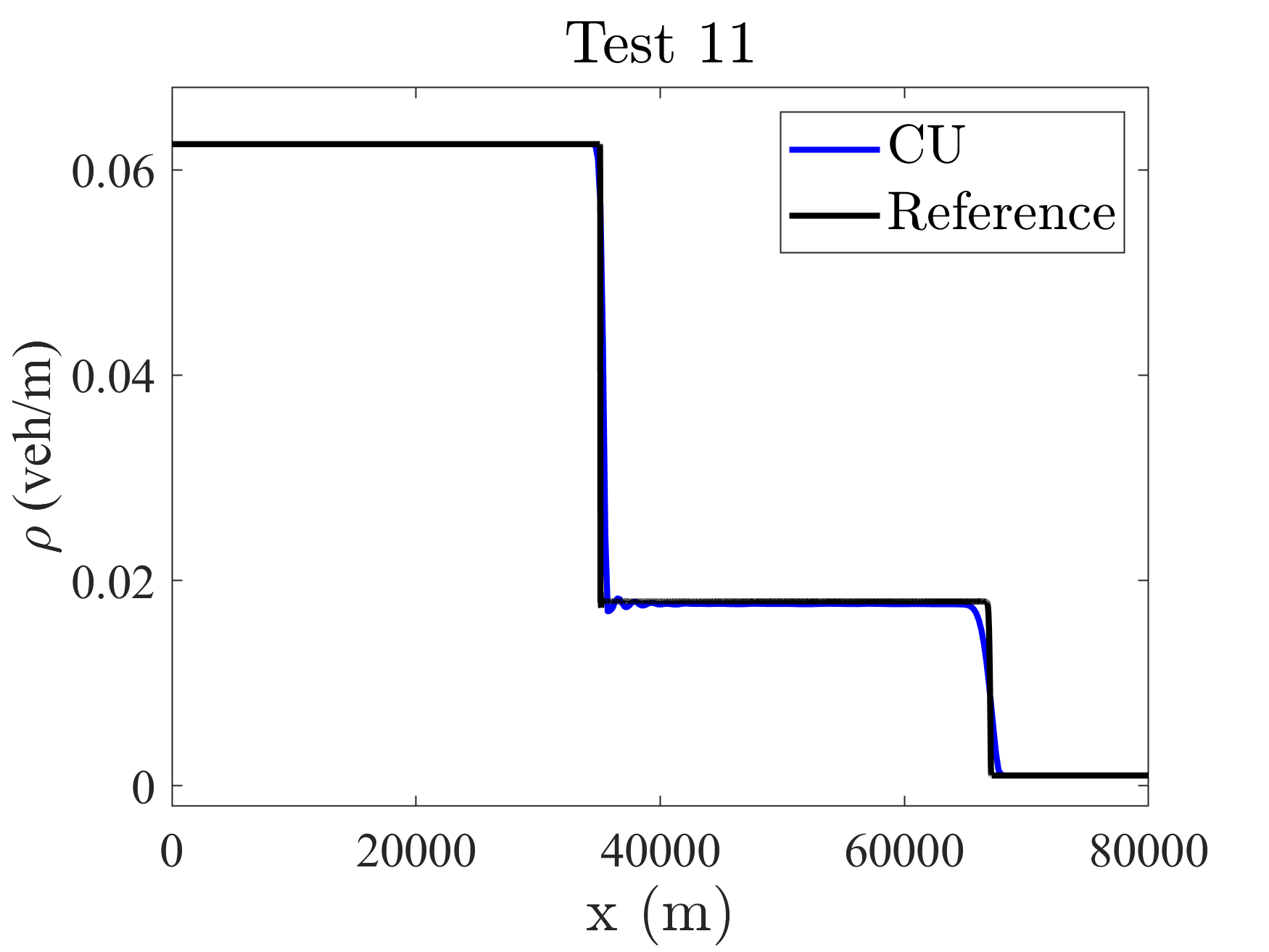}\hspace{1cm}
            \includegraphics[trim=0.1cm 0.1cm 0.7cm 0.2cm, clip, width=5.4cm]{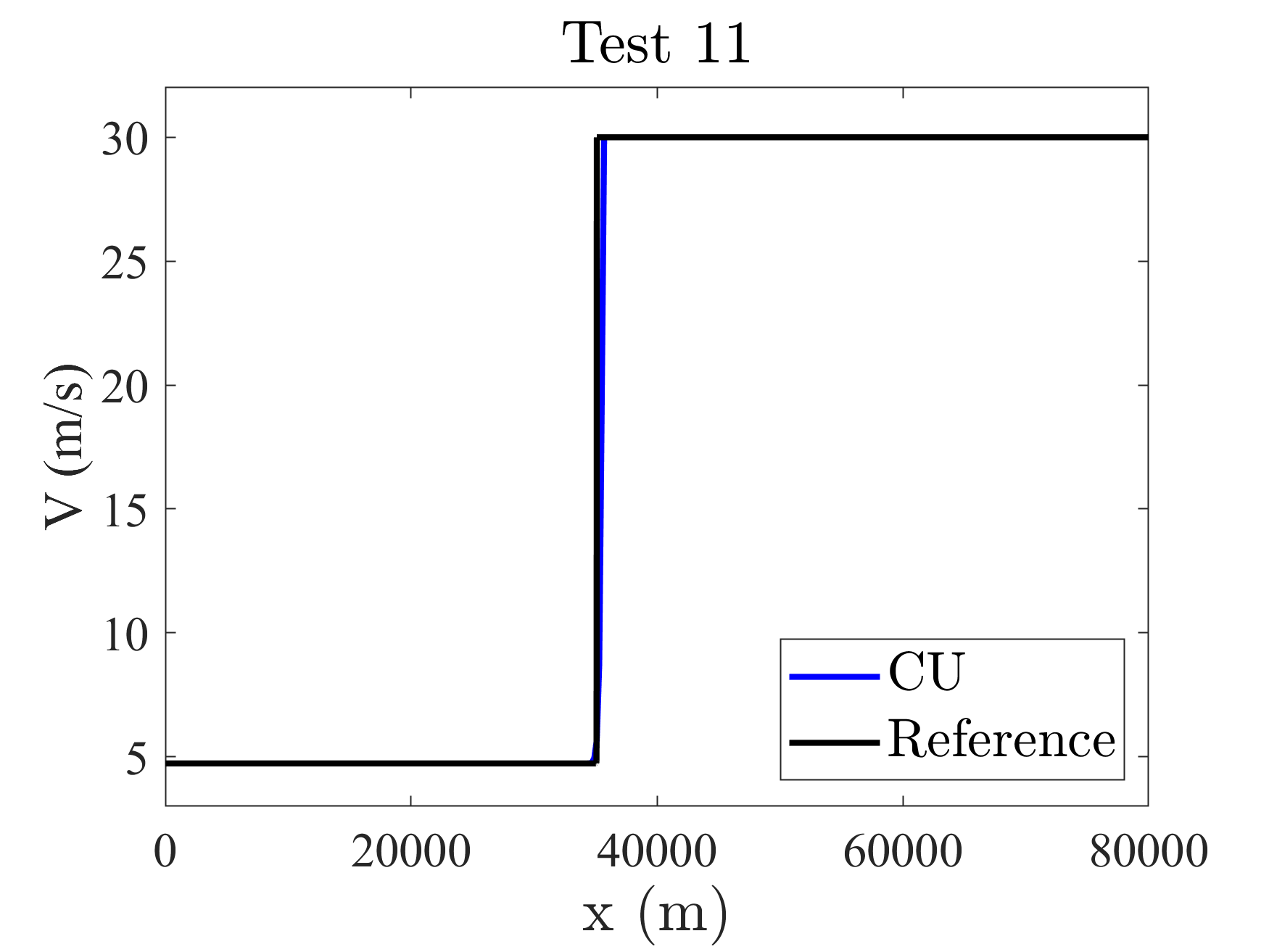}}
\vskip6pt
\centerline{\includegraphics[trim=0.1cm 0.1cm 0.7cm 0.2cm, clip, width=5.4cm]{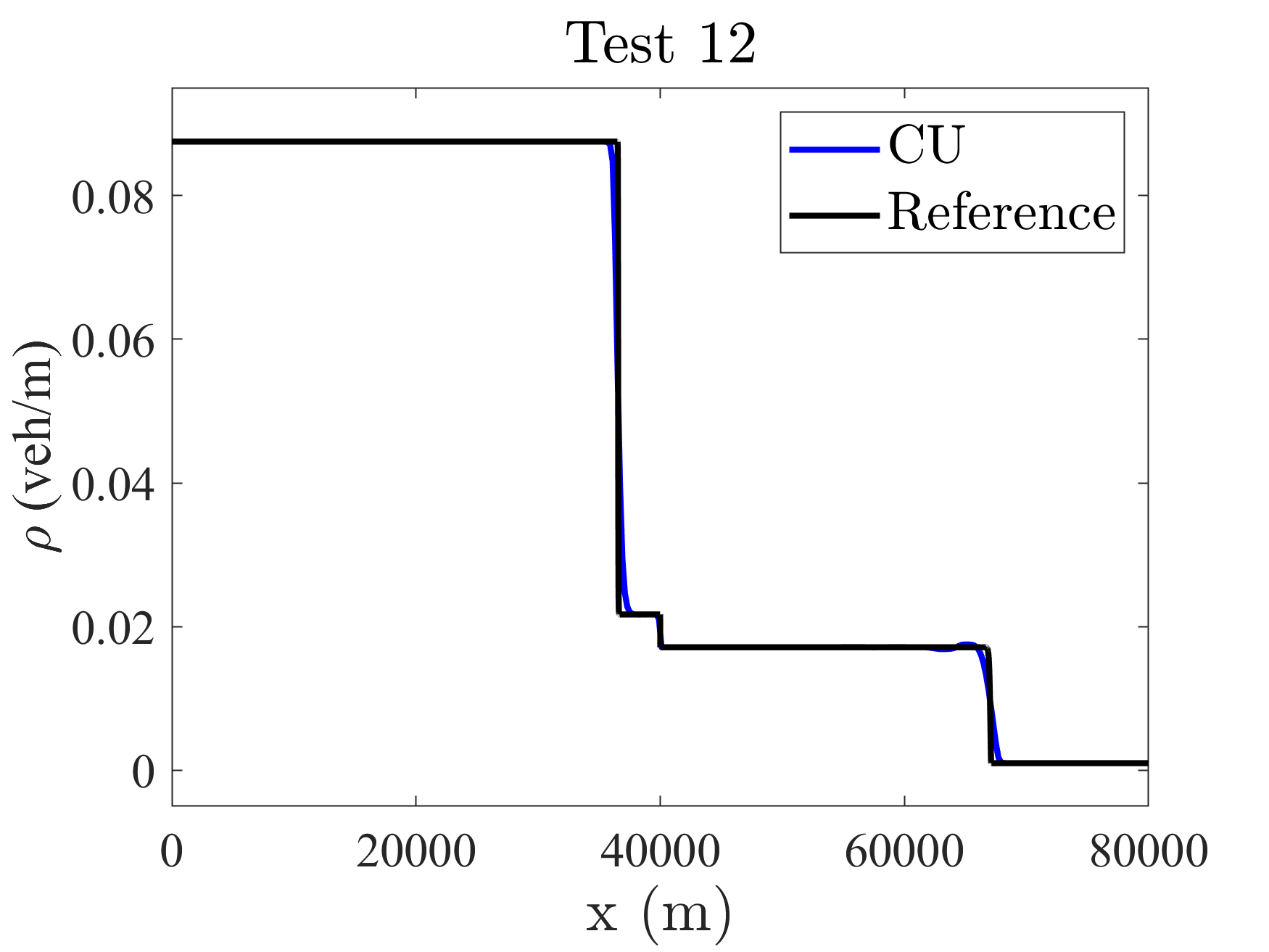}\hspace{1cm}
            \includegraphics[trim=0.1cm 0.1cm 0.7cm 0.2cm, clip, width=5.4cm]{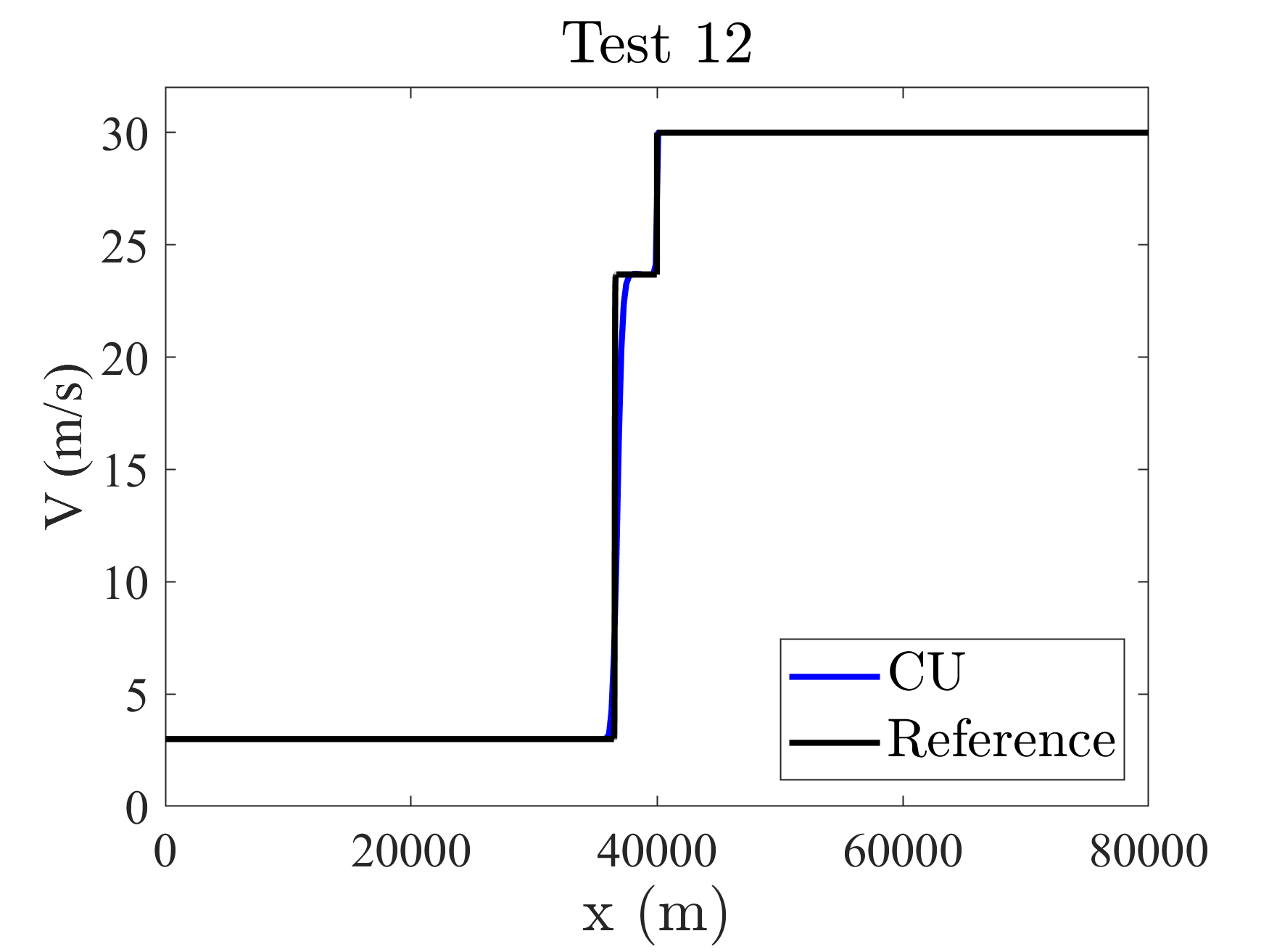}}
\caption{\sf Example 1, Tests 9--12: $\rho$ (left column) and $V$ (right column) at  time $T_{\rm final}=900$.\label{fig43}}
\end{figure}

We compute the solutions until the final time $T_{\rm final}=900$ by the proposed CU scheme on the computational domain $[0,80000]$ on the
uniform mesh with $\dx=200$ subject to the free boundary conditions. We present the obtained results in Figures \ref{fig41}--\ref{fig43}
together with the reference solution computed on a much finer mesh with $\dx=5$. As one can see, the developed CU scheme can capture the
solution structures of these Riemann problems in an non-oscillatory manner and the achieved resolution is high.

\begin{figure}[ht!]
\centerline{\includegraphics[trim=0.1cm 0.1cm 0.6cm 0.2cm, clip, width=5.4cm]{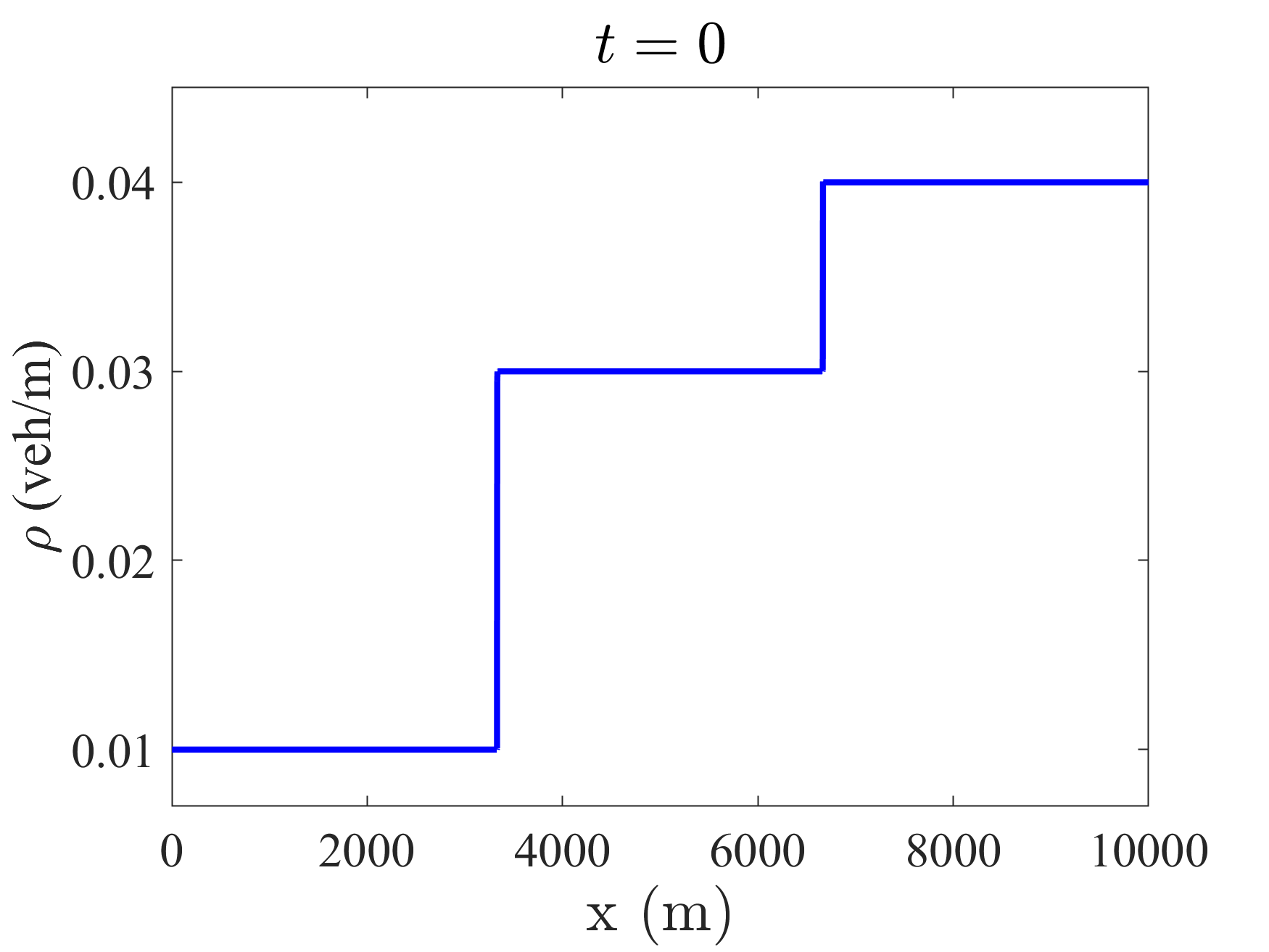}\hspace*{1cm}
            \includegraphics[trim=0.1cm 0.1cm 0.6cm 0.2cm, clip, width=5.4cm]{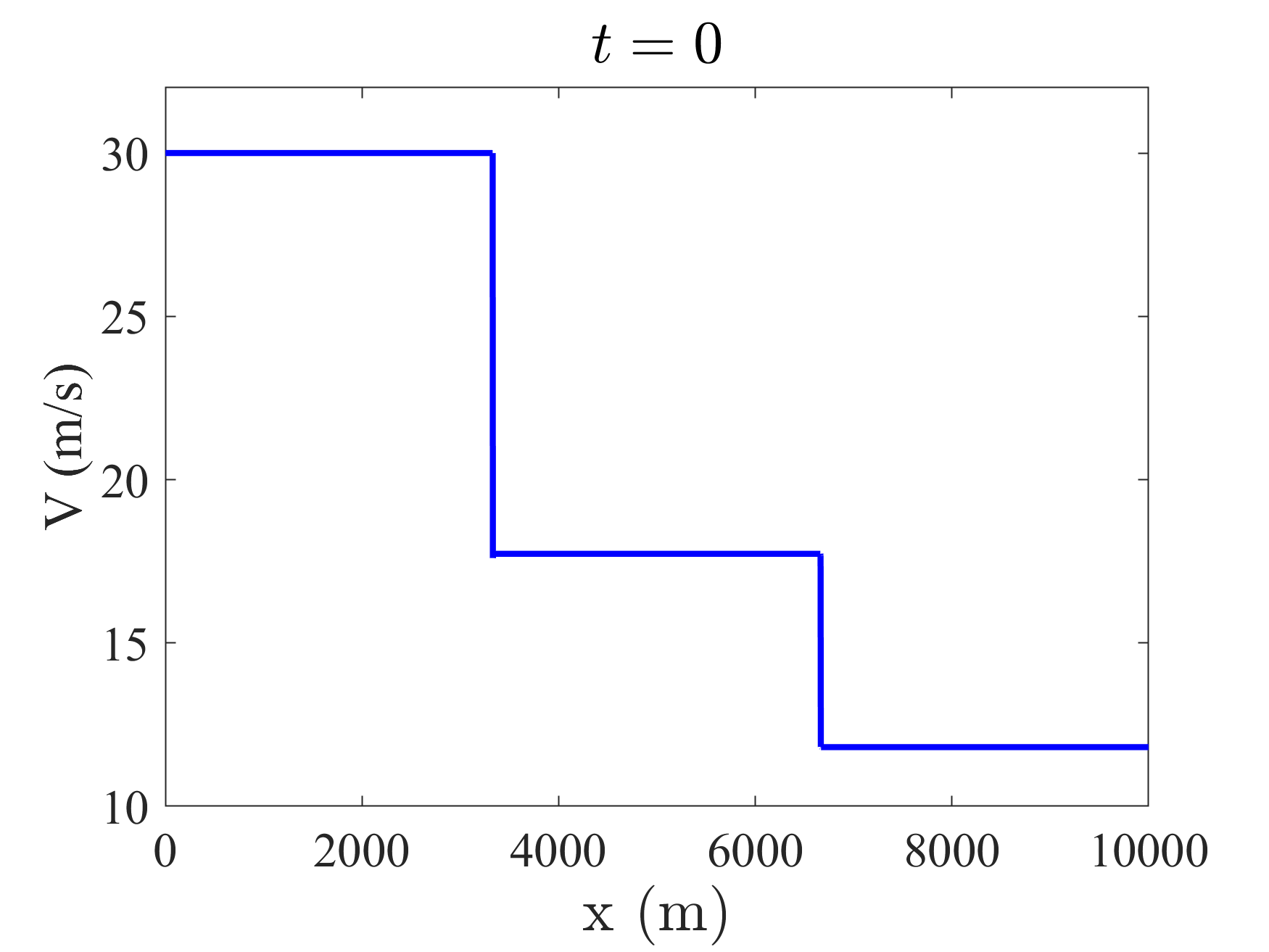}}
\vskip6pt
\centerline{\includegraphics[trim=0.1cm 0.1cm 0.6cm 0.2cm, clip, width=5.4cm]{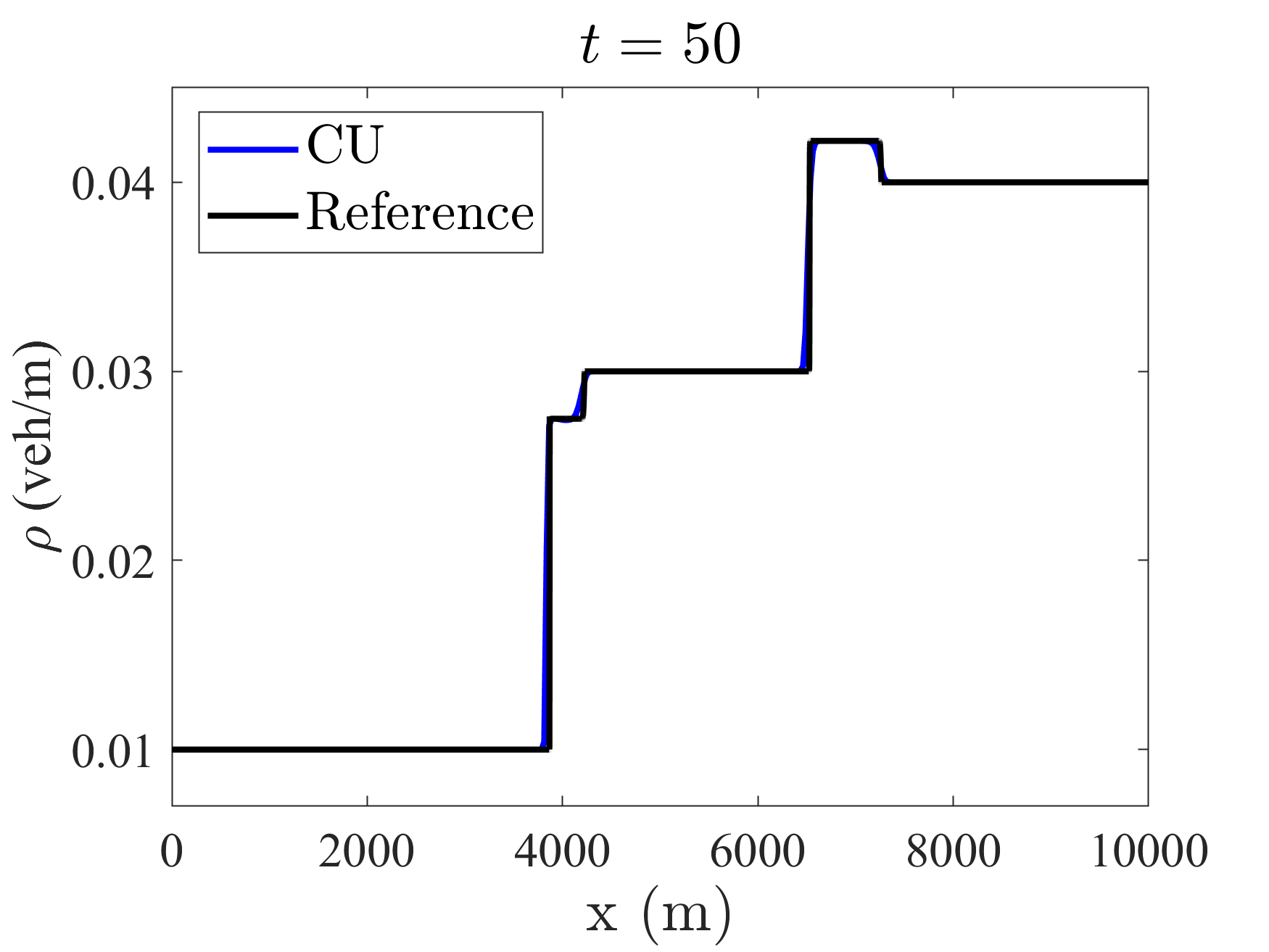}\hspace{1cm}
            \includegraphics[trim=0.1cm 0.1cm 0.6cm 0.2cm, clip, width=5.4cm]{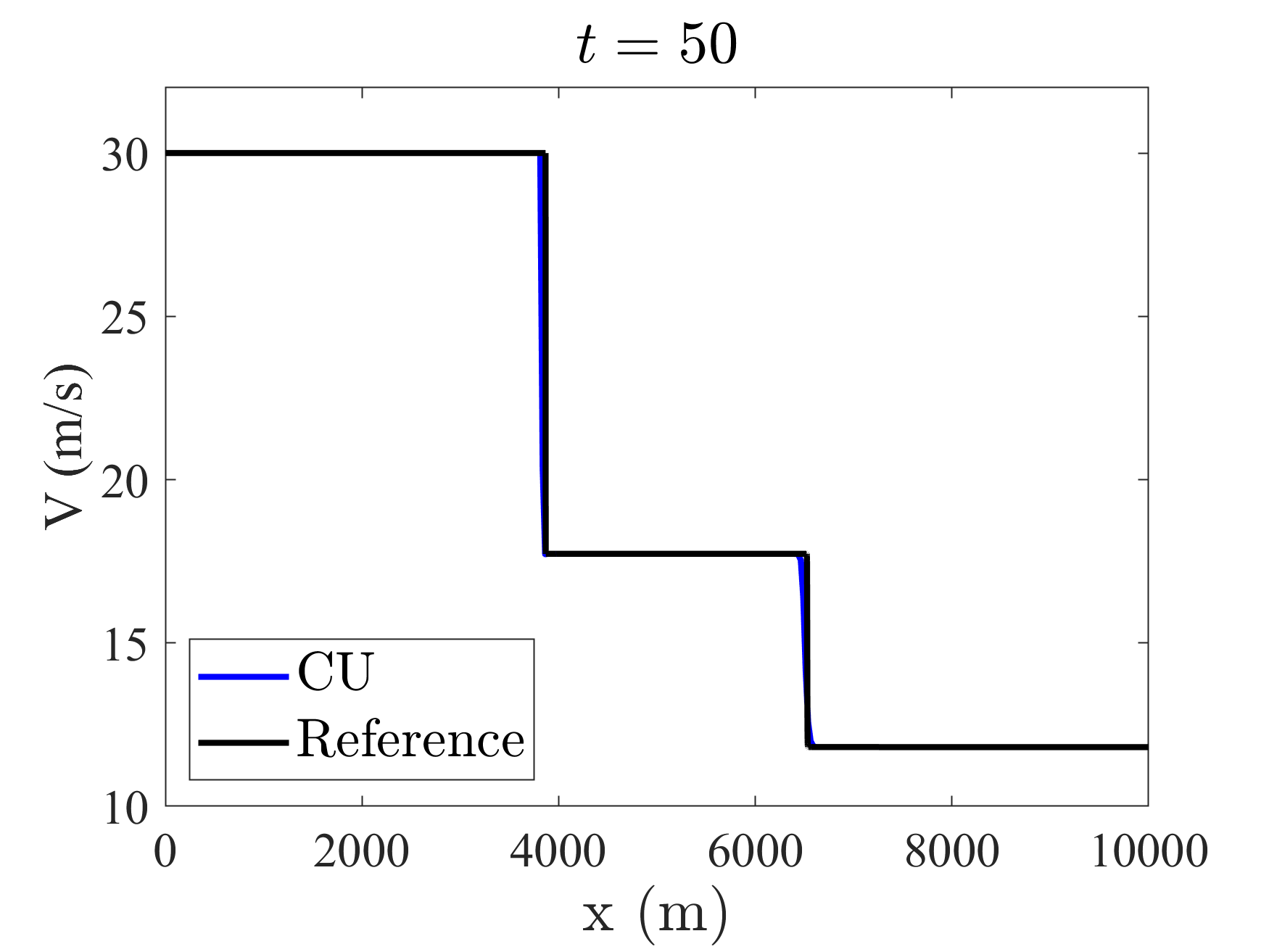}}
\vskip6pt
\centerline{\includegraphics[trim=0.1cm 0.1cm 0.6cm 0.2cm, clip, width=5.4cm]{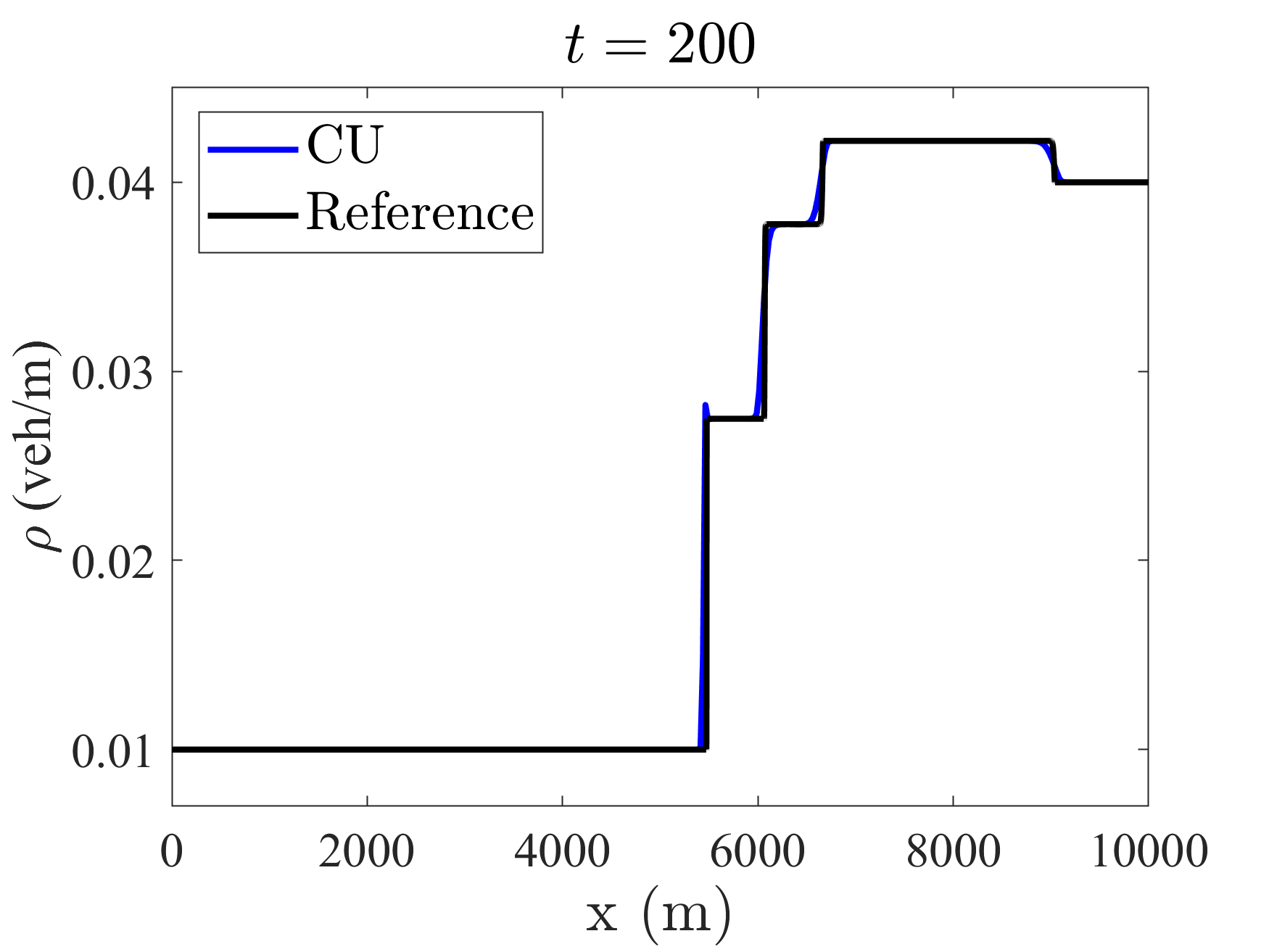}\hspace{1cm}
            \includegraphics[trim=0.1cm 0.1cm 0.6cm 0.2cm, clip, width=5.4cm]{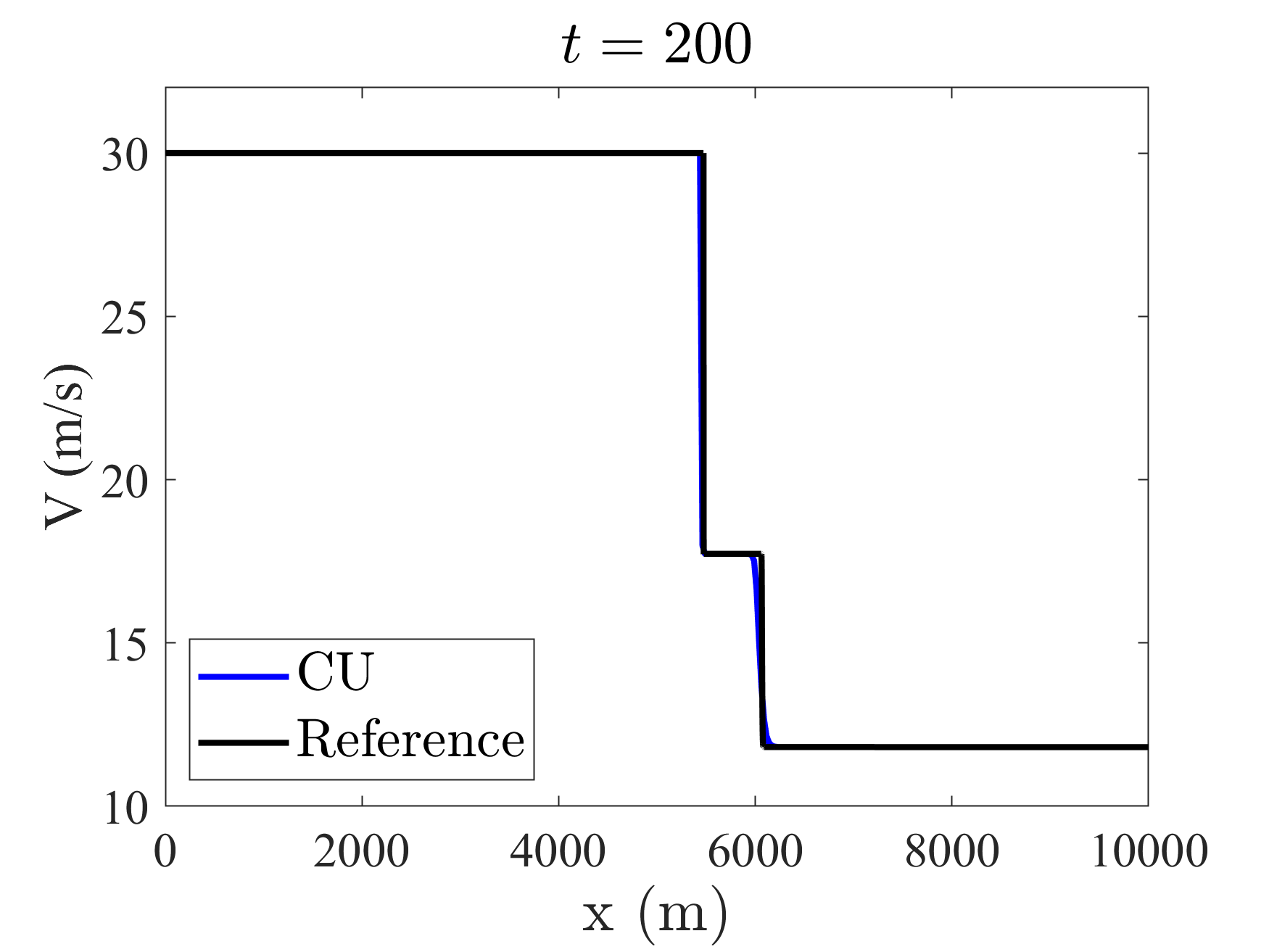}}
\vskip6pt
\centerline{\includegraphics[trim=0.1cm 0.1cm 0.6cm 0.2cm, clip, width=5.4cm]{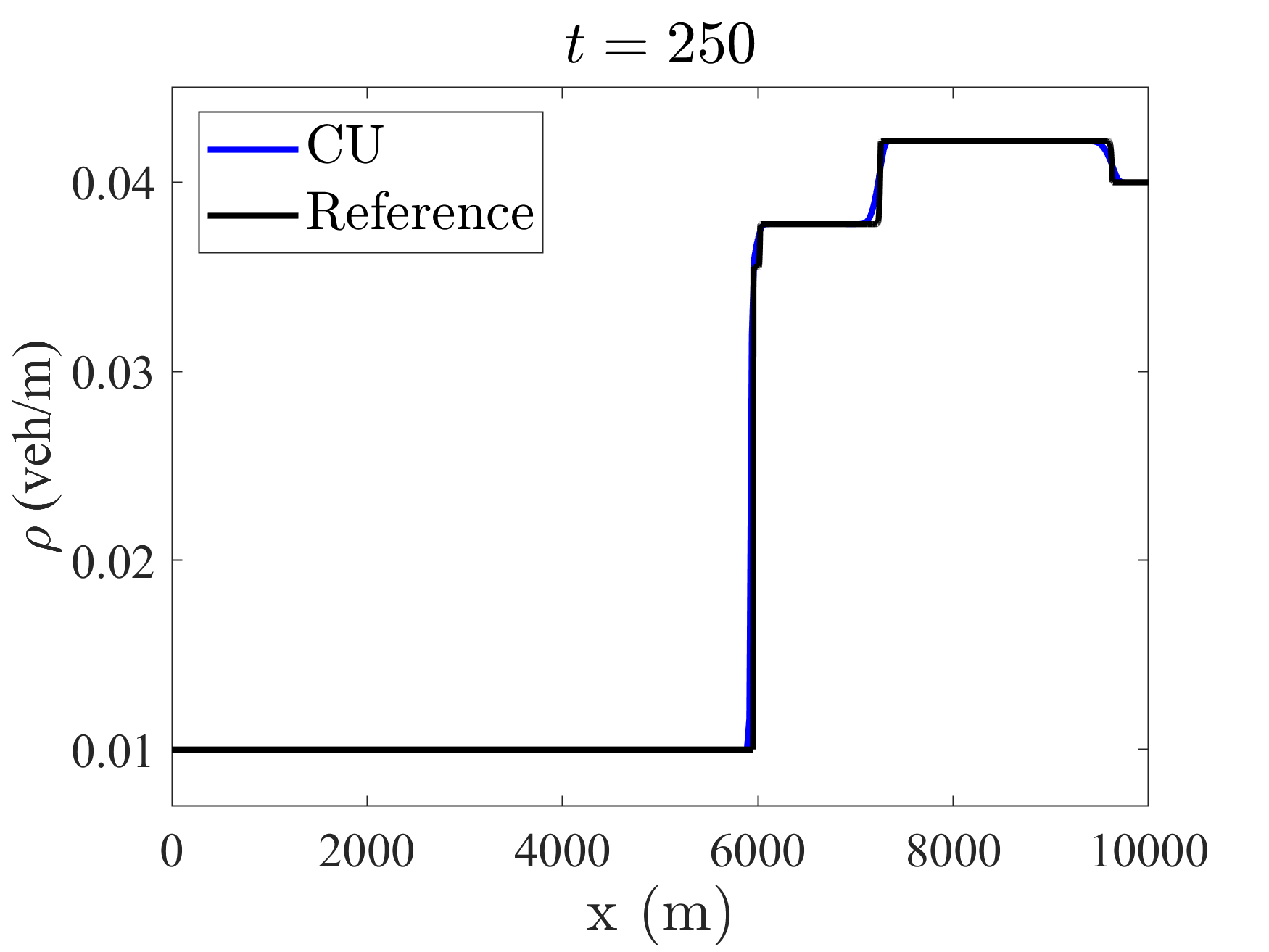}\hspace{1cm}
            \includegraphics[trim=0.1cm 0.1cm 0.6cm 0.2cm, clip, width=5.4cm]{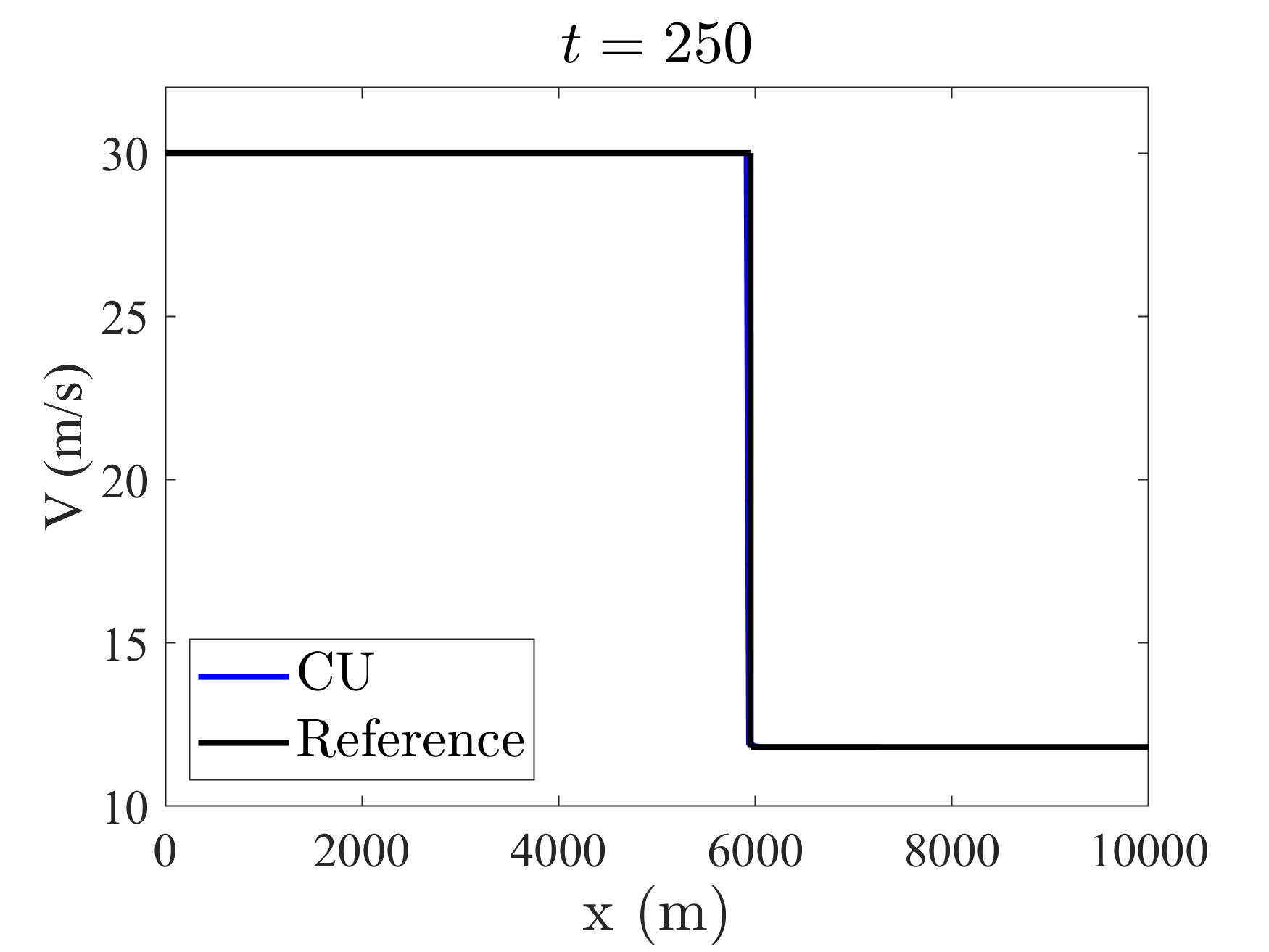}}
\caption{\sf Example 2: $\rho$ (left column) and $V$ (right column) at times $t=0$ (top row), $t=50$ (second row), $t=200$
(third row), and $t=250$ (bottom row).\label{fig46}}
\end{figure}

\paragraph*{Example 2.} In the second example, we take the computational domain $[0,L]$ with $L=10000$ and consider the following initial
conditions:
\begin{equation}
\rho(x,0)=\begin{cases}
0.01&\mbox{if}~x\le\frac{L}{3},\\
0.03&\mbox{if}~\frac{L}{3}<x\le\frac{2L}{3},\\
0.04&\mbox{otherwise},
\end{cases}\quad
V(x,0)=\begin{cases}
30&\mbox{if}~x\le\frac{L}{3},\\
17.729&\mbox{if}~\frac{L}{3}<x\le\frac{2L}{3},\\
11.812&\mbox{otherwise},
\end{cases}
\end{equation}
shown in Figure \ref{fig46} (top row). The corresponding values of the quantity $q-q^*$ are then
\begin{equation}
q(x,0)-q^*=\begin{cases}
-0.2800&\mbox{if}~x\le\frac{L}{3},\\
0.0546&\mbox{if}~\frac{L}{3}< x\le\frac{2L}{3},\\
0.0300&\mbox{otherwise}.
\end{cases}
\end{equation}

This test case is numerically challenging and physically complex due to the formation and propagation of multiple interacting compound waves
traveling in different directions. Specifically, the traffic dynamics vary across the domain: the region $x\le\frac{L}{3}$ exhibits
free-flowing traffic, which progresses to moderate congestion in the range $\frac{L}{3}<x\le\frac{2L}{3}$, and escalates to dense traffic
when $x>\frac{2L}{3}$. The initial conditions feature two discontinuities in speed and density profiles. As $q(x,0)-q^*\ne0$ for all $x$,
compound waves are anticipated at the interfaces of these discontinuities.

We compute the numerical results using the proposed CU scheme until the final time $T_{\rm final}=250$ on a uniform mesh with $\dx=25$
subject to the following Dirichlet boundary conditions imposed at the left end of the computational domain $x=0$:
\begin{equation}
\rho(0,t)=0.01,\quad q(0,t)=\frac{\rho(0,t)V_{\max}}{1-\rho(0,t)/\rho_{\max}},
\end{equation}
and free boundary conditions at the right end of the computational domain.

The computed $\rho$ and $V$ at times $t=50$, 200, and 250 are presented in Figure \ref{fig46} together with a reference solution calculated
on a finer mesh with $\dx=\frac{5}{4}$. The comparison underscores the scheme's robust capability to capture solution structures arising
from interacting waves. It maintains sharpness near contact discontinuities and shock waves, while ensuring non-oscillatory behavior.

\begin{figure}[ht!]
\centering
\begin{minipage}{0.5\textwidth}
\includegraphics[trim=0.1cm 0.1cm 0cm 0.2cm, clip, width=\linewidth]{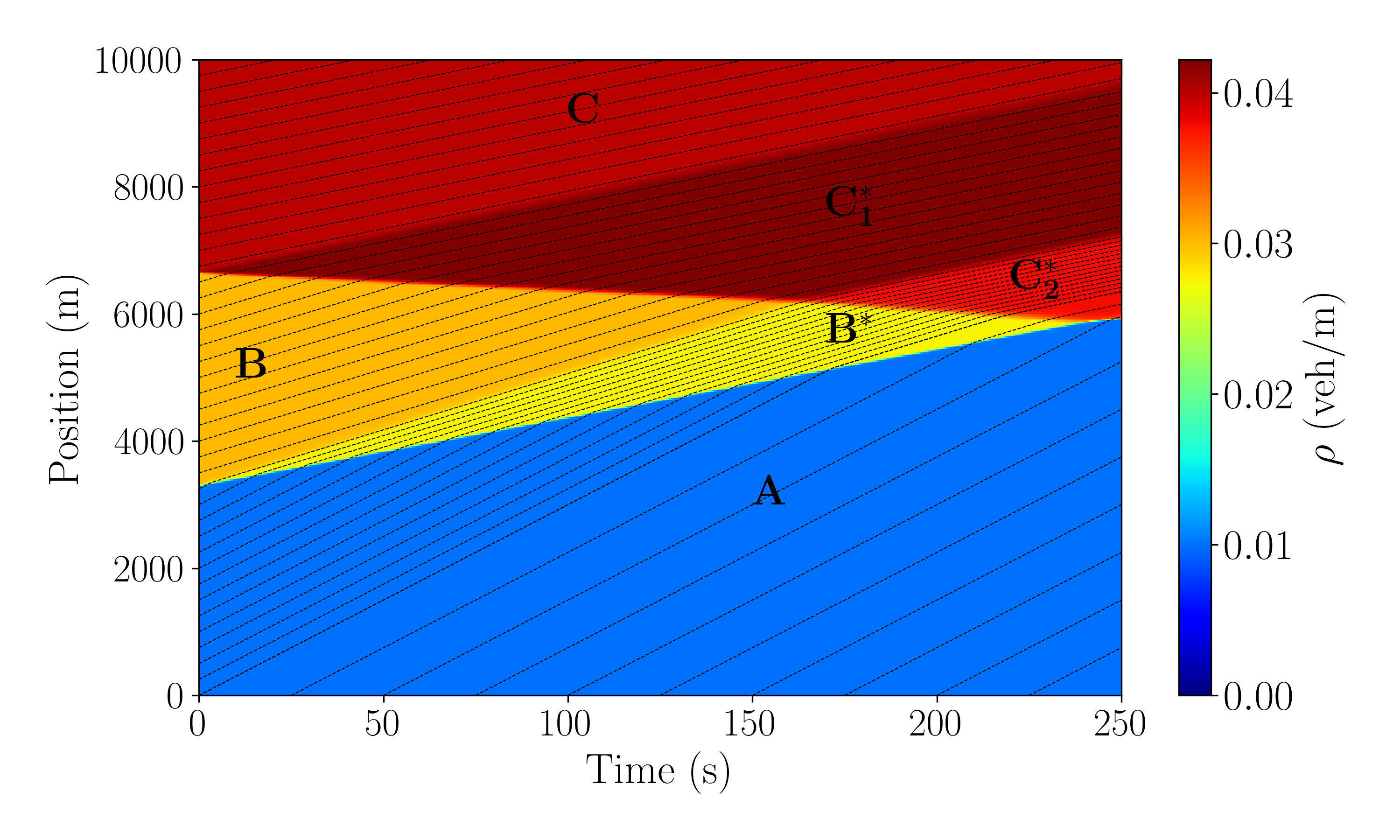}
\end{minipage}%
\begin{minipage}{0.5\textwidth}
\includegraphics[trim=0.1cm 0.1cm 0cm 0.2cm, clip, width=\linewidth]{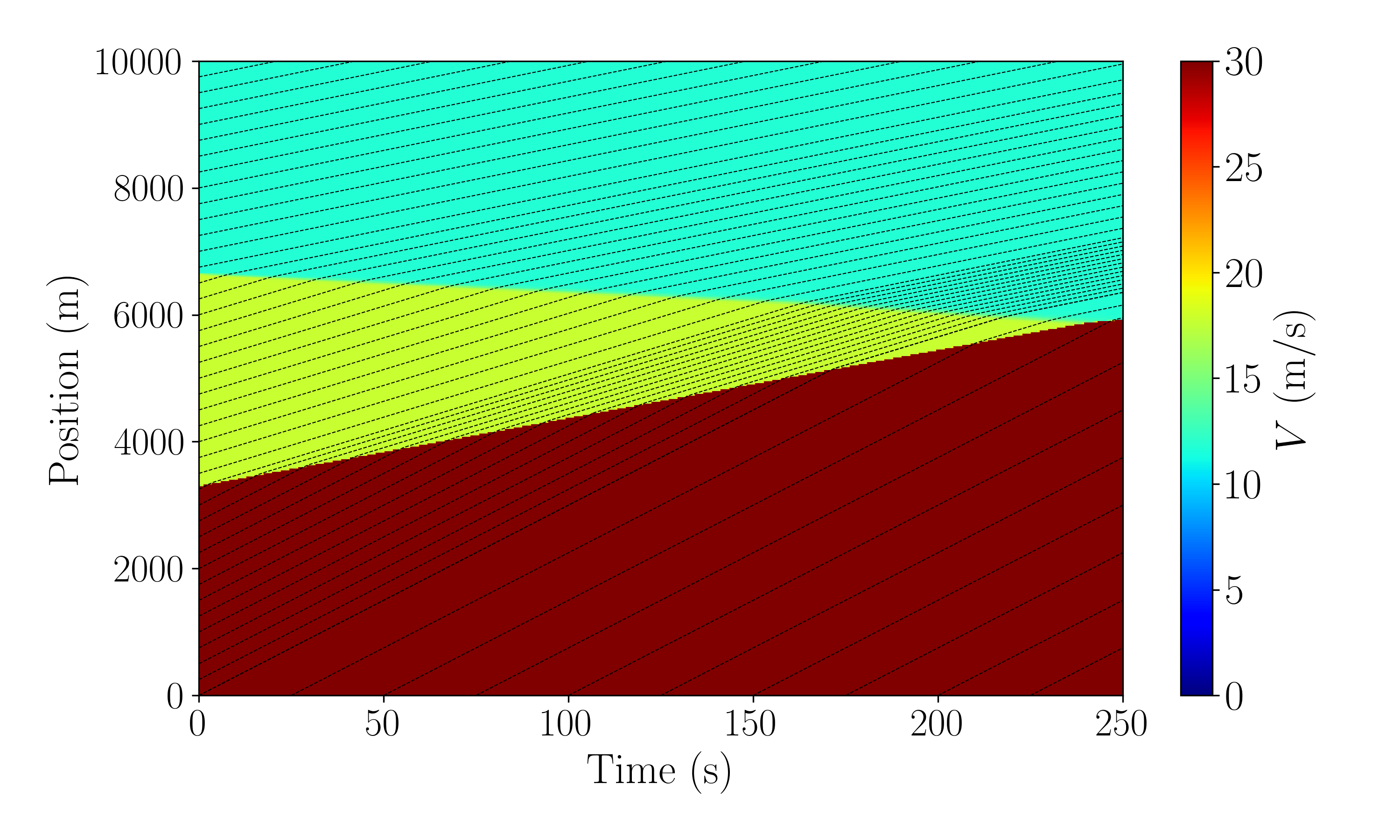}
\end{minipage}
\caption{\sf Example 2: Spatio-temporal evolution of $\rho$ (left) and $V$ (right) for $t\in [0,250]$.\label{Example3_ST}}
\end{figure}

Let us now examine the spatio-temporal evolution of traffic density and speed. In Figure \ref{Example3_ST}, we present
time-space diagrams for both variables, with representative vehicle trajectories superimposed---both from vehicles placed every 100$\,$m
along the road at the initial time and those released every 50$\,$s from the upstream boundary at $x=0$. These vehicle trajectories can be
viewed as phantom moving observers traveling with the local average traffic speed and adapting their dynamics to prevailing flow
conditions. The spatio-temporal analysis of traffic density indicates that it can be divided into six distinct regions (each characterized
by a constant density), which are labeled as A, B, ${\rm B}^*$, C, ${\rm C}^*_1$, and ${\rm C}^*_2$ in Figure \ref{Example3_ST} (left). Note
that the traffic speed remains invariant across Regions B, ${\rm B}^*$, and also across C, ${\rm C}^*_1$, and ${\rm C}^*_2$. We shall
elaborate on the waves arising at the interface between these regions.

In order to interpret the results physically, we consider the region defined by $\frac{L}{3}<x< \frac{2L}{3}$ in the initial conditions as
the reference point, referred to as Region B in Figure \ref{Example3_ST} (left). In Region B, $q(x,0)-q^*>0$, suggesting that the
inter-vehicular spacing is initially smaller than the equilibrium spacing, and consequently, the traffic density is higher than would be
expected with initial $q-q^*=0$. At time $t=0$, Region B is flanked upstream by a free-flow region ($x<\frac{L}{3}$, Region A) with very low
density, and downstream by a denser traffic condition ($x>\frac{2L}{3}$, Region C). As a result, compound waves form at both ends of Region
B. To facilitate discussion on these waves, it is important to clarify that henceforth, the endpoints of each region in the spatial domain
will be referred to as the right-end (downstream) and left-end (upstream), respectively.

At the left end of Region B, an intermediate state (Region ${\rm B}^*$) emerges and separates Regions A and B with a shockwave at the left
end of Region ${\rm B}^*$ (corresponding to $\lambda_1$ and traveling forward due to the lower density in Region A) and a contact
discontinuity at the right end (corresponding to $\lambda_2$). In Region ${\rm B}^*$, the traffic speed matches that of Region B, while the
density is intermediate between that of Regions A and B. Since the speed of the shockwave is smaller in magnitude compared to the speed of
the contact discontinuity at the interface between Regions ${\rm B}^*$ and B, Region ${\rm B}^*$ expands spatially over time. This expansion
continues until the contact discontinuity ceases to exists around $x=6200\,\mbox{m}$ and $t=160\,\mbox{s}$.

At the right end of Region B, another intermediate state, Region ${\rm C}^*_1$, emerges, separating Regions B and C. Region ${\rm C}^*_1$ is
characterized by a shockwave at the left end (corresponding to $\lambda_1$ and traveling backward due to congestion) and a contact
discontinuity at the right end (corresponding to $\lambda_2$). Within Region ${\rm C}^*_1$, the traffic speed is consistent with that of
Region C, while the density is intermediate between those in Regions B and C.

At a specific point in time and space (specifically, around $x=6200\,\mbox{m}$ and $t=160\,\mbox{s}$), the shockwave traveling backward at
the left end of Region ${\rm C}^*_1$ meets the contact discontinuity traveling forward at the right end of Region ${\rm B}^*$. This
interaction leads to the formation of Region ${\rm C}^*_2$, which is confined by a shockwave traveling backward and a contact discontinuity
traveling forward.

\begin{figure}[ht!]
\centerline{\includegraphics[trim=0.1cm 0.1cm 0.6cm 0.2cm, clip, width=5.4cm]{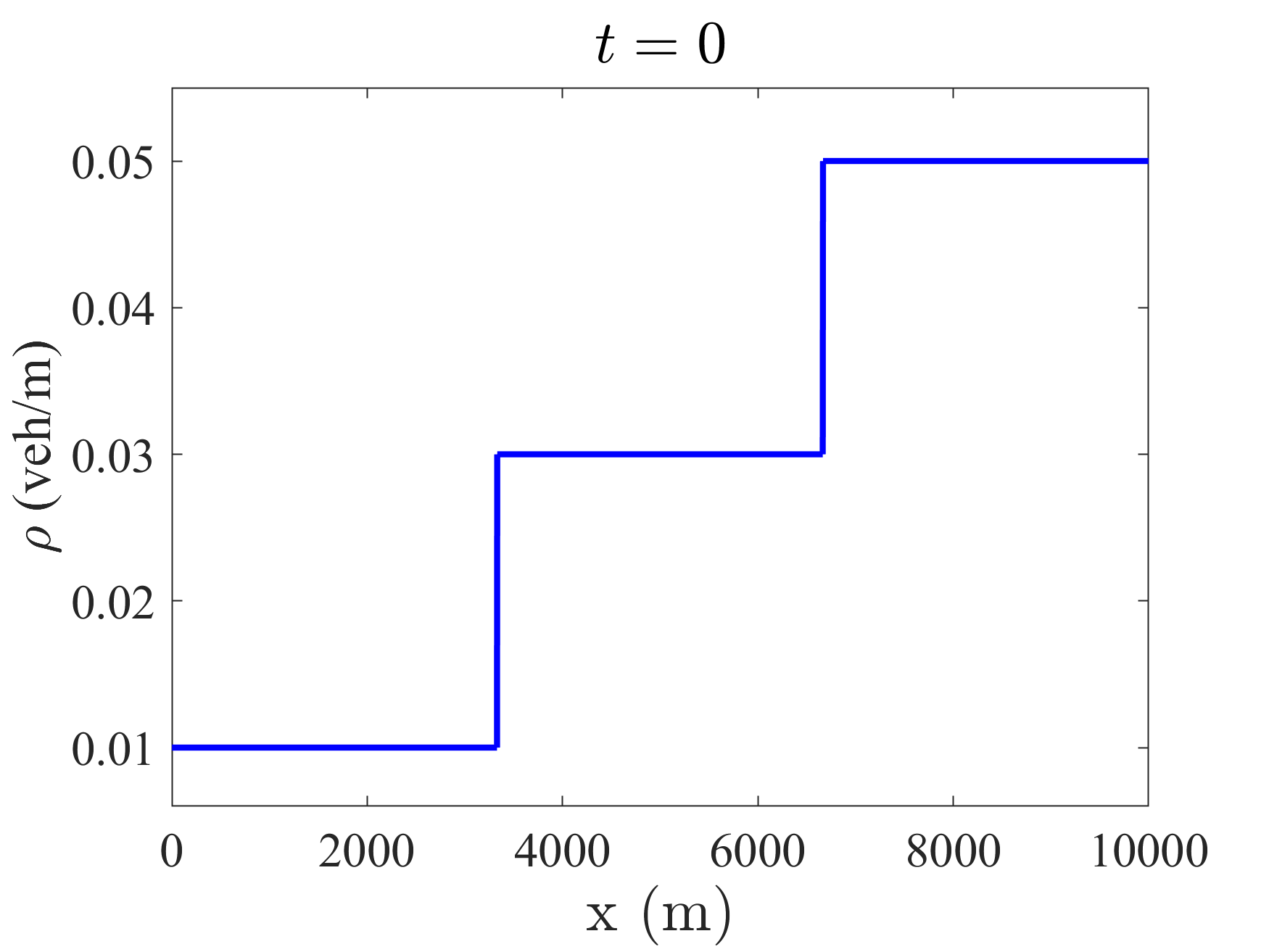}\hspace*{1cm}
            \includegraphics[trim=0.1cm 0.1cm 0.6cm 0.2cm, clip, width=5.4cm]{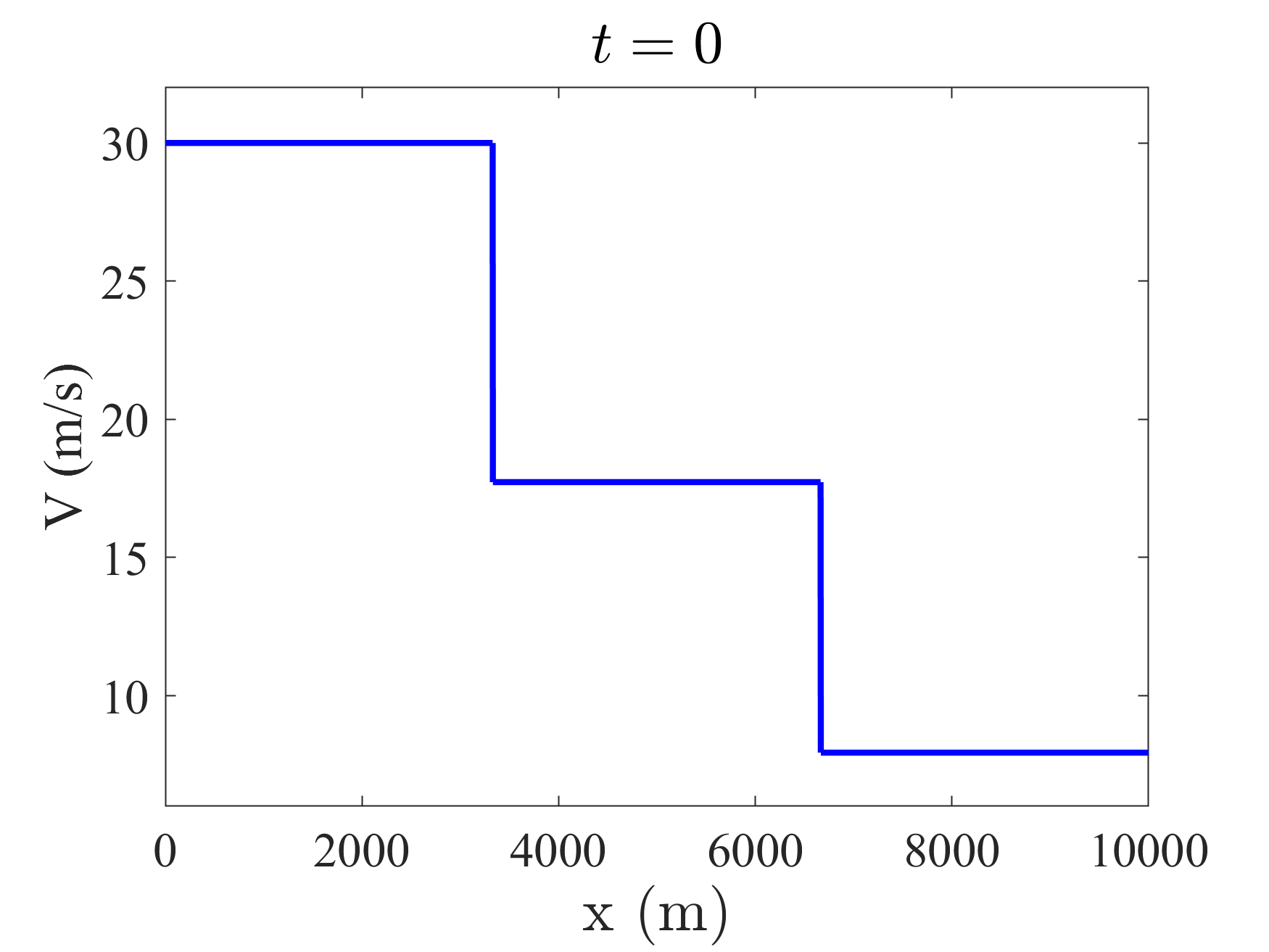}}
\vskip3pt           
\centerline{\includegraphics[trim=0.1cm 0.1cm 0.6cm 0.2cm, clip, width=5.4cm]{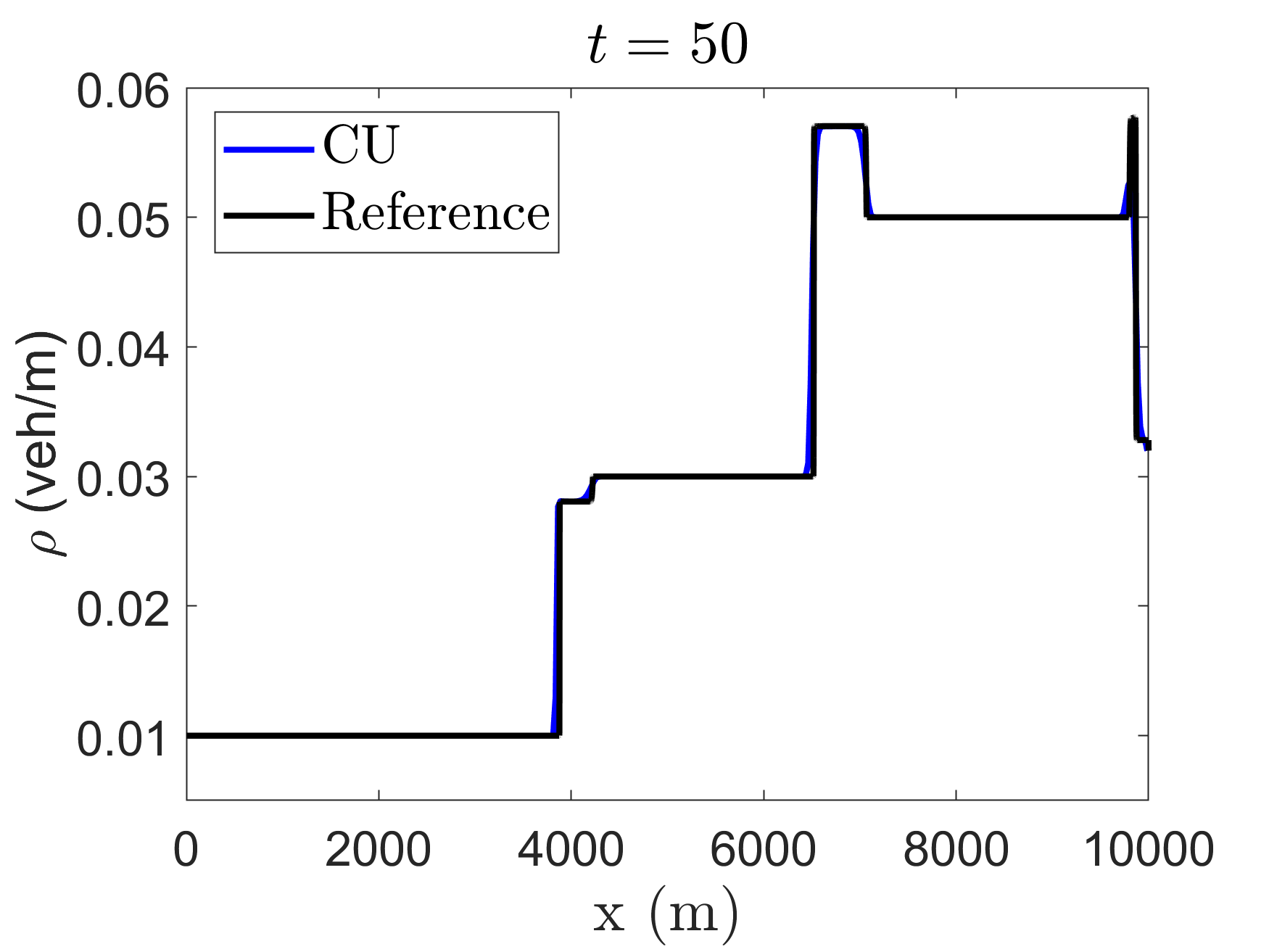}\hspace*{1cm}
            \includegraphics[trim=0.1cm 0.1cm 0.6cm 0.2cm, clip, width=5.4cm]{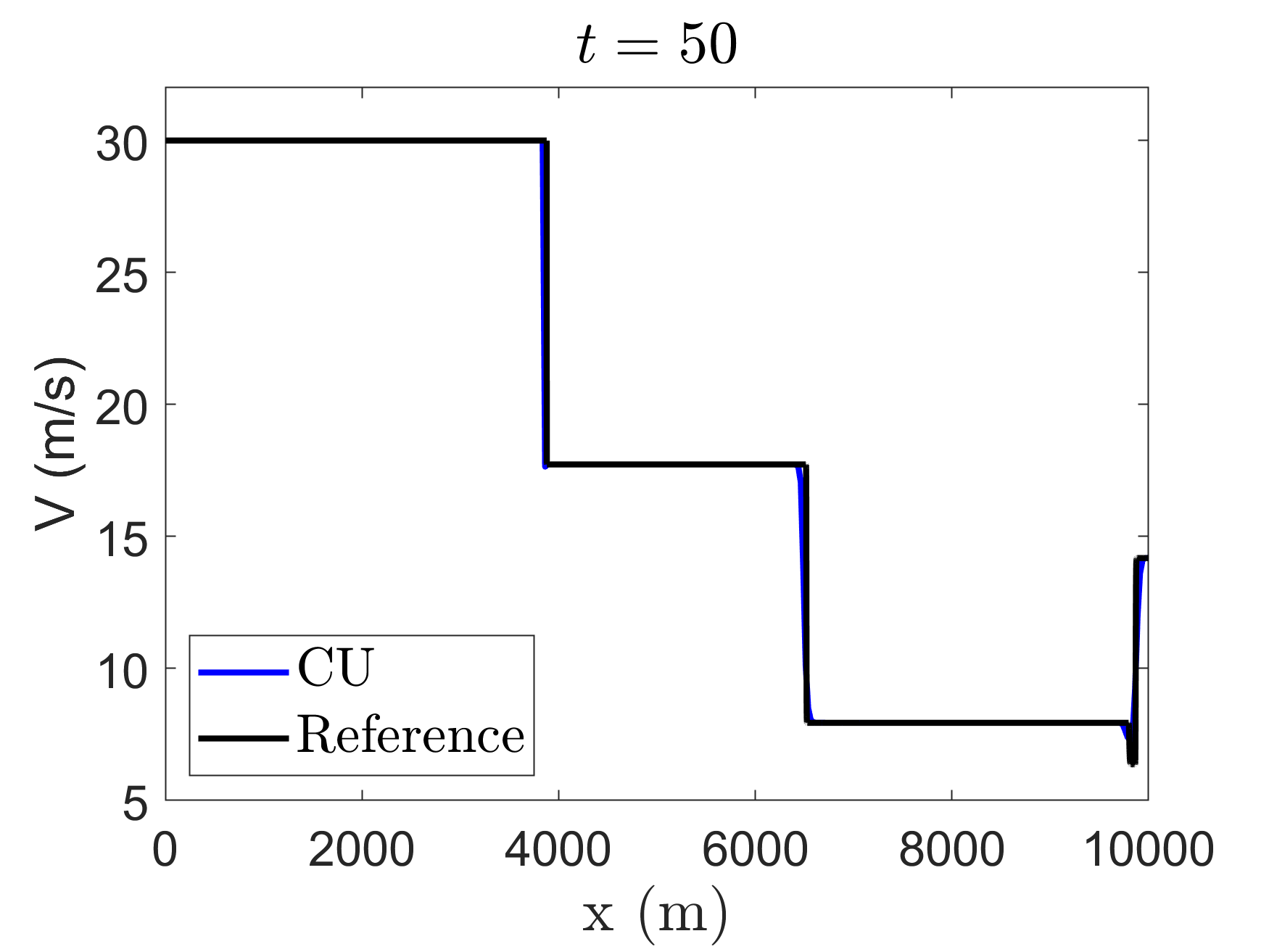}}
\vskip3pt
\centerline{\includegraphics[trim=0.1cm 0.1cm 0.6cm 0.2cm, clip, width=5.4cm]{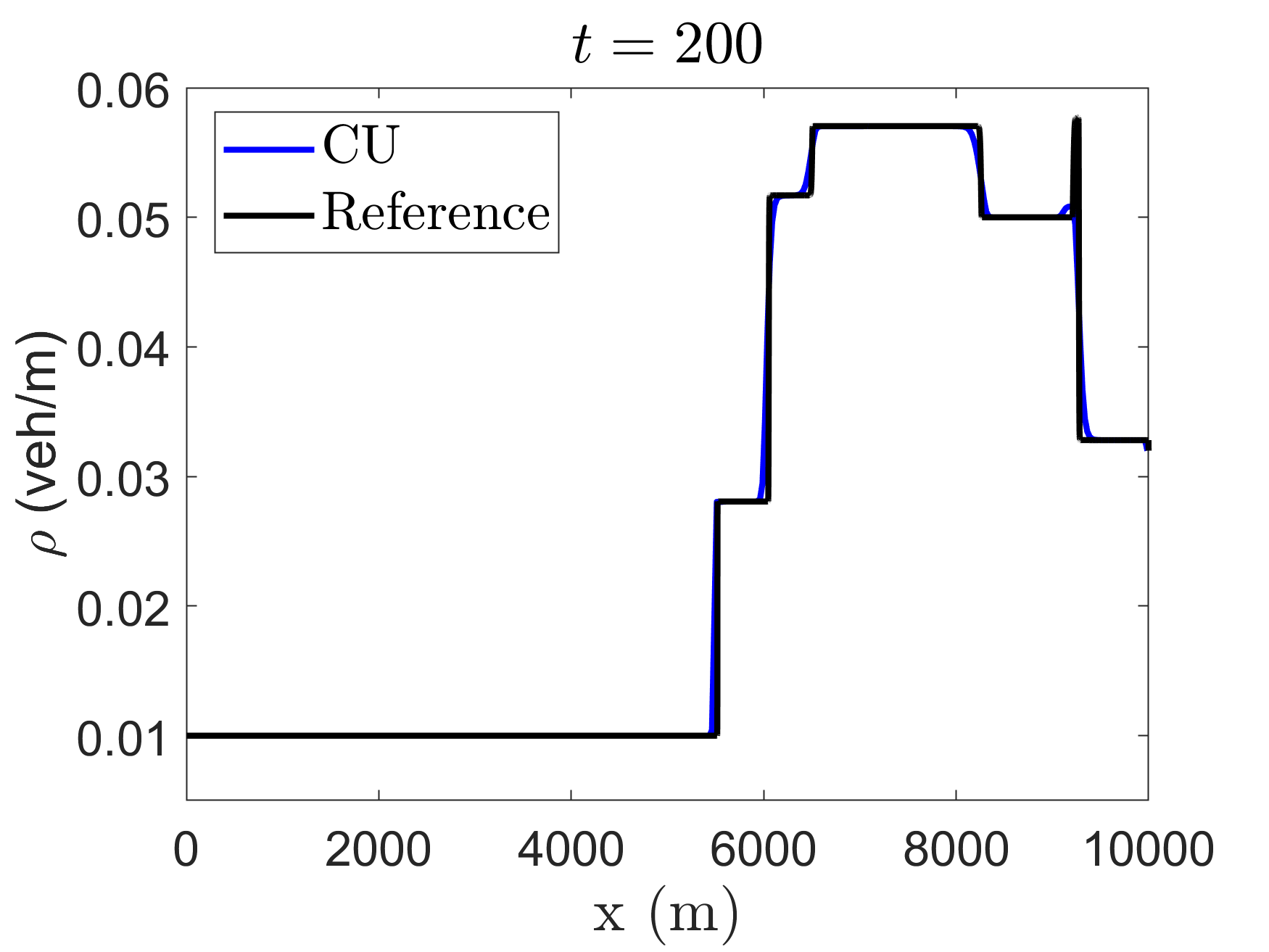}\hspace*{1cm}
            \includegraphics[trim=0.1cm 0.1cm 0.6cm 0.2cm, clip, width=5.4cm]{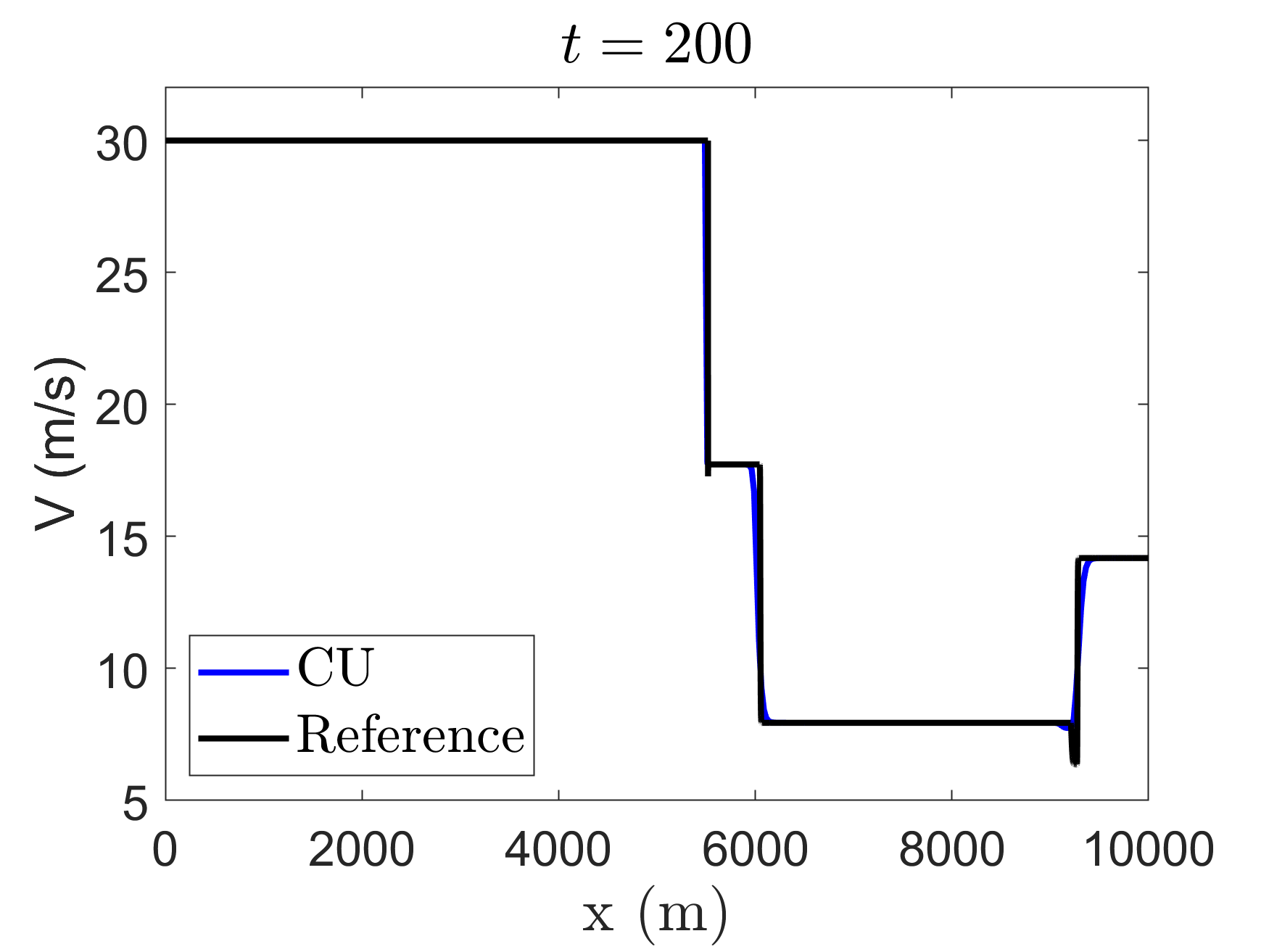}}
\vskip3pt
\centerline{\includegraphics[trim=0.1cm 0.1cm 0.6cm 0.2cm, clip, width=5.4cm]{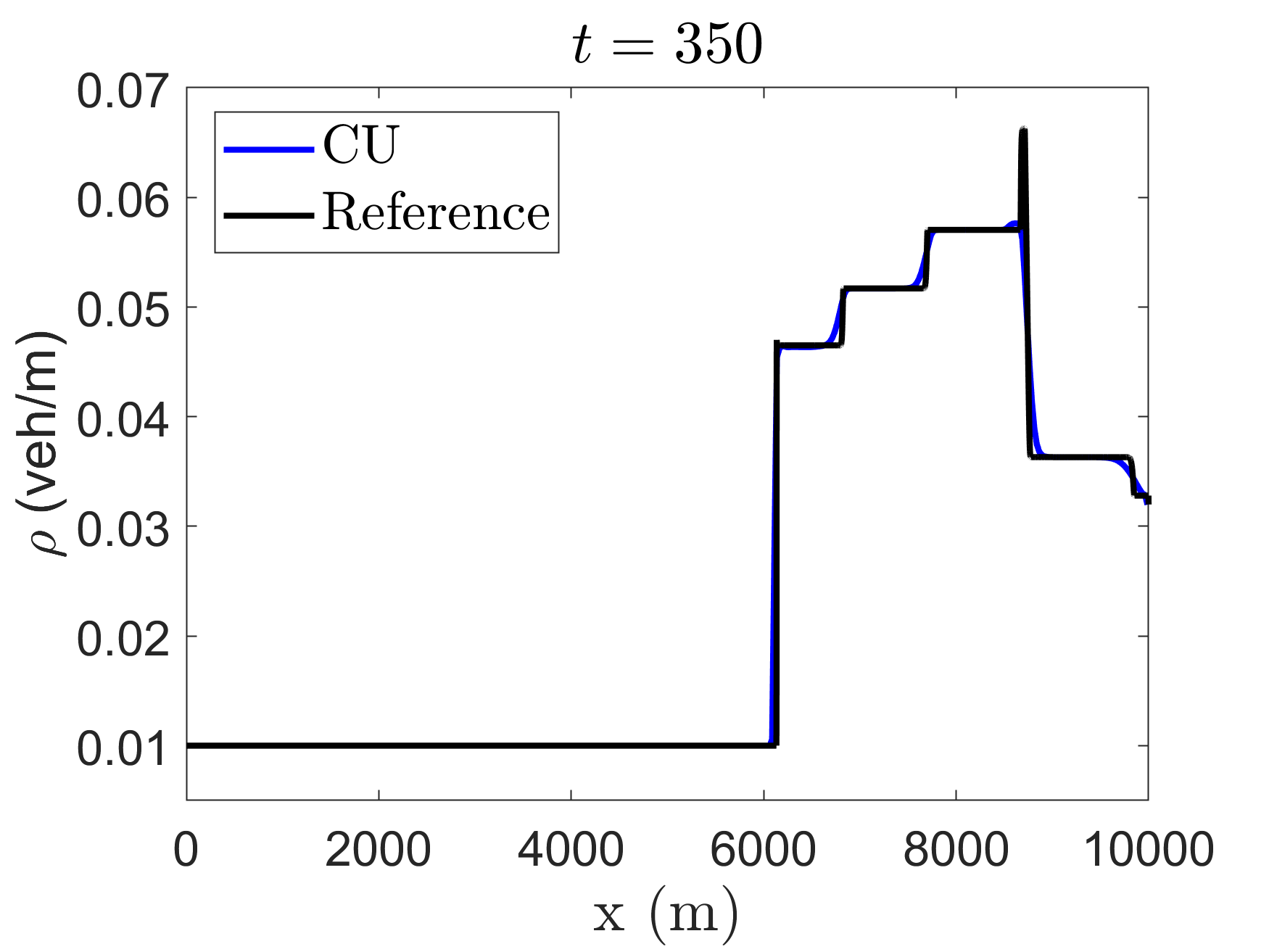}\hspace*{1cm}
            \includegraphics[trim=0.1cm 0.1cm 0.6cm 0.2cm, clip, width=5.4cm]{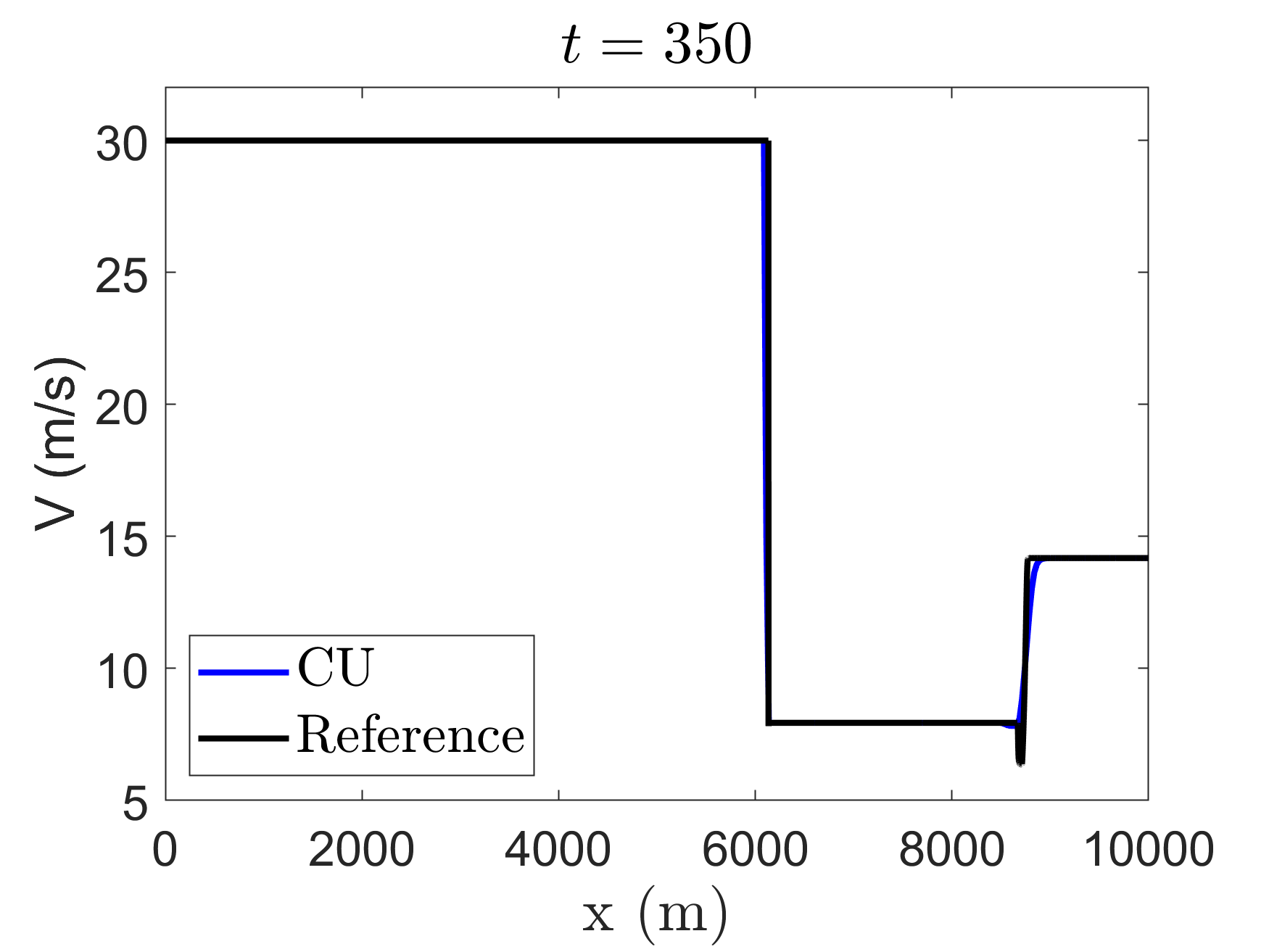}}
\vskip3pt
\centerline{\includegraphics[trim=0.1cm 0.1cm 0.6cm 0.2cm, clip, width=5.4cm]{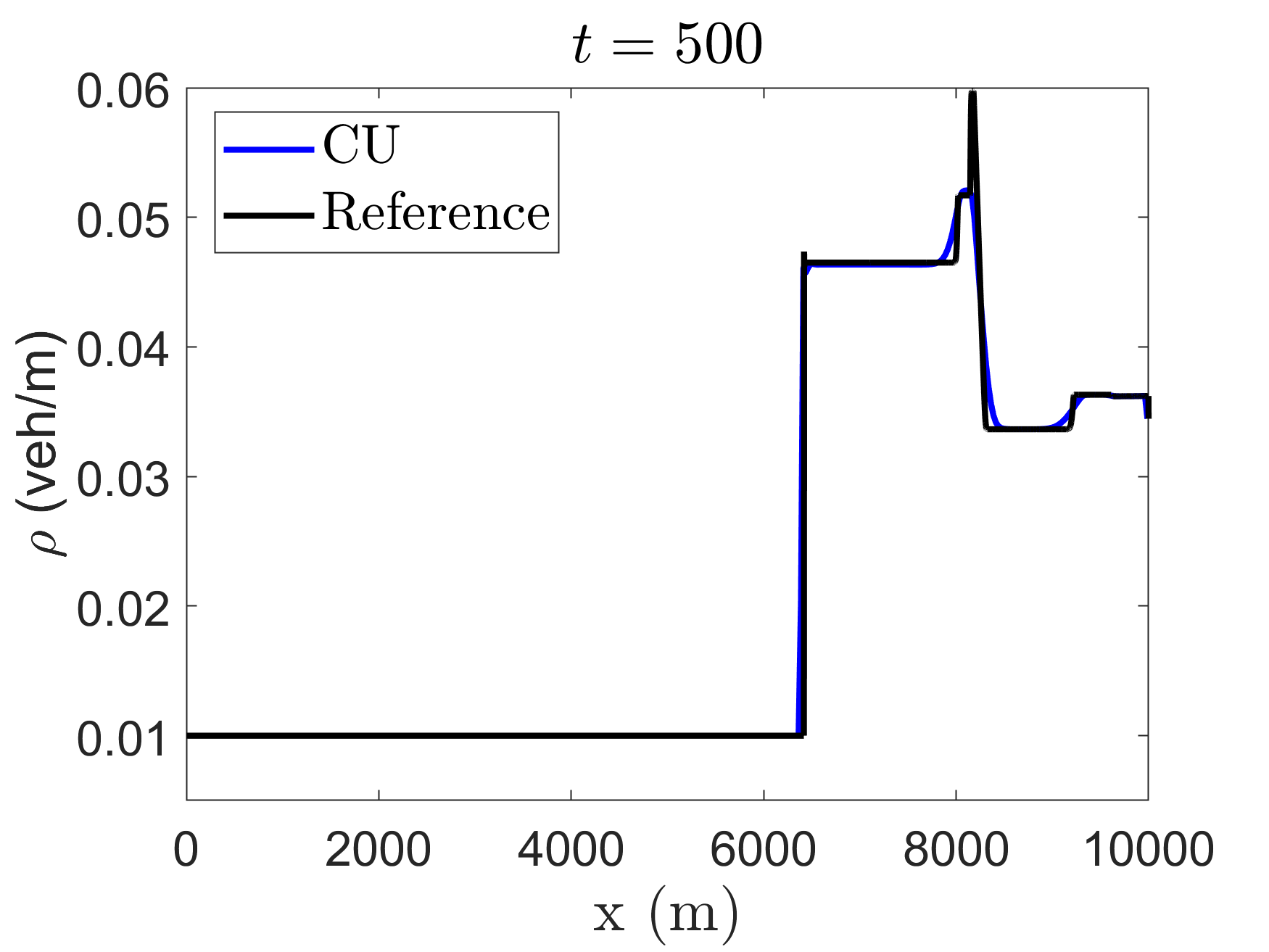}\hspace*{1cm}
            \includegraphics[trim=0.1cm 0.1cm 0.6cm 0.2cm, clip, width=5.4cm]{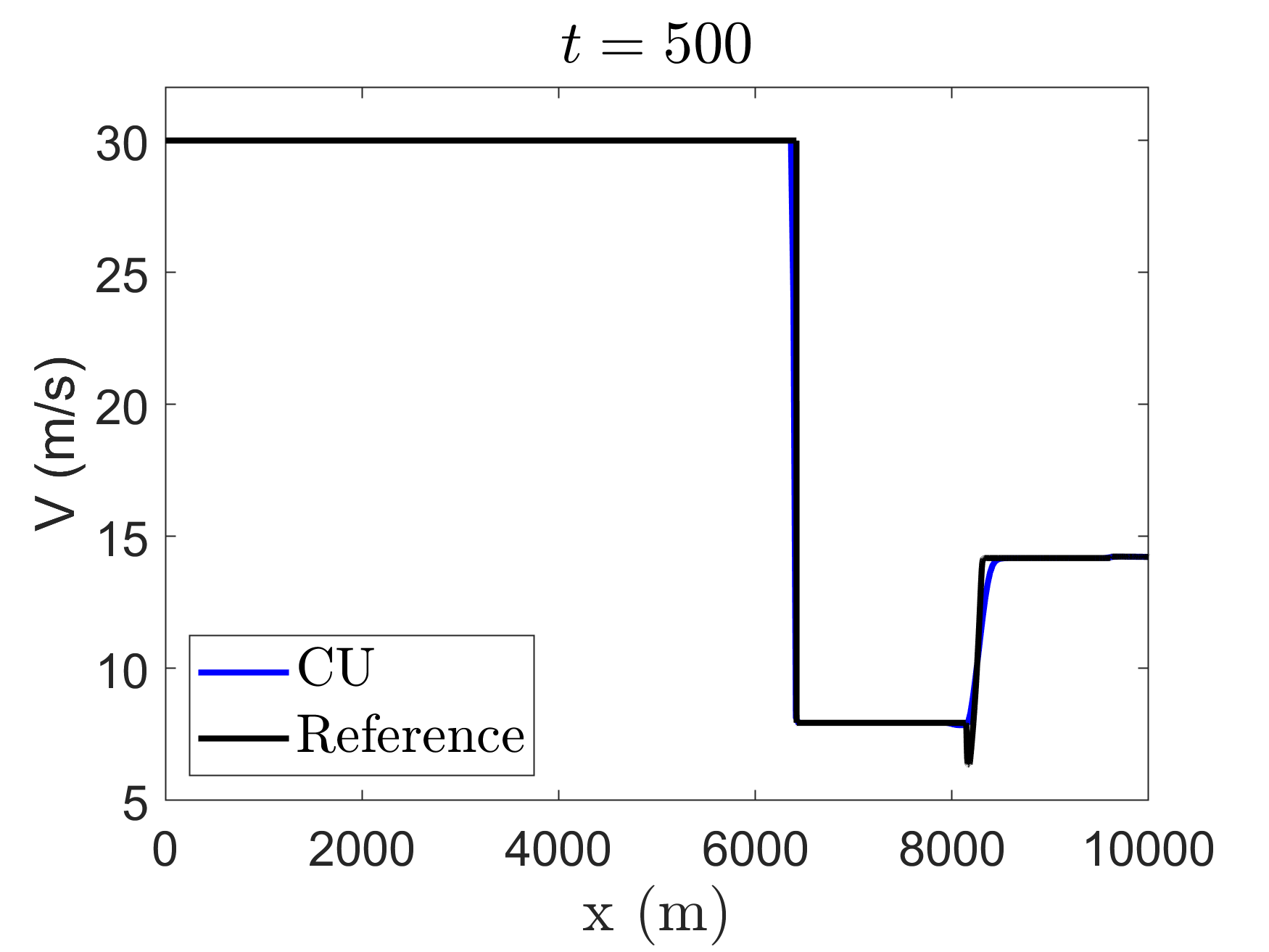}}
\caption{\sf Example 3: $\rho$ (left column) and $V$ (right column) at times $t=0$ (top row), $t=50$ (second row), $t=200$
(third row), $t=350$ (fourth row), and  $t=500$ (bottom row).\label{fig411}}
\end{figure}

\paragraph*{Example 3.} In this  example, we take  the computational domain $[0,L]$ with $L=10000$ and consider the following initial
conditions:
\begin{equation}
\rho(x,0)=\begin{cases}
0.01&\mbox{if}~x\le\frac{L}{3},\\
0.03&\mbox{if}~\frac{L}{3}<x\le\frac{2L}{3},\\
0.05&\mbox{otherwise},
\end{cases}\quad
V(x,0)=\begin{cases}
30&\mbox{if}~x\le\frac{L}{3},\\
17.729&\mbox{if}~\frac{L}{3}<x\le\frac{2L}{3},\\
7.941&\mbox{otherwise},\\
\end{cases}
\end{equation}
presented in Figure \ref{fig411} (top row). The corresponding values of the quantity $q-q^*$ are then
\begin{equation}
q(x,0)-q^*=\begin{cases}
-0.2800&\mbox{if}~x\le\frac{L}{3},\\
0.0546&\mbox{if}~\frac{L}{3}<x\le\frac{2L}{3},\\
-0.0225&\mbox{otherwise}.
\end{cases}
\end{equation}

In this example, the free boundary conditions are imposed at the left end of the computational domain, whereas at the right end of the
computational domain, we impose the following Dirichlet boundary conditions:
\allowdisplaybreaks
\begin{equation}
\resizebox{0.92\linewidth}{!}{$
\begin{aligned}
&\rho(L,t)=\begin{cases}
0.05+0.3\bigg[\cosh^{-2}\bigg(\dfrac{t-\frac{T_0}{2}}{W}\bigg)-\dfrac{1}{4}\cosh^{-2}\bigg(\dfrac{t-T_1-\frac{T_0}{2}}{W}\bigg)\bigg],
&~t\le\dfrac{2}{3}\,T_{\rm final},\\[2.5ex]
0.03+0.2\bigg[\cosh^{-2}\bigg(\dfrac{t-{T_0}}{W}\bigg)-\dfrac{1}{4}\cosh^{-2}\bigg(\dfrac{t-T_1-{T_0}}{W}\bigg)\bigg], 
&~\mbox{otherwise},
\end{cases}\\
&q(L,t)=\begin{cases}
\dfrac{21}{4}\rho_{\max}\Bigg(1+(a-1)\dfrac{\rho(L,t)}{\rho_{\max}}-\bigg[\bigg(a\dfrac{\rho(L,t)}{\rho_{\max}}\bigg)^{\!\!20}+
\bigg(1-\dfrac{\rho(L,t)}{\rho_{\max}}\bigg)^{\!\!20}\,\bigg]^{\frac{1}{20}}\Bigg),&~\rho(L,t)>\rho^f_{\rm cr}\\[2.5ex]
\rho(L,t)V_{\max}\bigg(1-\dfrac{\rho(L,t)}{\rho_{\max}}\bigg)^{\!\!-1},&~\mbox{otherwise}.
\end{cases}
\end{aligned}
$}
\end{equation}
Here, $W=201.25$, $T_0=1500$, $T_1=3000$, $T_{\rm final}=250$, and $a=\frac{30}{7}$.

While qualitatively similar to previous examples, this case is more complex due to a time-varying downstream boundary condition, leading to
the formation of additional wave types and regions. Another key difference is at the interface between the two congestion regions (initially
located at $x>\frac{L}{3}$), where upstream vehicles have smaller inter-vehicular spacing than at equilibrium ($q-q^*>0$), and downstream
vehicles have larger spacing ($q-q^*<0$). This results in multiple intermediate states and interacting waves, which are challenging to
capture in a non-oscillatory manner; see \cite{CKMZ23}.

We compute the numerical results using the proposed CU scheme until $t=T_{\rm final}$ on a uniform mesh with $\dx=25$. The computed values
of $\rho$ and $V$ at $t=50$, 200, 350, and 500 are shown in Figure \ref{fig411} together with a reference solution computed on a finer mesh
with $\dx=\frac{5}{4}$. As one can see, the scheme accurately captures solution structures from interacting waves, while staying sharp near
shock waves and contact discontinuities and remaining non-oscillatory around intermediate states.

\begin{figure}[ht!]
\centering
\begin{minipage}{0.5\textwidth}
\includegraphics[trim=0.1cm 0.1cm 0cm 0.2cm, clip, width=\linewidth]{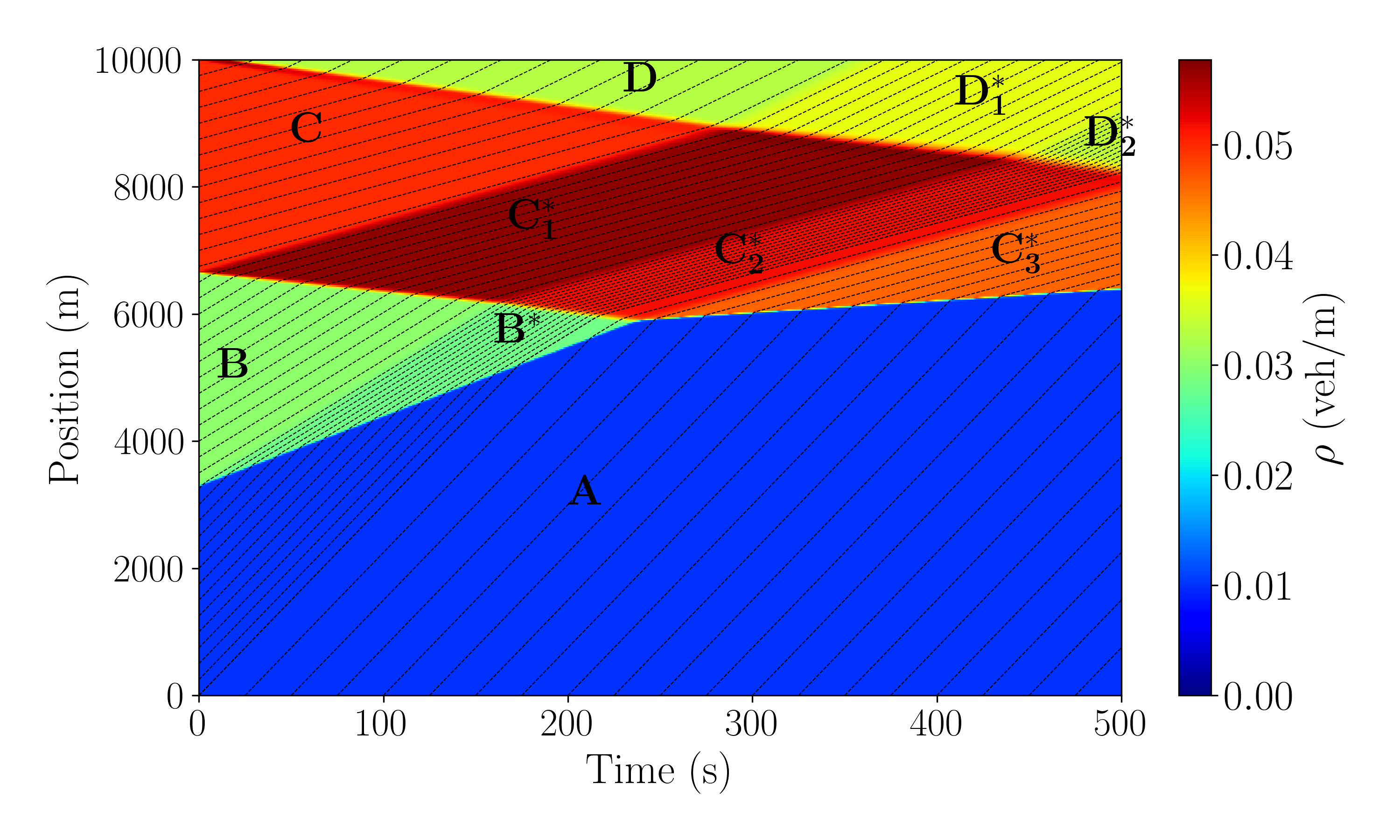}
\end{minipage}%
\begin{minipage}{0.5\textwidth}
\includegraphics[trim=0.1cm 0.1cm 0cm 0.2cm, clip, width=\linewidth]{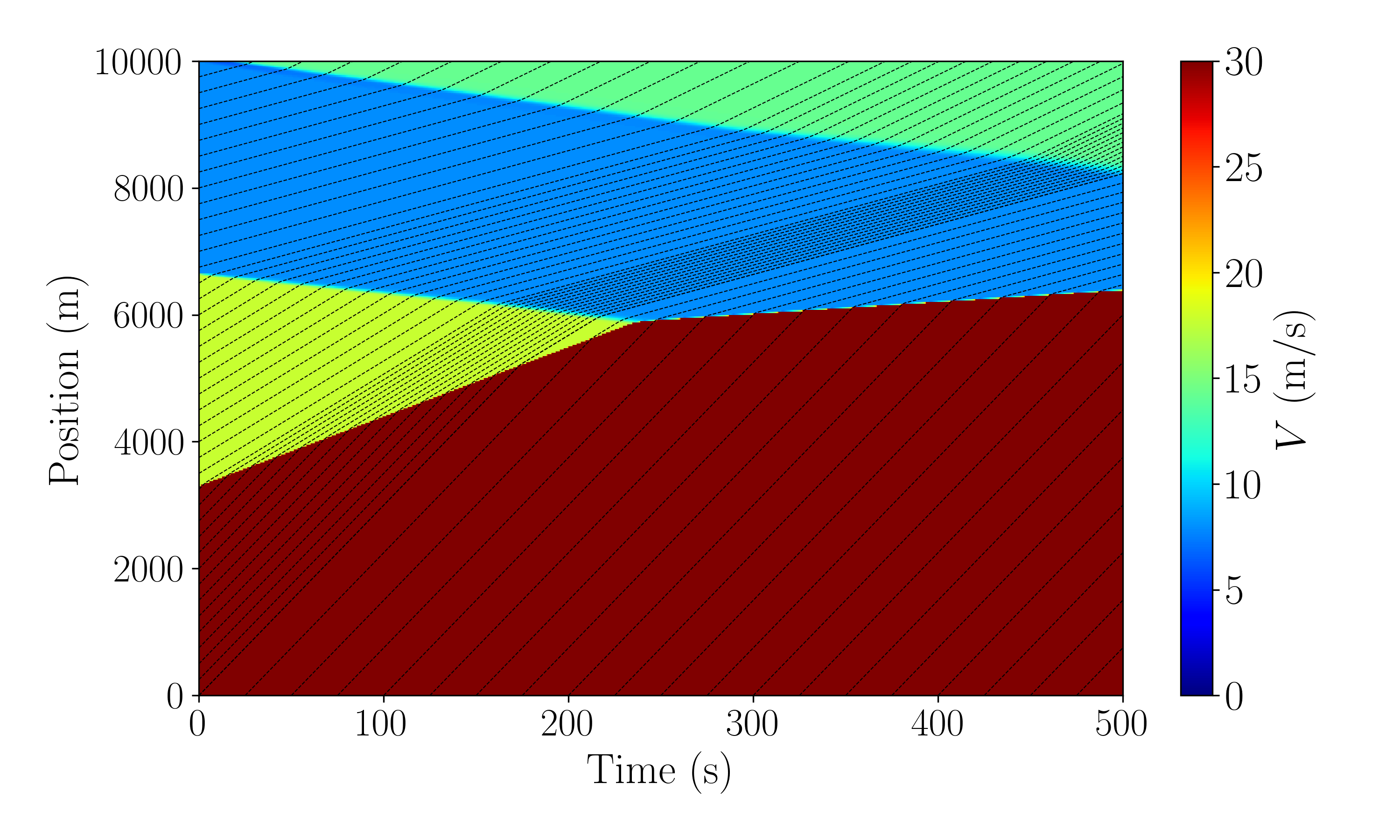}
\end{minipage}
\caption{\sf Example 3: Spatio-temporal evolution of $\rho$ (left) and $V$ (right)  for $t\in [0,500]$.\label{fig412}}
\end{figure}

To interpret the numerical results physically, we examine the spatio-temporal evolution of traffic density and speed in
Figure \ref{fig412}, where time-space diagrams are shown with superimposed vehicle trajectories similar to those in Figure
\ref{Example3_ST}. The comparison between speed and density graphs reveals that the solution domain can be divided into four regions where
speed is constant; see Figure \ref{fig412} (right). The density in these regions is, however, not constant except for the free-flow Region
A. The other regions with constant speed can be further divided into subregions with constant density, labeled as B, ${\rm B}^*_1$, C,
${\rm C}^*_1$, ${\rm C}^*_2$, ${\rm C}^*_3$, D, ${\rm D}^*_1$, and ${\rm D}^*_2$.

Let us begin a discussion on the physical interactions with the wave types arising at the initial discontinuities between the three main
Regions A, B, and C in the initial condition. The Riemann solution at the left end of Region B involves the formation and expansion of
Region ${\rm B}^*$ between Regions A and B. Region ${\rm B}^*$ is characterized by a shockwave moving upstream (backward) and a contact
discontinuity traveling forward (downstream) at a faster rate than the shockwave. Meanwhile, at $x=2L/3$ and $t=0$, one observes the
formation of a shockwave traveling backward, causing Region B to shrink, as well as a contact discontinuity traveling forward, which leads
to the formation of a high-density intermediate state (Region ${\rm C}^*_1$) between Regions B and C.

At approximately $x=6200\,\mbox{m}$ and $t=160\,\mbox{s}$, the contact discontinuity at the left end of Region B catches up with the
shockwave traveling backward at the right end of Region B, leading to the formation of Region ${\rm C}^*_2$ and the clearance of Region B.
Region ${\rm C}^*_2$ spans backward until around $x=5900\,\mbox{m}$ and $t=220\,\mbox{s}$, when Region ${\rm B}^*$ is cleared as the
shockwave at the left end of Region ${\rm B}^*$ meets the contact discontinuity at the right end of Region ${\rm B}^*$. This interaction
leads to the formation of Region ${\rm C}^*_3$, which arises as an intermediate state between another shockwave traveling upstream and a
contact discontinuity traveling downstream.

Finally, at $x=L$ and $t=0$, a shockwave arises traveling backward as the congested state in Region C discharges from the downstream, which
leads to the formation of another intermediate state between Region C and the downstream boundary (Region D). 

At approximately $x=9000\,\mbox{m}$ and $t=280\,\mbox{s}$, the shockwave traveling forward at the upper end of Region D meets the contact
discontinuity traveling forward at the lower end of Region C. This interaction forms region \({\rm D}^*\), characterized by a compound wave:
a contact discontinuity moving forward at the left end of region ${\rm D}^*_1$ and a shockwave moving upstream at a slightly slower rate at
the lower end of region ${\rm D}^*_1$.

Region ${\rm D}^*_1$ persists until approximately $x=8500\,\mbox{m}$ and $t=440\,\mbox{s}$, when the shockwave at the lower end of region
${\rm D}^*_1$ meets the contact discontinuity at the right end of region ${\rm C}^*_1$, forming region ${\rm D}^*_2$. This leads to another
compound wave with a shockwave traveling upstream at the lower end of region ${\rm D}^*_2$ and a contact discontinuity in region
${\rm D}^*_2$.

\begin{figure}[ht!]
\centerline{\includegraphics[trim=0.1cm 0.1cm 0.6cm 0.2cm, clip, width=5.4cm]{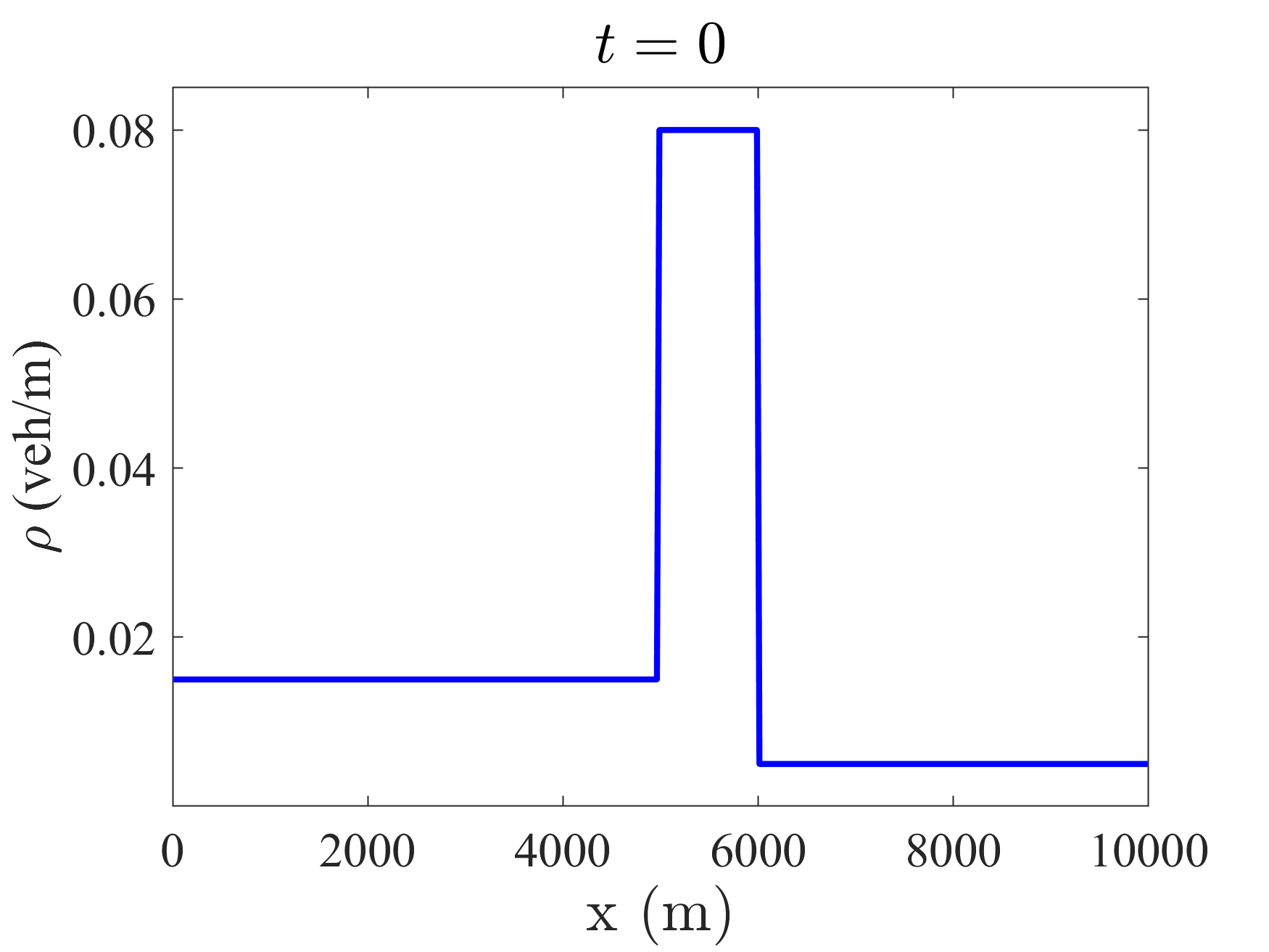}\hspace*{1cm}
            \includegraphics[trim=0.1cm 0.1cm 0.6cm 0.2cm, clip, width=5.4cm]{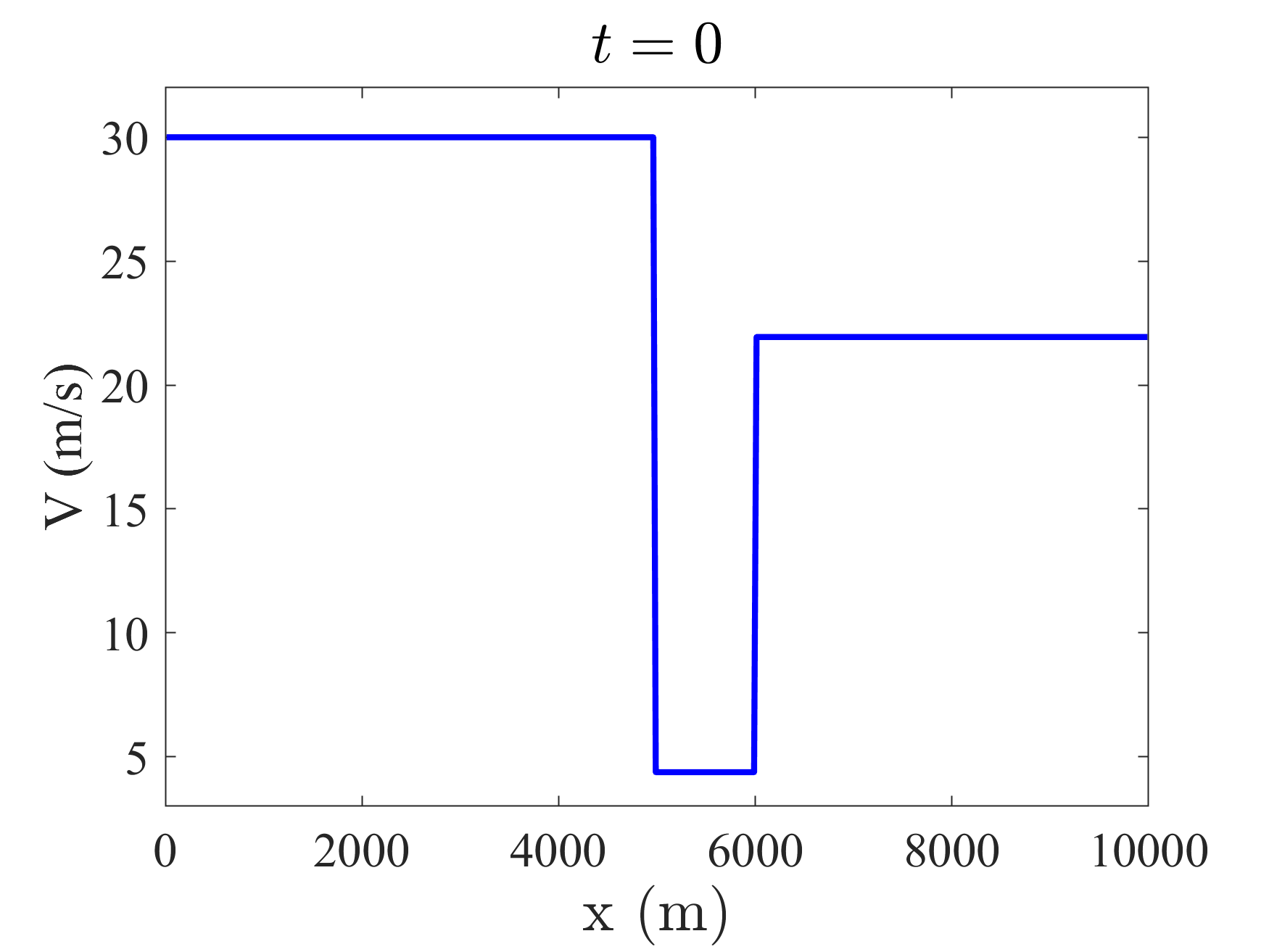}}
\caption{\sf Example 4: Initial conditions for $\rho$ (left) and $V$ (right).\label{fig413}}
\end{figure}

\paragraph*{Example 4.} In our final example, we aim to reconstruct a tangible real-world scenario using initial and boundary conditions
that are widely observable in real-world traffic. We take the computational domain $[0,L]$ with $L=10000$ and assume the following initial
conditions:
\begin{equation}
\rho(x,0)=\begin{cases}
0.015&\mbox{if}~x\le\frac{L}{2},\\
0.08&\mbox{if}~\frac{L}{2}<x\le\frac{3L}{5},\\
0.025&\mbox{otherwise},
\end{cases}\quad
V(x,0)=\begin{cases}
30&\mbox{if}~x\le\frac{L}{2},\\
4.375&\mbox{if}~\frac{L}{2}<x\le\frac{3L}{5},\\
21.94&\mbox{otherwise},
\end{cases}
\end{equation}
presented in Figure \ref{fig413}. The corresponding values of the quantity $q-q^*$ are then
\begin{equation}
q(x,0)-q^*=\begin{cases}
-0.1034&\mbox{if}~x\le\frac{L}{2},\\
0.1&\mbox{if}~\frac{L}{2}<x\le\frac{3L}{5},\\
0.0501&\mbox{otherwise}.
\end{cases}
\end{equation}

\begin{figure}[ht!]
\centerline{\includegraphics[trim=0.1cm 0.1cm 0.6cm 0.2cm, clip, width=5.4cm]{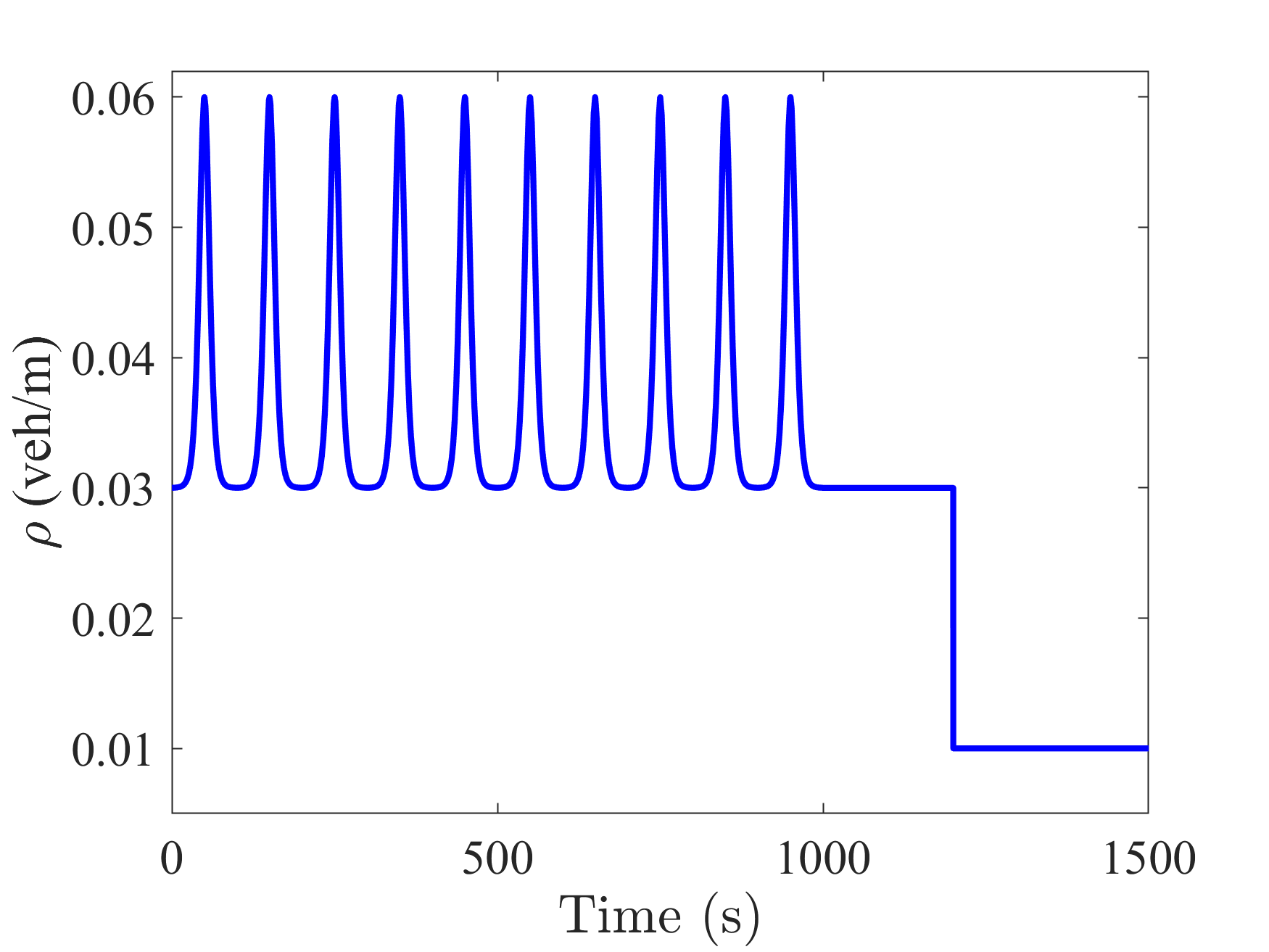}\hspace*{1cm}
            \includegraphics[trim=0.1cm 0.1cm 0.6cm 0.2cm, clip, width=5.4cm]{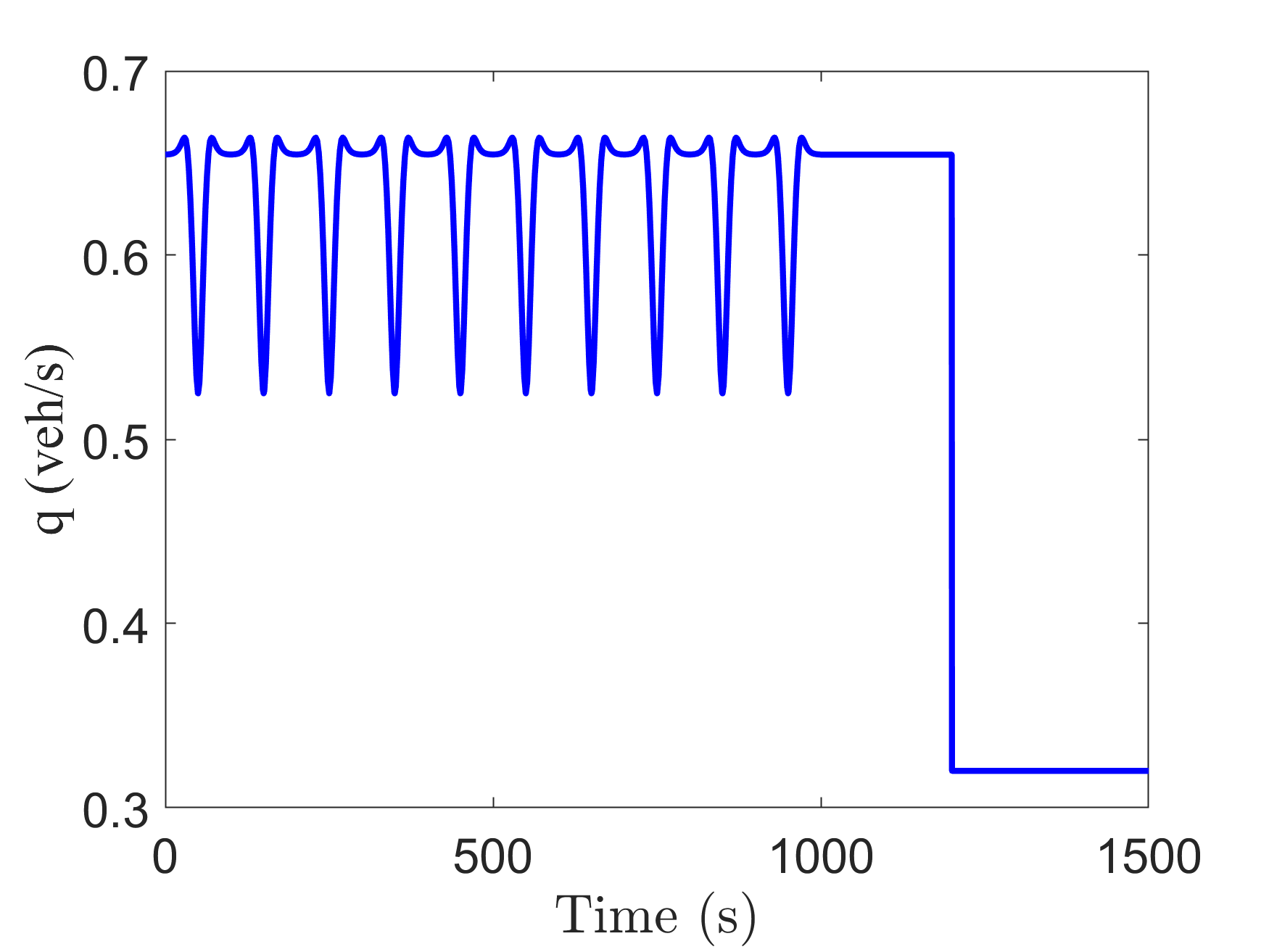}}
\caption{\sf Example 4: Boundary conditions for $\rho$ and $q$ at the right end of the computational domain.\label{fig413a}}
\end{figure}

The initial conditions in this scenario depict a strong bottleneck with highly congested traffic, represented as an initial queue in the
middle section of a freeway. At the left end of the computational domain, we impose a free upstream boundary condition, while at the right
end, we assume the downstream boundary conditions are time-dependent and periodic to reflect a stop-and-go traffic condition further
downstream. Note that in traffic flow, drivers mainly respond to conditions ahead, and the majority of traffic waves propagate backward.
Therefore, if the situation further downstream is congested, the outflow cannot exit. This example is practically relevant, as such
stop-and-go patterns are widely observable in real-world traffic. The downstream boundary condition is mathematically expressed as
\begin{equation}
\resizebox{0.91\linewidth}{!}{$
\begin{aligned}	
&\rho(L,t)=\begin{cases}
0.03+0.03\bigg[\cosh^{-2}\bigg(\dfrac{t-{T_0}}{W}\bigg)-2.05\cosh^{-2}\bigg(\dfrac{t-T_1-{T_0}}{W}\bigg)\bigg],&~t\le1000,\\
0.03,&~1000<t<1200,\\
0.01,&~\mbox{otherwise},
\end{cases}\\[1ex]
&q(L,t)=\begin{cases}
\dfrac{21}{4}\rho_{\max}\Bigg(1+(a-1)\dfrac{\rho(L,t)}{\rho_{\max}}-\bigg[\bigg(a\dfrac{\rho(L,t)}{\rho_{\max}}\bigg)^{\!\!20}+
\bigg(1-\dfrac{\rho(L,t)}{\rho_{\max}}\bigg)^{\!\!20}\,\bigg]^{\frac{1}{20}}\Bigg),&~t<1200,\\[2.5ex]
\rho(L,t)V_{\max}\bigg(1-\dfrac{\rho(L,t)}{\rho_{\max}}\bigg)^{\!\!-1},&~\mbox{otherwise},
\end{cases}
\end{aligned}
$}
\end{equation}
which are illustrated in Figure \ref{fig413a}. Here,  $W=10.25$, $T_0=\lfloor\frac{t}{100}\rfloor+50$, $T_1=3000$, and $a=\frac{30}{7}$.

This scenario is common in real-world traffic. On a homogeneous road without on- or off-ramps, we assume the initial queue in the middle can
represent traffic conditions after a crash that occurred at $x=3L/5$ before the simulation, with the crash scene persisting for a while.
This results in very slow queue discharge and light traffic conditions downstream in the area $x>3L/5$ until the crash scene is cleared at
$t=0$. Additionally, the downstream boundary condition involves stop-and-go traffic up until $t=1000\,\mbox{s}$, then transitions to a
mildly congested state between $t=1000\,\mbox{s}$ and $t=1200\,\mbox{s}$, followed by free-flow traffic again for $t>1200\,\mbox{s}$.
Stop-and-go traffic at the downstream boundary is assumed to result from highly periodic incoming flow through a hypothetical on-ramp
further downstream in the computational domain. Our aim is to investigate the spatio-temporal evolution of traffic flow in a complex
scenario involving multiple interacting bottlenecks. This test case provides a good basis for examining these interactions while
demonstrating the performance of the proposed CU scheme and highlighting the physical implications of the results.

\begin{figure}[ht!]
\centering
\begin{minipage}{0.5\textwidth}
\includegraphics[trim=0.1cm 0.1cm 0cm 0.2cm, clip, width=\linewidth]{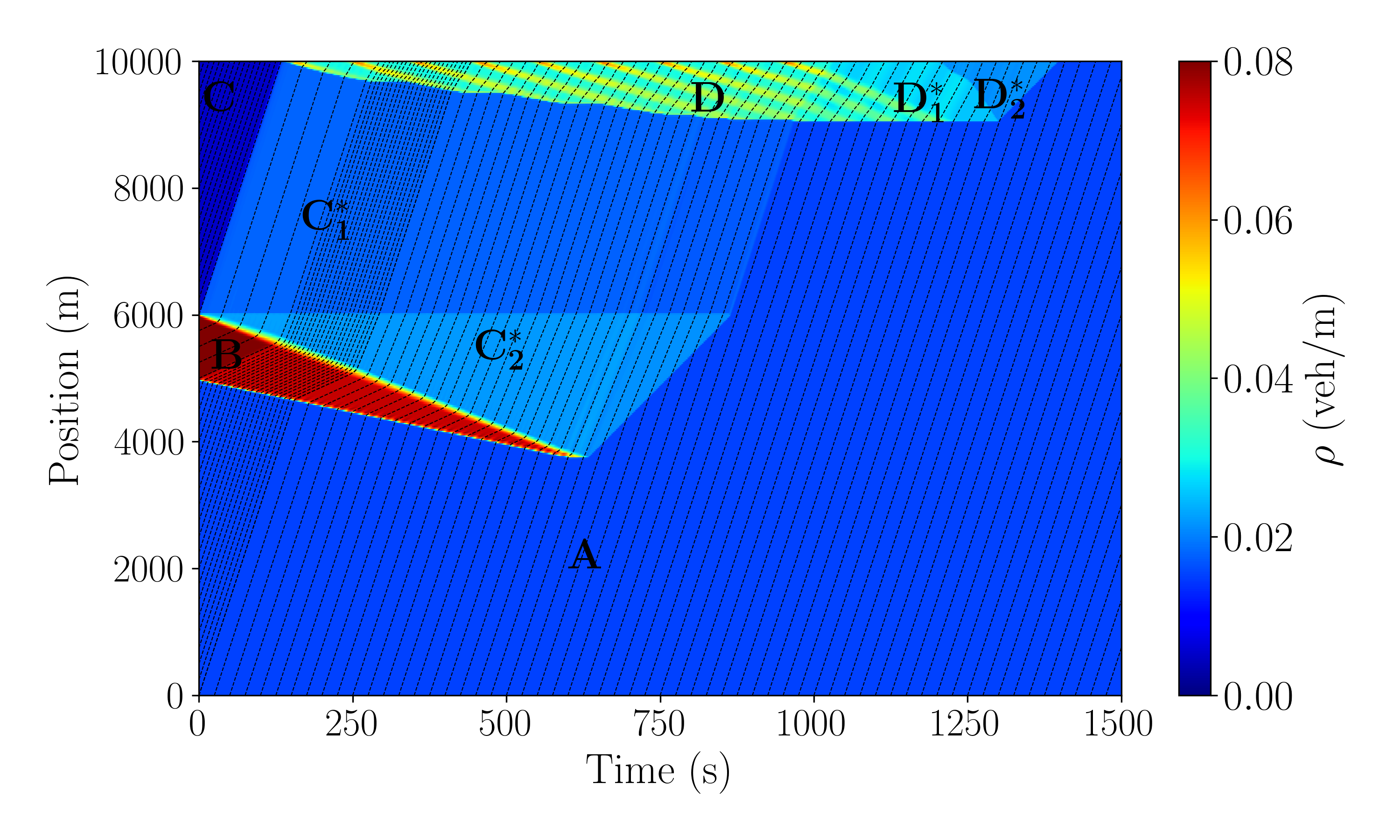}
\end{minipage}%
\begin{minipage}{0.5\textwidth}
\includegraphics[trim=0.1cm 0.1cm 0cm 0.2cm, clip, width=\linewidth]{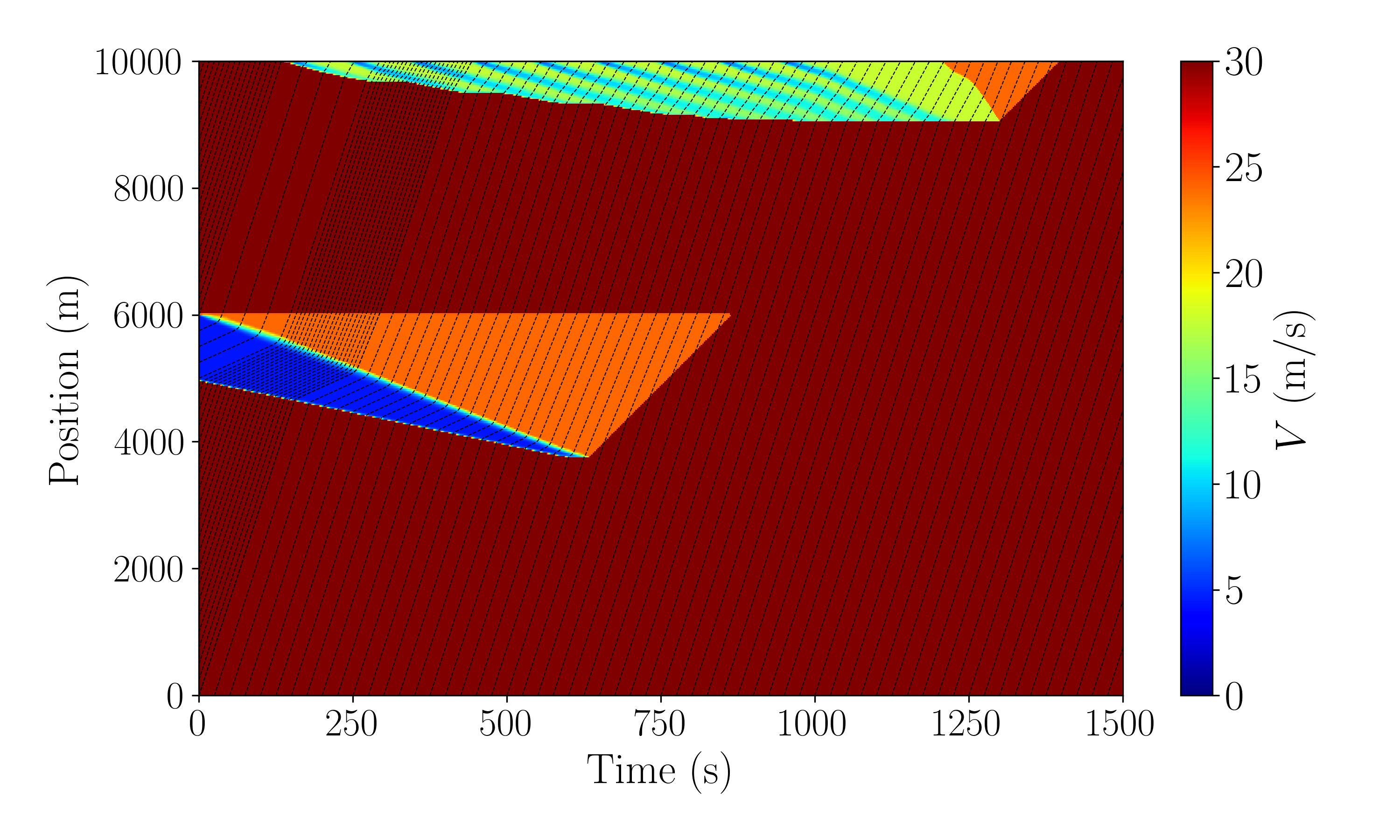}
\end{minipage}
\caption{\sf Example 4: Spatio-temporal evolution of $\rho$ (left) and $V$ (right) for $t\in [0,1500]$.\label{fig416}}
\end{figure}

We compute the numerical results until the final time $T_{\rm final}=1500$ on a uniform mesh with $\dx=25$. To stay focused, we present only
the spatio-temporal evolution of density and speed in Figure \ref{fig416}, which again shows that the proposed scheme can capture the
overall structures of the numerical solution in a non-oscillatory manner.

To interpret the numerical results physically, we examine the spatio-temporal evolution of traffic density and speed. In Figure
\ref{fig416}, time-space diagrams are shown with superimposed vehicle trajectories similar to those in Figures \ref{Example3_ST} and
\ref{fig412}. We consider several regions (A, B, ${\rm C}^*_1$, ${\rm C}^*_2$, D, ${\rm D}^*_1$, and ${\rm D}^*_2$) in the density panel,
which mark key phases of traffic evolution.

First, Region A corresponds to free-flow traffic entering from the upstream boundary. Second, we focus on the initial queue, highlighted as
Region B, which involves a strong bottleneck. Region B is characterized by two backward-traveling shocks at its upper and lower ends. The
shockwave at the lower end arises as free-flowing traffic encounters highly congested traffic, whereas the shockwave at the upper end arises
as vehicles leave the congestion region. Note that the upper-end shock travels backward faster as vehicles leaving Region B adjust their
speed to the maximum in the congested region $(V_{c+}=24\,\mbox{m/s}$), and, as a result, at around $x=3800\,\mbox{m}$ and
$t=600\,\mbox{s}$, the shockwaves meet, leading to the dissolution of Region B.

The upper end of Region B neighbors three regions with varying characteristics, as compound waves arise and interact at the interface
between the initial queue and free-flow traffic. First, Region C is characterized by a contact discontinuity at its right end. In this
region, vehicles maintain their free-flow initial condition. Due to the low traffic density, vehicles in Region C do not create
backward-propagating waves when reaching $x=L$, where slow-moving stop-and-go traffic occurs.

Next, Region ${\rm C}^*_2$ remains congested, with vehicles adapting their speed to the maximum in the congestion domain $\Omega_c$
($V_{c+}=24\,\mbox{m/s}$). Region ${\rm C}^*_2$ is characterized by intermediate congestion, confined by three shocks at its upper, lower,
and right ends. The shock at the lower end arises as vehicles leaving the highly congested bottleneck in Region B adjust their speed to the
maximum in the congestion domain. The shock at the right end forms as vehicles from the free-flow Region A enter and pass through
${\rm C}^*_2$. The shock at the upper end occurs during the transition from intermediate congestion to the free-flow domain, and as the
traffic states neighboring in the either side of the shock, are in a relatively comparable range.

Finally, in Region ${\rm C}^*_1$, free-flow traffic with a density close to $\rho^f_{\rm cr}$ persists. Region ${\rm C}^*_1$ is confined by
a shock from the lower end, contact discontinuities at the left and right ends, and a backward-traveling shock wave with varying speed at
the upper end, which is the interface between Regions ${\rm C}^*_1$ and D. This shock wave arises due to the high traffic flow in Region
${\rm C}^*_1$, and upon meeting the congested condition at the downstream boundary (see Figure \ref{fig413a}), it leads to the formation and
backward propagation of the oscillatory Region D. We note that the interface between Regions ${\rm C}^*_1$ and D has a varying shock speed
in the time-space diagram as the tiny cluster-like waves inside Region D also disperse and propagates forward. Region D continues to
propagate backward until shortly after Region ${\rm C}^*_1$ is cleared (roughly at $x=9300\,\mbox{m}$ and around
$t=950\,\mbox{s}$), at which point the last vehicles leaving Region ${\rm C}^*_2$ reach Region D. From this point onward, Region A with low
density meets Region D from behind, causing the interface between the two regions to become a stationary shock. This stationary shock
continues to remain as the downstream boundary condition (at $x=L$) is replaced with a steady congested state during
$1000\,\mbox{s}<t<1200\,\mbox{s}$. As the downstream boundary condition (at $x=L$) becomes free again for $t\ge1200\,\mbox{s}$, a new Region
${\rm D}^*_2$ forms as vehicles in Region ${\rm D}^*_1$ leave the congested state. A shock wave arises at the interface between Regions
${\rm D}^*_1$ and ${\rm D}^*_2$, traveling backward until it meets the free-flow Region A at approximately $x=8300\,\mbox{m}$ and
$t=1650\,\mbox{s}$. From this point, a forward-traveling shock wave arises, causing the dissolution of Region ${\rm D}^*_2$ as vehicles
adapt their dynamics to those in Region A.

\section{Conclusions}
Phase-transition hyperbolic traffic flow models meet most requirements for non-equilibrium traffic flow models, such as maintaining the
anisotropy property, maintaining consistency with the inverse correlation between driver speed and intervehicular spacing, satisfying zero
speed at maximum density, and distinctly representing dynamics in free-flow and congested phases. Phase-transition models incorporate
``time-gap'', a key element of automated vehicle control logic, as a conserved variable in their equations of state, giving them significant
potential for future applications in the era of connected and automated vehicles.

Despite their significant potential, phase-transition models have been largely overlooked in traffic flow theory, particularly in numerical
simulations and analysis. This gap is due to the mathematical complexity of these models, such as their discontinuous solution domains and
fluxes, which make them challenging to simulate numerically.

This paper presents the development of a second-order semi-discrete finite-volume central-upwind scheme for the phase-transition traffic
flow model. The developed scheme has been applied to several challenging numerical examples, demonstrating its ability to capture numerical
solution structures sharply and without oscillations. Using this method, the phase-transition model's performance was studied in various
real-world traffic scenarios, in which multiple traffic phases interact, leading to interacting bottlenecks and complex solution structures.
Such investigations were previously infeasible, but the proposed numerical method now facilitates the implementation of the model in complex
real-world traffic scenarios.

There are several future directions for this work. First, future studies can utilize the proposed numerical scheme in conjunction with
real-world traffic data to apply optimization-based parameter estimation of the phase-transition model and investigate the performance of
the calibrated model against real-world traffic phenomena. Second, the proposed numerical scheme solves traffic flow on homogeneous roads
without junctions or interchanges. Future work can extend this freeway-level solver by incorporating boundary coupling conditions at
junctions to implement the phase-transition model on complex networks with multiple interacting bottlenecks. Such approaches are crucial for
model-predictive optimization of traffic networks. These topics are currently being investigated in ongoing work by the authors.

\subsection*{Acknowledgment}
The work of S. Chu was supported in part by the DFG (German Research Foundation) through HE5386/19-3, 27-1. The work of A. Kurganov
was supported in part by NSFC grant 12171226 and W2431004.

\appendix
\section{Generalized Minmod Reconstruction}\label{appa}
In this appendix, we briefly describe a piecewise linear generalized minmod reconstruction \cite{lie03,Nessyahu90,Sweby84}.

Assume that the cell averages $\,\xbar\psi_j$ of a certain function $\psi(x)$ are given. We use them to reconstruct a second-order piecewise
linear interpolant
\begin{equation}
\widetilde\psi(x)=\,\xbar\psi_j+(\psi_x)_j(x-x_j),\quad x\in C_j,
\label{equ3.4}
\end{equation}
and then to compute the right- and left-sided point values of $\psi$ at the cell interfaces $x=x_\jph$:
\begin{equation}
\psi^-_\jph=\,\xbar\psi_j+\frac{\dx}{2}(\psi_x)_j,\quad\psi^+_\jph=\,\xbar\psi_{j+1}-\frac{\dx}{2}(\psi_x)_{j+1}.
\end{equation}
In order to ensure a non-oscillatory nature of this reconstruction, we compute the slopes $(\psi_x)_j$ in \eref{equ3.4} using a generalized
minmod limiter:
\begin{equation}
(\psi_x)_j={\rm minmod}\left(\theta\,\frac{\,\xbar\psi_j-\,\xbar\psi_{j-1}}{\dx},\,\frac{\,\xbar\psi_{j+1}-\,\xbar\psi_{j-1}}{2\dx},\,
\theta\,\frac{\,\xbar\psi_{j+1}-\,\xbar\psi_j}{\dx}\right),\quad\theta\in[1,2],
\label{equ3.5}
\end{equation}
where the minmod function is defined as
\begin{equation}
{\rm minmod}(z_1,z_2,\ldots):=\begin{cases}
\min_j\{z_j\}&\mbox{if}~z_j>0\quad\forall\,j,\\
\max_j\{z_j\}&\mbox{if}~z_j<0\quad\forall\,j,\\
0            &\mbox{otherwise}.
\end{cases}
\end{equation}
The parameter $\theta$ in \eref{equ3.5} can be used to control the oscillations: larger $\theta$'s correspond to sharper but, in general,
more oscillatory reconstructions.

\section{Local Characteristic Decomposition}\label{appb}
In order to suppress the oscillations in the piecewise linear reconstruction of $\rho$ and $q$ inside the congested domain, we apply the
generalized minmod limiter from Appendix \ref{appa} to the local characteristic variables. To this end, we first introduce the matrix
\begin{equation}
\begin{aligned}
{\widehat A}_\jph=A(\widehat\mU_\jph)=\begin{pmatrix}
-\dfrac{\widehat q_\jph}{\rho_{\max}}&\dfrac{\rho_{\max}-\widehat\rho_\jph}{\rho_{\max}}\\[2.5ex]
\dfrac{\widehat q_\jph(q^*-\widehat q_\jph)}{\widehat\rho_\jph^{\,2}}&
\dfrac{(q^*-2\widehat q_\jph)(\widehat\rho_\jph-\rho_{\max})}{\widehat\rho_\jph\rho_{\max}}
\end{pmatrix},
\end{aligned}
\end{equation}
where $A=\frac{\partial\mF}{\partial\mU}$, $\widehat\rho_\jph=\frac{\,\xbar\rho_j+\,\xbar\rho_{j+1}}{2}$ and
$\widehat q_\jph=\frac{\,\xbar q_j+\,\xbar q_{j+1}}{2}$, and construct the following matrix, which consists of the two eigenvectors of
${\widehat A}_\jph$:
\begin{equation}
\begin{aligned}
R_\jph=\begin{pmatrix}
\dfrac{\widehat\rho_\jph}{\widehat q_\jph-q^*}&\dfrac{\widehat\rho_\jph(\rho_{\max}-\widehat\rho_\jph)}{\widehat q_\jph\rho_{\max}}\\[2.5ex]
1&1\end{pmatrix}.
\end{aligned}
\end{equation}
 We then compute its inverse
\begin{equation}
\begin{aligned}
R^{-1}_\jph=\begin{pmatrix}
\dfrac{\widehat q_\jph(\widehat q_\jph-q^*)\rho_{\max}}{\widehat\rho_\jph(\widehat q_\jph\widehat\rho_\jph+
q^*(\rho_{\max}-\widehat\rho_\jph))}&\dfrac{(\widehat q_\jph-q^*)(\widehat\rho_\jph-\rho_{\max})}
{\widehat q_\jph\widehat\rho_\jph+q^*(\rho_{\max}-\widehat\rho_\jph)}\\[3.5ex]
\dfrac{\widehat q_\jph(q^*-\widehat q_\jph)\rho_{\max}}{\widehat\rho_\jph(\widehat q_\jph\widehat\rho_\jph+
q^*(\rho_{\max}-\widehat\rho_\jph))}&\dfrac{\widehat q_\jph\rho_{\max}}{\widehat q_\jph\widehat\rho_\jph+q^*(\rho_{\max}-\widehat\rho_\jph)}
\end{pmatrix},
\end{aligned}
\end{equation}
and introduce the local characteristic variables $\bm\Gamma$ in the neighborhood of $x=x_\jph$:
\begin{equation}
\bm\Gamma_k=R^{-1}_\jph\,\xbar\mU_k,\quad k=j-1,j,j+1,j+2.
\end{equation}
Finally, we compute the slopes $(\bm\Gamma_x)_j$ using the generalized minmod limiter \eref{equ3.5} applied to $\bm\Gamma$ in a
componentwise manner, evaluate the point values of $\bm \Gamma$:
\begin{equation}
\bm\Gamma^-_\jph=\bm\Gamma_j+\frac{\dx}{2}(\bm\Gamma_x)_j\quad\mbox{and}\quad
\bm\Gamma^+_\jph=\bm\Gamma_{j+1}-\frac{\dx}{2}(\bm\Gamma_x)_{j+1},
\end{equation}
and end up with obtaining the corresponding point values of $\mU$:
\begin{equation}
\mU^\pm_\jph=R_\jph\bm\Gamma^\pm_\jph.
\end{equation}

\bibliographystyle{siamnodash}
\bibliography{ref}

\begin{thebibliography}{10}

\bibitem{blandin2013phase}
{\sc S.~Blandin, J.~Argote, A.~M. Bayen, and D.~B. Work}, {\em Phase transition
  model of non-stationary traffic flow: Definition, properties and solution
  method}, Transport. Res. Part B-Meth, 52 (2013), pp.~31--55.

\bibitem{blandin2011general}
{\sc S.~Blandin, D.~B. Work, P.~Goatin, B.~Piccoli, and A.~M. Bayen}, {\em A
  general phase transition model for vehicular traffic}, SIAM J. Appl. Math.,
  71 (2011), pp.~107--127.

\bibitem{chalons2008godunov}
{\sc C.~Chalons and P.~Goatin}, {\em Godunov scheme and sampling technique for
  computing phase transitions in traffic flow modeling}, Interface. Free.
  Bound., 10 (2008), pp.~197--221.

\bibitem{CCHKL_22}
{\sc A.~Chertock, S.~Chu, M.~Herty, A.~Kurganov, and
  M.~Luk\'a\v{c}ov\'a-Medvi\v{d}ov\'a}, {\em Local characteristic decomposition
  based central-upwind scheme}, J. Comput. Phys., 473 (2023).
\newblock Paper No. 111718.

\bibitem{CKMZ23}
{\sc S.~Chu, A.~Kurganov, S.~Mohammadian, and Z.~Zheng}, {\em Fifth-order
  A-WENO path-conservative central-upwind scheme for behavioral non-equilibrium
  traffic models}, Commun. Comput. Phys., 33 (2023), pp.~692--732.

\bibitem{CKX_24}
{\sc S.~Chu, A.~Kurganov, and R.~Xin}, {\em New low-dissipation central-upwind
  schemes. {P}art {II}}, J. Sci. Comput., 103 (2025).
\newblock Paper No. 33.

\bibitem{Colombo02}
{\sc R.~M. Colombo}, {\em A {$2\times 2$} hyperbolic traffic flow model}, Math.
  Comput. Modelling, 35 (2002), pp.~141--163.

\bibitem{Colombo02a}
{\sc R.~M. Colombo}, {\em Hyperbolic phase transitions in traffic flow}, SIAM
  J. Appl. Math., 63 (2002), pp.~708--721.

\bibitem{CG14}
{\sc R.~M. Colombo and M.~Garavello}, {\em Phase transition model for traffic
  at a junction}, J. Math. Sci. (N.Y.), 196 (2014), pp.~30--36.

\bibitem{daganzo1995requiem}
{\sc C.~Daganzo}, {\em Requiem for second-order fluid approximations of traffic
  flow}, Transport. Res. B-Meth., 29 (1995), pp.~277--286.

\bibitem{don9}
{\sc W.~S. Don, D.-M. Li, Z.~Gao, and B.-S. Wang}, {\em A characteristic-wise
  alternative {WENO}-{Z} finite difference scheme for solving the compressible
  multicomponent non-reactive flows in the overestimated quasi-conservative
  form}, J. Sci. Comput., 82 (2020).
\newblock Paper No. 27, 24 pp.

\bibitem{Gottlieb11}
{\sc S.~Gottlieb, D.~Ketcheson, and C.-W. Shu}, {\em Strong stability
  preserving {R}unge-{K}utta and multistep time discretizations}, World
  Scientific Publishing Co. Pte. Ltd., Hackensack, NJ, 2011.

\bibitem{Gottlieb12}
{\sc S.~Gottlieb, C.-W. Shu, and E.~Tadmor}, {\em Strong stability-preserving
  high-order time discretization methods}, SIAM Rev., 43 (2001), pp.~89--112.

\bibitem{kerner2016failure}
{\sc B.~S. Kerner}, {\em Failure of classical traffic flow theories: Stochastic
  highway capacity and automatic driving}, Physica. A., 450 (2016),
  pp.~700--747.

\bibitem{Kurganov01}
{\sc A.~Kuganov, S.~Noelle, and G.~Petrova}, {\em Semidiscrete central-upwind
  schemes for hyperbolic conservation laws and {H}amilton-{J}acobi equations},
  SIAM J. Sci. Comput., 23 (2001), pp.~707--740.

\bibitem{Kurganov07}
{\sc A.~Kurganov and C.-T. Lin}, {\em On the reduction of numerical dissipation
  in central-upwind schemes}, Commun. Comput. Phys., 2 (2007), pp.~141--163.

\bibitem{Kurganov00}
{\sc A.~Kurganov and E.~Tadmor}, {\em New high-resolution semi-discrete central
  schemes for {H}amilton-{J}acobi equations}, J. Comput. Phys., 160 (2000),
  pp.~720--742.

\bibitem{KX_22}
{\sc A.~Kurganov and R.~Xin}, {\em New low-dissipation central-upwind schemes},
  J. Sci. Comput., 96 (2023).
\newblock Paper No. 56.

\bibitem{lie03}
{\sc K.-A. Lie and S.~Noelle}, {\em On the artificial compression method for
  second-order nonoscillatory central difference schemes for systems of
  conservation laws}, SIAM J. Sci. Comput., 24 (2003), pp.~1157--1174.

\bibitem{mohammadian2021benchmarking}
{\sc S.~Mohammadian, Z.~Zheng, Md.~M. Haque, and A.~Bhaskar}, {\em Performance
  of continuum models for realworld traffic flows: Comprehensive benchmarking},
  Transport. Res. B-Meth., 147 (2021), pp.~132--167.

\bibitem{mohammadian2023continuum}
{\sc S.~Mohammadian, Z.~Zheng, Md.~M. Haque, and A.~Bhaskar}, {\em Continuum
  modeling of freeway traffic flows: State-of-the-art, challenges and future
  directions in the era of connected and automated vehicles}, Commun. Transp.
  Res., 3 (2023).
\newblock Paper No. 100107, 25 pp.

\bibitem{Nessyahu90}
{\sc H.~Nessyahu and E.~Tadmor}, {\em Nonoscillatory central differencing for
  hyperbolic conservation laws}, J. Comput. Phys., 87 (1990), pp.~408--463.

\bibitem{Nonomura20}
{\sc T.~Nonomura and K.~Fujii}, {\em Characteristic finite-difference {WENO}
  scheme for multicomponent compressible fluid analysis: overestimated
  quasi-conservative formulation maintaining equilibriums of velocity,
  pressure, and temperature}, J. Comput. Phys., 340 (2017), pp.~358--388.

\bibitem{QZWZW2017}
{\sc Y.~Qian, J.~Zeng, N.~Wang, J.~Zhang, and B.~Wang}, {\em A traffic flow
  model considering influence of car-following and its echo characteristics},
  Nonlinear Dyn., 89 (2017), pp.~1099--1109.

\bibitem{Qiu02}
{\sc J.~Qiu and C.-W. Shu}, {\em On the construction, comparison, and local
  characteristic decomposition for high-order central {WENO} schemes}, J.
  Comput. Phys., 183 (2002), pp.~187--209.

\bibitem{Shu20}
{\sc C.-W. Shu}, {\em Essentially non-oscillatory and weighted essentially
  non-oscillatory schemes}, Acta Numer., 5 (2020), pp.~701--762.

\bibitem{Sweby84}
{\sc P.~K. Sweby}, {\em High resolution schemes using flux limiters for
  hyperbolic conservation laws}, SIAM J. Numer. Anal., 21 (1984),
  pp.~995--1011.

\bibitem{treiber2010three}
{\sc M.~Treiber, A.~Kesting, and D.~Helbing}, {\em Three-phase traffic theory
  and two-phase models with a fundamental diagram in the light of empirical
  stylized facts}, Transport. Res. B-Meth., 44 (2010), pp.~983--1000.

\bibitem{ZQLZX2023}
{\sc J.~Zeng, Y.~Qian, J.~Li, Y.~Zhang, and D.~Xu}, {\em Congestion and energy
  consumption of heterogeneous traffic flow mixed with intelligent connected
  vehicles and platoons}, Physica A., 609 (2023).
\newblock Paper No. 128331.

\bibitem{ZQYZX2022}
{\sc J.~Zeng, Y.~Qian, F.~Yin, L.~Zhu, and D.~Xu}, {\em A multi-value cellular
  automata model for multi-lane traffic flow under lagrange coordinate},
  Comput. Math. Organ. Th., 28 (2022), pp.~178--192.

\bibitem{ZHANG2023128556}
{\sc J.~Zhang, Y.~Qian, J.~Zeng, X.~Wei, and H.~Li}, {\em Hybrid
  characteristics of heterogeneous traffic flow mixed with electric vehicles
  considering the amplitude of acceleration and deceleration}, Physica A., 614
  (2023).
\newblock Paper No. 128556.

\end{thebibliography}
\end{document}